\documentclass[a4paper,12pt]{amsart}
\usepackage{tabularx}
\usepackage[dvips]{graphicx}
\usepackage{amsmath}
\usepackage{amssymb}
\usepackage{amsbsy}
\usepackage{amsgen}
\usepackage{amsopn}
\usepackage{amscd}
\usepackage{amstext}
\usepackage{amsthm}

\theoremstyle{definition}
\newtheorem{definition}{Definition}[section]

\theoremstyle{plain}
\newtheorem{theorem}{Theorem}[section]
\newtheorem{lemma}{Lemma}[section]
\newtheorem{corollary}{Corollary}[section]
\newtheorem{proposition}{Proposition}[section]

\theoremstyle{remark}
\newtheorem{remark}{Remark}[section]

\theoremstyle{definition}
\newtheorem{example}{Example}[section]

\newtheorem{conv}{Conventions}

\begin{document}

\title{Khovanov homology and words}
\author{Tomonori Fukunaga and Noboru Ito}
\maketitle

\begin{abstract}
This paper is concerned with nanowords, a generalization of links, introduced by Turaev.  It is shown that the system of bigraded homology groups is an invariant of nanowords by introducing a new notion $\mathcal{U}_{L}$.  This paper gives two examples which show the independence of this invariant from some of Turaev's homotopy invariants.  
\end{abstract}

\section{Introduction}
In this paper, an alphabet is a finite set and a letter means its element. 
A word on an alphabet
$\mathcal{A}$ is a finite sequence of letters in $\mathcal{A}$, 
and a phrase is a finite sequence of words on $\mathcal{A}$.  If each letter in the alphabet appears exactly twice in 
the word (respectively phrase), then we call this word a Gauss word
(respectively a Gauss phrase).\par 
In the papers \cite{turaev1} and \cite{turaev2}, 
V. Turaev introduced the theory of topology of words and phrases
(see also \cite{turaev3}).  The theory is a combinatorial extension of the 
theory of virtual knots and links. Let $\alpha$ be an alphabet endowed 
with an involution $\tau$. An $\alpha$-alphabet is a pair  
an alphabet $\mathcal{A}$ and a map $|\cdot|$ from $\mathcal{A}$ to
$\alpha$. We call this map $|\cdot|$ a projection. Then 
a nanoword (respectively nanophrase) over $\alpha$ is 
a pair an $\alpha$-alphabet $\mathcal{A}$ and 
a Gauss word (respectively a Gauss phrase) on $\mathcal{A}$.\par
Turaev defined an equivalence relation which is called $S$-homotopy 
on nanophrases for a subset $S$ of $\alpha^{3}$. Two nanophrases 
are $S$-homotopic each other if and only if they are related by a 
finite sequence of isomorphisms and $S$-homotopy moves (i), (ii), (iii),
and inverse moves of $S$-homotopy moves
(definitions of an isomorphism of 
nanophrases and $S$-homotopy moves is given in Section \ref{turaevswords}).\par
In the paper \cite{turaev2}, Turaev gave some geometric meanings of 
the theory of topology of words and phrases. Let $\alpha_{\ast}$ be a set 
consisting of four letters $a_{+}$, $a_{-}$, $b_{+}$, and $b_{-}$.
Moreover let $\tau_{\ast}$ be an involution on $\alpha_{\ast}$
and $S_{\ast}$ be a subset of $\alpha_{\ast}^{3}$ which is defined 
in Section 2. 
Then Turaev proved that the set of the $S_{\ast}$-homotopy classes of 
nanophrases over $\alpha_{\ast}$ is one to one corresponds to 
the set of the stable equivalence classes of 
ordered pointed $k$-component curves.\par
In the same paper, Turaev considered nanophrases over $\alpha_{1}$ with 
an involution $\tau_{1}$ up to $S_{1}$-homotopy (we call 
elements of this set pseudolinks, see Section 2 of \cite{turaev2} for more 
details). 
Turaev constructed the Jones polynomial for pseudolinks.\par   
The purpose of this paper is to construct Khovanov homology
for pseudolinks and we prove $S_{1}$-homotopy invariance of the homology.  As an application of the Khovanov homology for pseudolinks to
nanophrases over an arbitrary alphabet, we construct a new invariant for nanophrases by using the Khovanov homology
for pseudolinks. Moreover, we show that the new invariant is independent of the homotopy invariants $\lambda$ 
and characteristic sequences for nanowords
which were defined in \cite{turaev1}.\par
\begin{figure}
\centerline{\includegraphics[width=6cm]{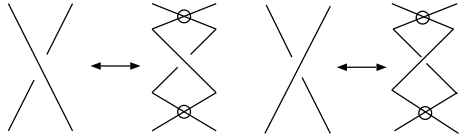}}
\caption{Virtualization.}\label{fig1}
\end{figure}

In \cite{manturov2}, Manturov defines Khovanov homology with $\mathbb{Z}/2\mathbb{Z}$ coefficients of virtual knots and observes that the homology he defined is invariant under virtualizations (see Figure \ref{fig1}).  
Since we can view pseudolinks as virtual links modulo virtualizations, Manturov's Khovanov homology is an invariant of pseudolinks.  We can see that our homology $KH^{i, j}$ $=$ $H^{i}(C^{*, j}, d)$ is the same as that of Manturov by considering the isomorphism between Viro's definition \cite{viro} and Bar-Natan's definition \cite{bar-natan} of Khovanov homology.  However, our construction has some benefits over Manturov's construction as follows.  First, it is easy to calculate our homology as an invariant of long virtual strings because there is a natural bijection from pseudolinks without shifts to long virtual strings using word theory.  Second, our construction makes it easier to calculate our new invariants of nanophrases over an arbitrary alphabet.  Third, our proof of invariance of the homology is simple because the proof is given by explicit chain homotopy maps and retractions on complexes.\par  
This paper is organized as follows.
In Section 2, we review the theory of topology of words and 
phrases. In Section 3, we introduce Kauffman-type states of 
pseudolinks, and we define the Jones polynomial for pseudolinks
by using Kauffman-type states of pseudolinks.  After that we show that this definition is equivalent to the definition of the Jones polynomial for pseudolinks introduced by Turaev in 
\cite{turaev2}.  In Section 4, we define Khovanov homology with $\mathbb{Z}/2\mathbb{Z}$ coefficient, and in Section 5, we prove $S_{1}$-homotopy invariance of the homology which is constructed 
in Section 4. In Section 6, we discuss an application of 
the Khovanov homology for pseudolinks to the theory of 
topology of nanowords and nanophrases over an arbitrary alphabet.    
\section{Turaev's theory of words}\label{turaevswords}
\subsection{Nanowords and Nanophrase}
For our preliminary discussions, we define {\it nanophrases} and their $S$-{\it homotopy} in the same manner as that in Turaev's original paper \cite[Section 2]{turaev1}, \cite[Section 2]{turaev2}, Gibson's paper \cite[Section 2]{andrew3}, or Fukunaga's paper \cite[Section 2.1]{fukunaga1}; these papers provide a detailed description of the terminology of nanophrases.  

An {\it alphabet} is a finite set and {\it letters} are its elements.    
For an alphabet $\alpha$, an $\alpha$-alphabet $\mathcal{A}$ is a set where every element $A$ of $\mathcal{A}$ has a projection $|~|:$ $A\mapsto |A|\in\alpha$.  
A {\it word of length $n \ge 1$} in an alphabet $\mathcal{A}$ is a mapping $w:$ $\hat{n} \to {\mathcal A}$, where $\hat{n}$ $=$ $\{ i \in {\mathbb{N}}~|~1 \le i \le n \}$.  Such a word is encoded by the sequence $w(1)w(2) \cdots w(n)$.  By definition, there exists a unique word $\emptyset$ of length $0$.  We define an {\it opposite word} by writing the letters of a word $w$ in the opposite order.  For example, if $w$ $=$ $abc$, then $w^{-}$ $=$ $cba$.  A word $w:$ $\hat{n} \to {\mathcal A}$ is a {\it Gauss word} in an alphabet $\mathcal{A}$ if each element of $\mathcal A$ is the image of precisely two elements of $\hat{n}$ or $w$ is $\emptyset$.  A {\it Gauss phrase} in an alphabet $\mathcal{A}$ is a sequence of words $x_1$, $x_2, \dots$, $x_m$ in $\mathcal{A}$ denoted by $x_1|x_2|\dots|x_m$ such that $x_1 x_2 \cdots x_m$ is a Gauss word in $\mathcal{A}$.  We call $x_{i}$ the $i$th component of the Gauss phrase.  In particular, if a Gauss phrase has only one component, that component is a Gauss word.  A {\it nanoword} $(\mathcal A, w)$ over $\alpha$ is a pair (an $\alpha$-{\it alphabet} $\mathcal A$, a Gauss word in the alphabet $\mathcal A$).  For a nanoword $(\mathcal{A}, w=w_{1}w_{2}\cdots w_{k})$ over $\alpha$ consisting of subwords $w_{i}$ $(1 \le i \le k)$ of $w$, a {\it nanophrase} of length $k \ge 0$ over $\alpha$ is defined as $(\mathcal{A}, w_{1}|w_{2}|\cdots|w_{k})$.   Whenever possible, $(\mathcal{A}, w_{1}|w_{2}|\cdots|w_{k})$ is indicated by simple symbols: $w_{1}|w_{2}|\cdots|w_{k}$, $(\mathcal{A}, P)$, or $P$.  
We call $w_{i}$ the $i$th component of the nanophrase.  
An arbitrary nanoword $w$ over $\alpha$ yields a nanophrase $w$ of length $1$.  However, we distinguish between nanowords and nanophrases of length $1$.  By definition, there exists a unique nanophrase of length $0$.  Note the fact that $\emptyset$ is not a nanophrase of length $0$ (see \cite[Subsection 6.1]{turaev2}.  Turaev did not differentiate between nanowords and nanophrases of length $1$).  We denote the nanophrase of length $0$ by $\emptyset$.  Note that we distinguish the nanophrase $\emptyset|\emptyset | \dots |\emptyset$ of length $k$ from that $\emptyset|\emptyset | \dots |\emptyset$ of length $l$ if $k$ $\neq$ $l$.  

An {\it isomorphism} of $\alpha$-alphabets $\mathcal{A}_1$, $\mathcal{A}_2$ is a bijection $f:$ $\mathcal{A}_1$ $\to$ $\mathcal{A}_2$ such that $|A|=|f(A)|$ for an arbitrary $A \in \mathcal{A}_1$.  Two nanophrases $(\mathcal{A}_{1}, p_{1}= w_{1}|w_{2}|\cdots|w_{k})$ and $(\mathcal{A}_{2}, p_{2}= w'_{1}|w'_{2}|\cdots|w'_{k'})$ over $\alpha$ are {\it isomorphic} if $k = k'$ and there exists an isomorphism of $\alpha$-alphabets $f:$ $\mathcal{A}_1$ $\to$ $\mathcal{A}_2$ such that $w'_{i}=fw_{i}$ for every $i$ $\in \{1, 2, \cdots, k\}$.  
  
\subsection{Homotopy of nanophrases}
To define a homotopy of nanophrases, we consider a finite set $\alpha$ with an involution $\tau:$ $\alpha \to \alpha$ and a subset $S$ $\subset \alpha \times  \alpha \times \alpha$.  We call the triple $(\alpha, \tau, S)$ {\it homotopy data}.  
Turaev defined an $S$-{\it homotopy} as follows (see \cite[Section 2.2]{turaev2}, \cite[Section 2.1]{fukunaga1}, \cite[Section 2]{andrew3}).  
\begin{definition}
Let $(\alpha, \tau, S)$ be homotopy data.  
Two nanowords $(\mathcal{A}_1, w_{1})$ and $(\mathcal{A}_2, w_{2})$ are $S$-homotopic if one nanophrase is changed into the other by the finite sequence of the isomorphisms and the following three type deformations (1)--(3), called {\it homotopy moves}, and their inverses.  The relation $S$-homotopy is denoted by $\simeq_{S}$.  

(H1) Replace ($\mathcal{A}$, $xAAy$) with ($\mathcal{A} \setminus \{A\}$, $xy$) for $\mathcal{A}$, and $x$, $y$ are words in $\mathcal{A} \setminus \{A\}$ that possibly include the character $|$ such that $xy$ is a Gauss phrase.  

(H2) Replace $(\mathcal{A}, xAByBAz)$ with $(\mathcal{A} \setminus \{A, B\}, xyz)$ if $A, B \in \mathcal{A}$ with $\tau(|A|)=|B|$ where $x$, $y$, $z$ are words in $\mathcal{A} \setminus \{A, B\}$ that possibly include the character $|$ such that $(xyz)$ is a Gauss phrase.  

(H3) Replace $(\mathcal{A}, xAByACzBCt)$ with $(\mathcal{A}, xBAyCAzCBt)$ for $(|A|, |B|, |C|)$ $\in S$, where $x$, $y$, $z$, and $t$ are words in $\mathcal{A}$ that possibly include the character $|$ such that $(xyzt)$ is a Gauss phrase.  
\end{definition}

We note the following two lemmas from \cite[Lemma 2.1, Lemma 2.2]{turaev2} (see \cite[Lemma 2.4, Lemma 2.5]{fukunaga1}).  

\begin{lemma}
Let $(\alpha, \tau, S)$ be homotopy data and $\mathcal A$ be an $\alpha$-alphabet.  
Let $A$, $B$, and $C$ be distinct letters in $\mathcal A$ and let $x$, $y$, $z$, and $t$ be words that possibly include the character $|$ in the alphabet $\mathcal A$ $\setminus$ $\{A, B, C\}$ such that $(xyzt)$ is a Gauss phrase in this alphabet.  Then, 

(i) $(\mathcal{A}, xAByCAzBCt)$ $\simeq_{S}$ $(\mathcal{A}, xBAyACzCBt)$~{\text{for}}~$(|A|, \tau(|B|), |C|) \in S$; 

(ii) $(\mathcal{A}, xAByCAzCBt)$ $\simeq_{S}$ $(\mathcal{A}, xBAyACzBCt)$~{\text{for}}~$(\tau(|A|), \tau(|B|), |C|) \in S$; and

(iii) $(\mathcal{A}, xAByACzCBt)$ $\simeq_{S}$ $(\mathcal{A}, xBAyCAzBCt)$~{\text{for}}~$(\tau(|A|), |B|, |C|) \in S$.  
\end{lemma}

\begin{lemma}\label{abab}
Suppose that $S$ $\cap$ $(\alpha \times \{b\} \times \{b\})$ $\neq$ $\emptyset$ for all $b \in \alpha$.  Let $(\mathcal{A}, xAByABz)$ be a nanophrase over $\alpha$ with $|B| = \tau(|A|)$, where $x$, $y$, and $z$ are words that possibly include the character $|$ in the alphabet $\mathcal{A} \setminus \{A, B\}$ such that $xyz$ is a Gauss phrase in this alphabet.  Then, 
$(\mathcal{A}, xAByABz)$ $\simeq_{S}$ $(\mathcal{A} \setminus \{A, B\}, xyz)$.  
\end{lemma}

\begin{definition}
Let $\alpha$ be a finite set.  Consider an involution $\nu:$ $\alpha$ $\to$ $\alpha$ called the {\it shift involution}.  
The $\nu$-shift of a nanoword $(\mathcal{A}, w:\hat{n} \to \mathcal{A})$ over $\alpha$ is the nanoword $(\mathcal{A'}, w':\hat{n} \to \mathcal{A'})$ obtained by steps (1)--(3):
(1) Let $\mathcal{A}$ $:=$ $(\mathcal{A} - \{A\})$ $\cup$ $\{A_{\nu}\}$, where $A_{\nu}$ is a letter not belonging to $\mathcal{A}$.  

(2) The projection $\mathcal{A}' \to \alpha$ extends the given projection $\mathcal{A}-\{A\} \to \alpha$ by $|A_{\nu}|=\nu(|A|)$.  

(3) The word $w'$ in the alphabet $\mathcal{A}'$ is defined by $w' = xA_{\nu}yA_{\nu}$ for $w = AxAy$.  
\end{definition}

We define $\nu$-shifts and $\nu$-permutations of words in a nanophrase $P$ $=$ $(\mathcal{A},$ $w_{1}|w_{2}|$ $\cdots$ $|w_{k})$ over $\alpha$ and define $\mathcal{P}(\alpha, S, \nu)$ in the following manner as in \cite[Subsection 6.2]{turaev2}.  

Fix a homotopy data $(\alpha, \tau, S)$ and a shift involution in $\alpha$.  
\begin{definition}
For $i$ $=$ $1, \dots, k$, the $i$th $\nu$-{\it shift} of a nanophrase $P$ moves the first letter, say $A$, of $w_{i}$ to the end of $w_{i}$, keeping $|A| \in \alpha$ if $A$ appears in $w_{i}$ only once and applying $\nu$ if $A$ appears in
 $w_{i}$ twice.  All other words in $P$ are preserved.  
\end{definition}

\begin{definition}
Given two words $u$, $v$ on an $\alpha$-alphabet $\mathcal{A}$, consider the mapping $\mathcal{A} \to \alpha$ sending $A \in \mathcal{A}$ to $\nu(|A|) \in \alpha$ if $A$ appears both in $u$ and $v$ and sending $A$ to $|A|$ otherwise.  The set $\mathcal{A}$ with this projection to $\alpha$ is an $\alpha$-alphabet denoted by $\mathcal{A}_{u \cap v}$.  For $i$ $=$ $1, \dots, k-1$, the $\nu$-{\it permutation} of the $i$th and $(i+1)$st words transforms a nanophrase $P$ $=$ $(\mathcal{A}, w_{1}|w_{2}|\cdots|w_{k})$ into the nanophrase $(\mathcal{A}, w_{1}|w_{2}|\cdots|w_{i-1}|w_{i+1}|w_{i}|w_{i+2}|\cdots|w_{k})$.  The operation is involutive.  The $\nu$-permu-tations define an action of the symmetric group $S_{k}$ on the set of nanophrases of length $k$.  
\end{definition}

$\mathcal{P}(\alpha, S, \nu)$ denotes the set of nanophrases over $\alpha$ quotiented by the equivalence relation generated by $S$-homotopy, $\nu$-permutations and, $\nu$-shifts on words.  

Turaev defined {\it pseudolinks} in the following manner as in \cite[Subsection 7.1]{turaev2}.  
\begin{definition}
Let $\alpha_{1}$ $=$ $\{-1, 1\}$ with involution $\tau$ permuting $1$ and $-1$ and let $S_{1} \subset \alpha_{1} \times \alpha_{1} \times \alpha_{1}$ consists of the following six triples: $(1, 1, 1)$, $(1, 1, -1)$, $(-1, 1, 1)$, $(-1, -1,- 1)$, and $(-1, -1, 1)$, $(1, -1, -1)$.  Let $\nu$ $=$ $\operatorname{id}$.  
Nanophrases in $\mathcal{P}(\alpha_{1}, S_{1}, \operatorname{id})$ are called {\it pseudolinks}.  
\end{definition}

\begin{remark}
Let $\alpha_{*}$ be a set consisting of $4$ distinct elements $a_{+}, a_{-},  b_{+},$ and $b_{-}$ with involution $\tau:$ $a_{\pm} \mapsto b_{\mp}$.  Let $S_{*}$ $=\{(a_{\pm},$ $a_{\pm}, a_{\pm}),$ $(a_{\pm}, a_{\pm}, a_{\mp}),$ $(a_{\mp}, a_{\pm}, a_{\pm}),$ $(b_{\pm}, b_{\pm}, b_{\pm}),$ $(b_{\pm}, b_{\pm}, b_{\mp}),$ $(b_{\mp}, b_{\pm},$ $b_{\pm})\}$.  A projection $\alpha_{*} \to$ $\alpha_{1}:=\{1, -1\};$ $a_{+}, b_{+} \mapsto 1$ and $a_{-}, b_{-} \mapsto -1$ induces a surjective mapping ${\mathcal P}(\alpha_{*}, S_{*}, \nu)$ $\to$ ${\mathcal P}(\alpha_{1}, S_{1}, \operatorname{id})$.  
\end{remark}

In the last part of this section, we describe the notation $\mathcal{A}_{w}$ as in \cite[Subsection 6.2]{turaev2} and the notation $P_{w}$ as in \cite[Subsection 8.2]{turaev2}.  

\begin{definition}
For a word $w$, $\mathcal{A}_{w}$ denotes the same alphabet $\mathcal{A}$ with a new projection $|\cdots|_{w}$ to $\alpha$ defined as follows: for $A \in \mathcal{A}$, set $|A|_{w}$ $=$ $\tau(|A|)$ if $A$ occurs once, $|A|_{w}$ $=$ $\nu(|A|)$ if $A$ occurs twice, and $|A|_{w}$ $=$ $|A|$ otherwise.  
For a phrase $P$ in an $\alpha_{1}$-alphabet $\mathcal{A}$ and a word $w$ on $\mathcal{A}$, $P_{w}$ denotes the same phrase on the $\alpha_{1}$-alphabet $\mathcal{A}_{w}$.  
\end{definition}

\section{Jones polynomial for pseudolinks}\label{jones-by-word}

Turaev defined the Jones polynomial for pseudolinks by using recursive relations for the bracket polynomial of nanophrases over $\alpha_{*}$ \cite[Section 8]{turaev2}.  
In this section, we present a state sum representation of the Jones polynomial for pseudolinks.

\begin{definition}
For every pseudolink $P$ $=$ $(\mathcal{A},$ $w_{1}|w_{2}|\cdots|w_{k})$, we assign a sign $-1$ or $1$ to each letter $A$ and call the sign the {\it marker} of $A$, denoted by $mark(A)$.  Let a {\it state} $s$ of $P$ be $P$ with their markers for all the elements of $\mathcal{A}$.  
\end{definition}

For an arbitrary pseudolink $P$ assigned with state $s$, we consider the following deformation ($\ast$): 

\begin{equation*}
(\ast)  \begin{cases}
          & w_1|\cdots|AxAy|\cdots|w_{k} \to \begin{cases}
          & w_1|\cdots|x|y|\cdots|w_{k}~\text{if}~{\rm mark}(A) = |A|\\
          & (w_1|\cdots|x^{-}y|\cdots|w_{k})_{x}~\text{if}~{\rm mark}(A) = -|A|
       \end{cases}
       \\
          & w_1|\cdots|Ax|Ay|\cdots|w_{k} \to \begin{cases}
          & w_1|\cdots|xy|\cdots|w_{k}~\text{if}~{\rm mark}(A) = |A|\\
          & (w_1|\cdots|x^{-}y|\cdots|w_{k})_{x}~\text{if}~{\rm mark}(A) = -|A|.  
       \end{cases}
       \end{cases}
\end{equation*}

A pseudolink $\emptyset|\cdots|\emptyset$ is obtained by repeating these deformations from $P$.  We denote the length of this pseudolink $\emptyset|\cdots|\emptyset$ by $|s|$. 

\begin{definition}
We denote a letter $A$ with $|A| = 1$ and ${\rm mark}(A) = +1$ 
(respectively ${\rm mark}(A) = -1$) by $A_{+}$ (respectively $A_{-}$),
and we denote a letter $A$ with $|A|=-1$ and ${\rm mark}(A) = +1$ 
(respectively ${\rm mark}(A) = -1$) by $\overline{A}_{+}$ 
(respectively $\overline{A}_{-}$).
\end{definition}
\begin{example}
Consider $P$ $=$ $ABAB$ with $|A|$ $=$ $|B|$ $=$ $1$.  If mark($A$) $=$ $1$ and mark($B$) $=$ $-1$, $P$ is represented as $A_{+}B_{-}A_{+}B_{-}$ and 
\begin{equation}
\begin{split}
A_{+}B_{-}A_{+}B_{-} &\stackrel{(\ast)}{\rightarrow} B_{-}|B_{-}\\
 &\stackrel{(\ast)}{\rightarrow} \emptyset.  
\end{split}
\end{equation} 
If $P$ has mark($A$) $=$ $1$ and mark($B$) $=$ $-1$, 
\begin{equation}
\begin{split}
A_{-}B_{+}A_{-}B_{+} &\stackrel{(\ast)}{\rightarrow} \overline{B}_{+}\overline{B}_{+}\\ &\stackrel{(\ast)}{\rightarrow} \emptyset.  
\end{split}
\end{equation}
\end{example}

\begin{example}\label{3LetterEx}
Let us add two more examples.  
\begin{equation}
\begin{split}
A_{+}\overline{B}_{+}A_{+}C_{+}\overline{B}_{+}C_{+} &\stackrel{(\ast)}{\rightarrow} \overline{B}_{+}|C_{+}\overline{B}_{+}C_{+}\\ &\stackrel{(\ast)}{\rightarrow} C_{+}C_{+}\\ &\stackrel{(\ast)}{\rightarrow}\emptyset|\emptyset.  
\end{split}
\end{equation}
\begin{equation}
\begin{split}
A_{-}\overline{B}_{-}A_{-}C_{+}\overline{B}_{-}C_{+} &\stackrel{(\ast)}{\rightarrow}{B_{-}}C_{+}{B_{-}}C_{+}\\ &\stackrel{(\ast)}{\rightarrow}\overline{C}_{+}\overline{C}_{+}\\ &\stackrel{(\ast)}{\rightarrow}\emptyset.  
\end{split}
\end{equation}
\end{example}

%%%%%%%%%%%%%%%%%%%%%%%%%%%%%%%%%%%%%%%%%%%%%%%%%%%%%%%%%%%%%%%%%%%%%%%%%%%%%%%%%%%%
\begin{lemma}\label{well-def}
$|s|$ is well defined. In other words, $|s|$ does not depend on the order in which letters are deleted.  
\end{lemma}
\begin{proof}
By the definition of the deformation, 
it is sufficient to consider the cases which do not 
contain an overline.
For such cases, we obtain Table \ref{table1}.  
\begin{table}
\caption{All the types of the nanophrases that should be checked are arranged in the left column.  For each nanophrase in the same and left line, the nanophrase is got after we delete a letter $A$ and then $B$ in each line of the center column.  The counterparts of the center column that are got after we delete $B$ and $A$ in the left column.}
{\begin{tabular}{|c|c|c|}
\hline
cases & deleting $A$ then $B$ &  deleting $B$ then $A$\\
\hline
$A_{+}xA_{+}yB_{+}zB_{+}t$ &  $x|z|ty$ &  $z|x|yt$ \\
$A_{-}xA_{-}yB_{+}zB_{+}t$ &  $z|tx^{-}y$ &  $z|x^{-}yt$ \\
$A_{+}xA_{+}yB_{-}zB_{-}t$ & $z|x|yt$ & $x|z|ty$ \\
$A_{-}xA_{-}yB_{-}zB_{-}t$ &  $z^{-}tx^{-}y$ & $ x^{-}yz^{-}t$ \\
$A_{+}xB_{+}yA_{+}zB_{+}t$  &  $yxtz$  & $ zyxt$ \\
$A_{-}xB_{+}yA_{-}zB_{+}t$ & $z^{-}xty^{-}$ &  $y^{-}z^{-}xt$ \\
$A_{+}xB_{-}yA_{+}zB_{-}t$ & $x^{-}y^{-}tz$  & $t^{-}yxz^{-}$ \\
$A_{-}xB_{-}yA_{-}zB_{-}t$ &  $x^{-}z|ty^{-}$ &  $yt|xz^{-}$\\
$A_{+}xA_{+}y|B_{+}zB_{+}t$ & $ x|y|z|t$ &  $x|y|z|t$ \\
$A_{-}xA_{-}y|B_{+}zB_{+}t$ & $x^{-}y|z|t$  &  $ x^{-}y|z|t$ \\
$A_{+}xA_{+}y|B_{-}zB_{-}t$ & $x|y|z^{-}t$ & $x|y|z^{-}t$  \\
$A_{-}xA_{-}y|B_{-}zB_{-}t$ & $x^{-}y|z^{-}t$  & $ x^{-}y|z^{-}t$ \\
$A_{+}xB_{+}y|A_{+}zB_{+}t$ & $yz|tx$  &  $ xt|zy$  \\
$A_{-}xB_{+}y|A_{-}zB_{+}t$ & $z^{-}xty^{-}$  & $ t^{-}x^{-}zy$ \\
$A_{+}xB_{-}y|A_{+}zB_{-}t$ &  $z^{-}y^{-}tx$  &  $t^{-}yzx^{-}$  \\
$A_{-}xB_{-}y|A_{-}zB_{-}t$ &  $ x^{-}z|ty^{-}$  &  $ y^{-}t|zx^{-}$  \\
$A_{+}x|A_{+}yB_{+}zB_{+}t$ & $z|txy$  & $ xyt|z$  \\
$A_{-}x|A_{-}yB_{+}zB_{+}t$ & $z|tx^{-}y$ & $ x^{-}yt|z$ \\
$A_{+}x|A_{+}yB_{-}zB_{-}t$ & $z^{-}txy$ & $ xyz^{-}t$ \\
$A_{-}x|A_{-}yB_{-}zB_{-}t$ & $z^{-}tx^{-}y$ & $x^{-}yz^{-}t$ \\
$A_{+}x|A_{+}y|B_{+}zB_{+}t$ & $ xy|z|t$  & $ xy|z|t$  \\
$A_{-}x|A_{-}y|B_{+}zB_{+}t$ & $x^{-}y|z|t$ &  $x^{-}y|z|t$ \\ 
$A_{+}x|A_{+}y|B_{-}zB_{-}t$ & $xy|z^{-}t$ & $xy|z^{-}t$ \\
$A_{-}x|A_{-}y|B_{-}zB_{-}t$ & $ x^{-}y|z^{-}t$ & $ x^{-}y|z^{-}t$ \\
$A_{+}x|B_{+}y|A_{+}zB_{+}t$ & $ txzy$ & $xzyt$  \\
$A_{-}x|B_{+}y|A_{-}zB_{+}t$ & $ tx^{-}zy$ & $x^{-}zyt$  \\
$A_{+}x|B_{-}y|A_{+}zB_{-}t$ & $xzy^{-}t$  & $y^{-}txz$  \\
$A_{-}x|B_{-}y|A_{-}zB_{-}t$ & $ z^{-}xt^{-}y$  & $ t^{-}yz^{-}x$  \\
$A_{+}x|A_{+}y|B_{+}z|B_{+}t$ & $xy|zt$  & $xy|zt$  \\
$A_{-}x|A_{-}y|B_{+}z|B_{+}t$ & $x^{-}y|z|t$ & $x^{-}y|z|t$ \\
$A_{+}x|A_{+}y|B_{-}z|B_{-}t$ & $xy|z^{-}t$ & $xy|z^{-}t$ \\
$A_{-}x|A_{-}y|B_{-}z|B_{-}t$ & $x^{-}y|z^{-}t$ & $x^{-}y|z^{-}t$ \\
\hline
\end{tabular}}
\label{table1}
\end{table}

For example, consider the case $A_{+}xA_{+}yB_{+}zB_{+}t$.  
If we delete $A$ first, then
\begin{eqnarray*}
A_{+}xA_{+}yB_{+}zB_{+}t &\longrightarrow& x|yB_{+}zB_{+}t \\
                        &\longrightarrow& x|B_{+}zB_{+}ty \\
                        &\longrightarrow& x|z|ty \\  
\end{eqnarray*}
If we delete $B$ first, then
\begin{eqnarray*}
A_{+}xA_{+}yB_{+}zB_{+}t &\longrightarrow& B_{+}zB_{+}tA_{+}xA_{+}y \\
                        &\longrightarrow& z|tA_{+}xA_{+}y \\
                        &\longrightarrow& z|A_{+}xA_{+}yt \\  
                        &\longrightarrow& z|x|yt
\end{eqnarray*}  
Thus $|s|$ does not depend on the order of 
deletion of letters.\\

\end{proof}

\begin{remark}
The deformation $(\ast)$ corresponds to smoothing crossings of link diagrams in the following figures (see \cite[Page 320, Figure 1]{viro}).  
\begin{figure}[htbp]
\begin{center}
\begin{minipage}{300pt}
\begin{picture}(40,40)
\put(0,0){\line(1,1){30}}
\qbezier(30,0)(30,0)(20,10)
\qbezier(10,20)(0,30)(0,30)
\qbezier(40,10)(45,12)(50,14)
\qbezier(50,14)(50,12)(50,10)
\qbezier(50,10)(55,12)(60,14)
\put(60,14){\vector(4,1){0}}
\linethickness{2pt}
\put(15,5){\circle*{5}}
\put(15,25){\circle*{5}}
\put(15,5){\line(0,1){20}}
\end{picture}
\qquad
\begin{picture}(40,40)
\qbezier(0,0)(20,15)(0,30)
\qbezier(30,0)(10,15)(30,30)
\end{picture}
\qquad\qquad
\begin{picture}(40,40)
\put(0,0){\line(1,1){30}}
\qbezier(30,0)(30,0)(20,10)
\qbezier(10,20)(0,30)(0,30)
\qbezier(40,10)(45,12)(50,14)
\qbezier(50,14)(50,12)(50,10)
\qbezier(50,10)(55,12)(60,14)
\put(60,14){\vector(4,1){0}}
\linethickness{2pt}
\put(5,15){\circle*{5}}
\put(25,15){\circle*{5}}
\put(5,15){\line(1,0){20}}
\end{picture}
\qquad
\begin{picture}(40,40)
\qbezier(0,0)(15,20)(30,0)
\qbezier(0,30)(15,10)(30,30)
\end{picture}
\end{minipage}
\label{marker}
\caption{Smoothing of a diagram according to thick segments corresponding to markers.  }
\end{center}
\end{figure}
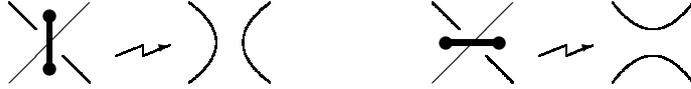
\end{remark}

\begin{definition}
For an arbitrary pseudolink $P$ and state $s$ of $P$, we define $[P]$, $[P|s]$ $\in {\mathbb{Z}}[t, u, d]$ by 
\begin{align}
[P|s] &:= t^{\sharp\{\text{positive marker}\}}u^{\sharp\{\text{negative marker}\}}d^{|s|-1}, \\ 
[P] &:= \sum_{s} [P|s].  
\end{align}
\end{definition}

\begin{proposition}
The polynomial $[P]$ is invariant under an $S_{1}$-homotopy move (H2) for an arbitrary pseudolink $P$ if and only if $u = t^{-1}$ and $d = -t^{2}-t^{-2}$.  
\end{proposition}
\begin{proof}
Consider a nanophrase $P= P_{1}|ABxBAy|P_{2}$ with $|A|=+$ and $|B|=-$, 
where $x$ and $y$ are words not including the character $``|"$ 
Then,
\begin{eqnarray*} 
[P_{1}|ABxBAy|P_{2}]&=&t[P_{1}|BxB|y|P_{2}]+s[(P_{1}|Bx^{-}By|P_{2})_{x}]\\
                     &=&(t^{2}+tsd+s^{2})[(P_{1}|x^{-}|y|P_{2})_{x}]+
                           st[P_{1}|xy|P_{2}]   
\end{eqnarray*}
Thus if $[P]$ does not change by the second homotopy move, 
then $t^{2}+tsd+s^{2}=0$ and $st=1$.
In other words $s=t^{-1}$ and $d=-t^{2}-t^{-2}$.\par
Converse is easily verified  by the above equation.
\end{proof}

\begin{remark}
Substituting $t^{-1}$ for $u$ and $-t^{2}-t^{-2}$ for $d$, we have \[[P]=\sum_{s}t^{\sigma(s)}(-t^{2}-t^{-2})^{|s|-1},\] 
where $\sigma(s)$ $:=$ $\sharp \{$ positive marker $\}$ $-$ $\sharp \{$ negative marker $\}$.  
\end{remark}

\begin{proposition}
$[P]$ is invariant under an $S_{1}$-homotopy move (H3) 
for an arbitrary pseudolink $P$.  
\end{proposition}
\begin{proof}
First, we consider the case of 
$(\epsilon(A),\epsilon(B),\epsilon(C)) = (\pm,\pm,\pm)$.
Consider the third homotopy move
$$P_{1}|ABxACyBCz|P_{2} \longrightarrow P_{1}|BAxCAyCBz|P_{2}.$$
Then, 
\begin{eqnarray*}
[P_{1}|ABxACyBCz|P_{2}] &=& t^{3\epsilon(A)}[P_{1}|xy|z|P_{2}]\\
                          &+&
(2t^{\epsilon(A)}+t^{-3\epsilon(A)}+t^{-\epsilon(A)}(-t^{2}-t^{-2})
)[P_{1}|zx^{-}y^{-}|P_{2}]\\
                           &+&
t^{\epsilon(A)}[P_{1}|x^{-}y|z|P_{2}]\\ 
                           &+&
t^{-\epsilon(A)}[P_{1}|xy^{-}z|P_{2}]\\
&+&
t^{-\epsilon(A)}[P_{1}|x^{-}yz|P_{2}]
\end{eqnarray*}
and
\begin{eqnarray*}
[P_{1}|BAxCAyCBz|P_{2}] &=& t^{3\epsilon(A)}[P_{1}|xy|z|P_{2}]\\
                          &+&\!\!
(2t^{\epsilon(A)}+t^{-3\epsilon(A)}+t^{-\epsilon(A)}(-t^{2}-t^{-2})
)[P_{1}|x^{-}y|z|P_{2}]\\
                           &+&
t^{\epsilon(A)}[P_{1}|x^{-}y^{-}z|P_{2}]\\ 
                           &+&
t^{-\epsilon(A)}[P_{1}|xy^{-}z|P_{2}]\\
&+&
t^{-\epsilon(A)}[P_{1}|z^{-}y^{-}x|P_{2}].
\end{eqnarray*}
Note that 
$$2t^{\epsilon(A)}+t^{-3\epsilon(A)}+t^{-\epsilon(A)}(-t^{2}-t^{-2})
  =t^{\epsilon(A)}$$
Therefore 
$[P_{1}|ABxACyBCz|P_{2}]$ is equal to $[P_{1}|BAxCAyCBz|P_{2}]$.\par
Consider the third homotopy move
$$P_{1}|ABx|ACyBCz|P_{2} \longrightarrow P_{1}|BAx|CAyCBz|P_{2}.$$
Then 
\begin{eqnarray*}
[P_{1}|ABx|ACyBCz|P_{2}] &=& t^{3\epsilon(A)}[P_{1}|xzy|P_{2}]\\
                          &+&
(2t^{\epsilon(A)}+t^{-3\epsilon(A)}+t^{-\epsilon(A)}(-t^{2}-t^{-2})
)[P_{1}|zx^{-}y^{-}|P_{2}]\\
                           &+&
t^{-\epsilon(A)}[P_{1}|x^{-}|y^{-}z|P_{2}]\\ 
                           &+&
t^{\epsilon(A)}[P_{1}|xzy^{-}|P_{2}]\\
&+&
t^{-\epsilon(A)}[P_{1}|x^{-}yz|P_{2}]
\end{eqnarray*}
and
\begin{eqnarray*}
[P_{1}|BAx|CAyCBz|P_{2}] &=& t^{3\epsilon(A)}[P_{1}|xyz|P_{2}]\\
                          &+&
(2t^{\epsilon(A)}+t^{-3\epsilon(A)}+t^{-\epsilon(A)}(-t^{2}-t^{-2})
)[P_{1}|z^{-}x^{-}y|P_{2}]\\
                           &+&
t^{-\epsilon(A)}[P_{1}|x^{-}yz|P_{2}]\\ 
                           &+&
t^{-\epsilon(A)}[P_{1}|x^{-}|y^{-}z|P_{2}]\\
&+&
t^{\epsilon(A)}[P_{1}|x^{-}y^{-}z|P_{2}].
\end{eqnarray*}
Therefore, $[P_{1}|ABx|ACyBCz|P_{2}]$ is equal to $[P_{1}|BAx|CAyCBz|P_{2}]$.\par
Consider the third homotopy move
$$P_{1}|ABxACy|BCz|P_{2} \longrightarrow P_{1}|BAxCAy|CBz|P_{2}.$$
Then 
\begin{eqnarray*}
[P_{1}|ABxACy|BCz|P_{2}] &=& t^{3\epsilon(A)}[P_{1}|xzy|P_{2}]\\
                          &+&
(2t^{\epsilon(A)}+t^{-3\epsilon(A)}+t^{-\epsilon(A)}(-t^{2}-t^{-2})
)[P_{1}|xy^{-}z|P_{2}]\\
                           &+&
t^{-\epsilon(A)}[P_{1}|x^{-}y^{-}z|P_{2}]\\ 
                           &+&
t^{\epsilon(A)}[P_{1}|x^{-}yz|P_{2}]\\
&+&
t^{-\epsilon(A)}[P_{1}|yx^{-}|z|P_{2}]
\end{eqnarray*}
and
\begin{eqnarray*}
[P_{1}|BAx|CAyCBz|P_{2})] &=& t^{3\epsilon(A)}[P_{1}|xyz|P_{2}]\\
                          &+&
(2t^{\epsilon(A)}+t^{-3\epsilon(A)}+t^{-\epsilon(A)}(-t^{2}-t^{-2})
)[P_{1}|x^{-}yz|P_{2}]\\
                           &+&
t^{-\epsilon(A)}[P_{1}|x^{-}y^{-}z|P_{2}]\\ 
                           &+&
t^{-\epsilon(A)}[P_{1}|y^{-}x|z|P_{2}]\\
&+&
t^{\epsilon(A)}[P_{1}|xy^{-}z|P_{2}].
\end{eqnarray*}
Therefore $[P_{1}|ABxACy|BCz|P_{2}]$ is equal to $[P_{1}|BAxCAy|CBz|P_{2}]$.\par
Consider the third homotopy move
$$P_{1}|ABx|ACy|BCz|P_{2} \longrightarrow P_{1}|BAx|CAy|CBz|P_{2}.$$
Then 
\begin{eqnarray*}
[P_{1}|ABx|ACy|BCz|P_{2}] &=& t^{3\epsilon(A)}[P_{1}|y|zx|P_{2}]\\
                          &+&
(2t^{\epsilon(A)}+t^{-3\epsilon(A)}+t^{-\epsilon(A)}(-t^{2}-t^{-2})
)[P_{1}|y^{-}zx|P_{2}]\\
                           &+&
t^{-\epsilon(A)}[P_{1}|x^{-}|zy^{-}|P_{2}]\\ 
                           &+&
t^{\epsilon(A)}[P_{1}|xzy^{-}|P_{2}]\\
&+&
t^{-\epsilon(A)}[P_{1}|yx^{-}|z|P_{2}]
\end{eqnarray*}
and
\begin{eqnarray*}
[P_{1}|BAx|CAy|CBz|P_{2}] &=& t^{3\epsilon(A)}[P_{1}|yz|x|P_{2}]\\
                          &+&
(2t^{\epsilon(A)}+t^{-3\epsilon(A)}+t^{-\epsilon(A)}(-t^{2}-t^{-2})
)[P_{1}|y^{-}xz|P_{2}]\\
                           &+&
t^{-\epsilon(A)}[P_{1}|zy^{-}|x^{-}|P_{2}]\\ 
                           &+&
t^{-\epsilon(A)}[P_{1}|x^{-}y|z|P_{2}]\\
&+&
t^{\epsilon(A)}[P_{1}|z^{-}yx^{-}|P_{2}].
\end{eqnarray*}
Therefore $[P_{1}|ABx|ACy|BCz|P_{2}]$ is equal to $[P_{1}|BAx|CAy|CBz|P_{2}]$.\par
The cases of $(\epsilon(A),\epsilon(B),\epsilon(C)) = (\mp,\pm,\pm)$ 
and $(\epsilon(A),\epsilon(B),\epsilon(C)) = (\pm,\pm,\mp)$ are
proved in a similar way as the above case.
\end{proof}

\begin{proposition}
For an arbitrary pseudolink $P$, the Jones polynomial $J(P)$ for pseudolinks is given as
\begin{equation}\label{jones}
J(P) = (-t)^{-3w(P)}\sum_{s:\text{states}} t^{\sigma(s)} (- t^{2} - t^{-2})^{|s|-1}, 
\end{equation}
where
$w(P)$ $=$ $\sum_{\text{letters}~A~\text{in}~P}|A|$.  
\end{proposition}

\begin{remark}
The Jones polynomial $J(P)$ of a pseudolink $P$ is given by using recursive relations for the bracket polynomial of nanophrases over $\alpha_{*}$ \cite[Section 8]{turaev2}.  The existence of $J(P)$ can be confirmed from geometrical objects (links).  Here, we give this well-definedness by Lemma \ref{well-def} and (\ref{jones}) using only pseudolinks.  
\end{remark}

Definition \ref{enhanced} of {\it enhanced states} is given in the same manner as that in \cite[Page 326, Subsection 4.3]{viro}.  
\begin{definition}\label{enhanced}
An {\it enhanced state} $S$ of pseudolink $P$ implies a collection of markers constituting a state $s$ of $P$ enhanced by an assignment of a plus or minus sign to each of the components $\emptyset|\cdots|\emptyset$.  (Recall that $\emptyset|\cdots|\emptyset$ is obtained by deformations ($\ast$).)  We denote $\emptyset$ with a positive marker $+$ by $\emptyset_{+}$ and $\emptyset$ with a negative marker $-$ by $\emptyset_{-}$.  
\end{definition}

\begin{definition}\label{astast}
We rewrite the deformation ($\ast$) as follows: 
\begin{equation*}
(\ast\ast)  \begin{cases}
          &\!\!\!\! w_1|\cdots|AxAy|\cdots|w_{k} \to \begin{cases}
          &\!\!\!\! w_1|\cdots|ax|ay|\cdots|w_{k}~\text{if}~{\rm mark}(A) = |A|\\
          &\!\!\!\! (w_1|\cdots|ax^{-}ay|\cdots|w_{k})_{x}~\text{if}~{\rm mark}(A) = -|A|
       \end{cases}
       \\
          &\!\!\!\! w_1|\cdots|Ax|Ay|\cdots|w_{k} \to \begin{cases}
          &\!\!\!\! w_1|\cdots|axay|\cdots|w_{k}~\text{if}~{\rm mark}(A) = |A|\\
          &\!\!\!\! (w_1|\cdots|ax^{-}ay|\cdots|w_{k})_{x}~\text{if}~{\rm mark}(A) = -|A|.  
       \end{cases}
       \end{cases}
\end{equation*}
where $a$ is a reminder put on the place of deleting letter $A$ in the case of ($\ast$).  We define $a$ as a letter of a nanophrase where $|a|$ is $|A|$.  
A pseudolink $a_{1}^{1} \cdots a_{n_{1}}^{1}|$ $a_{2}^{1} \cdots a_{n_2}^{2}|\cdots $$|a_{1}^{k'} \cdots a_{n_{k'}}^{k'}$ is given by repeating these deformations ($\ast\ast$) from $P$.  The pseudolinks represents an enhanced state $S$ and then the pseudolink is denoted by ${\emptyset_{\epsilon_1}}^{B_{1}^{1} \cdots B_{m_{1}}^{1}}| {\emptyset_{\epsilon_2}}^{B_{1}^{2} \cdots B_{m_{2}}^{2}}|\cdots| {\emptyset_{\epsilon_{k'}}}^{B_{1}^{k'} \cdots B_{m_{k'}}^{k'}}$, where $B_{1}^{i} \cdots B_{m_{i}}^{i}$ is a word obtained by arranging all the distinct letters in $\{A_{1}^{i}, \ldots, A_{n_{i}}^{i}\}$ in any desired order.  The symbol $B_{y}^{x}$ implies that the letter of $P$ belongs to the $x$-th component of an enhanced state of $P$ and the letter is the $y$-th letter of the $x$-th component by the desired order.  
\end{definition}

\begin{example}
Example \ref{3LetterEx} is rewritten by using Definition \ref{astast}.  
\begin{equation}
\begin{split}
A_{+}\overline{B}_{+}A_{+}C_{+}\overline{B}_{+}C_{+} &\stackrel{(\ast\ast)}{\rightarrow} a\overline{B}_{+}|aC_{+}\overline{B}_{+}C_{+} = \overline{B}_{+}a|\overline{B}_{+}C_{+}aC_{+} \\&\stackrel{(\ast\ast)}{\rightarrow} babC_{+}aC_{+} = C_{+}aC_{+}bab \\&\stackrel{(\ast\ast)}{\rightarrow} ca|cbab = {\emptyset}^{AC}|{\emptyset}^{ABC}.  
\end{split}
\end{equation}
\begin{equation}
\begin{split}
A_{-}\overline{B}_{-}A_{-}C_{+}\overline{B}_{-}C_{+} &\stackrel{(\ast\ast)}{\rightarrow} a{B_{-}}aC_{+}{B_{-}}C_{+} = {B_{-}}aC_{+}{B_{-}}C_{+}a \\ &\stackrel{(\ast\ast)}{\rightarrow} b\overline{C}_{+}ab\overline{C}_{+}a = \overline{C}_{+}ab\overline{C}_{+}ab \\ &\stackrel{(\ast\ast)}{\rightarrow} cbacab = \emptyset^{ABC}.  
\end{split}
\end{equation}
\end{example}

\begin{definition}
For an arbitrary enhanced state $S$ of pseudolink $P$, let
\begin{align}
i(S) &:= \frac{w(P) - \sigma(S)}{2}, \\
\tau(S) &:= \sharp\{\emptyset_{+}~\text{in}~P_{S}\} - \sharp\{\emptyset_{-}~\text{in}~P_{S}\}, \\
j(S) &:= - \frac{\sigma(S) + 2 \tau (S) - 3w(P)}{2} \in {\mathbb{Z}}.  
\end{align}
\end{definition}

Let $s$ be a state of a pseudolink $P$, $S$ be an enhanced state of $P$, and $\hat{J}(P)$ $=$ $(- t^{2} - t^{-2})J(P)$.  
By using the above notations, we have 
\begin{align}
\hat{J}(P) &= (-t)^{-3w(P)}\sum_{\text{states $s$}} t^{\sigma(s)} (- t^{2} - t^{-2})^{|s|}\\
&= (-t)^{-3w(P)}\sum_{\text{enhanced states $S$}} t^{\sigma(S)} (- t^{2})^{\tau(S)}\\
&= \sum_{\text{enhanced states $S$}}(-1)^{w(P)+\tau(S)}~t^{-2 j(S)}\\
&= \sum_{\text{enhanced states $S$}}(-1)^{\frac{w(P) - \sigma(S)}{2}} q^{j(S)} \quad (q = -t^{-2})\\
&= \sum_{\text{enhanced states $S$}}(-1)^{i(S)} q^{j(S)}.  
\end{align}

\begin{remark}\label{s1tos0}
Let $\alpha_{0}$ be a set $\{-1, 1\}$ with an involution $\tau_{0}:$ $\pm 1 \mapsto \mp 1$ and $S_{0}$ be $\{($$-1,$ $-1,$ $-1),$ $(1,$ $1,$ $1)\}$.  Note that every $S_{1}$-homotopy invariant of pseudolinks is an $S_{0}$-homotopy invariant of nanophrases over $\alpha_{0}$ because $S_{0}$ $\subset$ $S_{1}$.  
\end{remark}

\begin{corollary}\label{jones-s0}
$J(P)$ and $\hat{J}(P)$ are $S_{0}$-homotopy invariants for nanophrases $P$ over $\alpha_{0}$.  
\end{corollary}

\section{The Khovanov homology for pseudolinks}

\begin{definition}
For an arbitrary pseudolink $P$, let 
$C(P)$ be a free abelian group generated by the enhanced states of $P$.  
We define the subgroup $C^{i, j}(P)$ of $C(P)$ by
\[C^{i,j}(P) := \left\langle S : \text{enhanced states}~|~j(S) = j, i(S) = i \right\rangle~(i,~j \in \mathbb{Z}).  \]
\end{definition}

\begin{remark}
The Jones polynomial is given as \[\hat{J}(P) = \sum_{j= - \infty}^{\infty}q^{j}\sum_{i= - \infty}^{\infty}(-1)^{i}{\rm rk}C^{i,j}(P).  \]
\end{remark}

Let us define the differential $d$ of bidegree $(1, 0)$ as follows: 
\[d(S) = \sum_{\text{enhanced states}~T} (S:T) T.  \]
In other words, for two arbitrary enhanced states $S$ and $T$, 
we define incidence numbers $(S:T)$.  We define the differential in the manner similar to that in \cite[Section 5]{viro}.  
Assume that the order of letters in the alphabet of a pseudolink $P$ is given.  

\begin{definition}
The incidence number $(S:T)$ is zero unless the markers of $S$ and $T$ differ at only one letter of $P$; this letter is called {\it the different part between $S$ and $T$}.  The marker of $S$ is positive and that of $T$ is negative at this different part.  
If $(S:T) \neq 0$, the different part between $S$ and $T$ satisfies one of the six cases (\ref{diff1})--(\ref{diff6}) in the following: 

\begin{align}\label{diff1}
S &: \emptyset_{\epsilon_{1}}|\cdots|\emptyset_{\epsilon_{l-1}}|{\emptyset_{-}}^{\cdots A \cdots}|{\emptyset_{-}}^{\cdots A \cdots}|\emptyset_{\epsilon_{l+1}}|\cdots|\emptyset_{\epsilon_{k}}\\ \nonumber
&\begin{minipage}{20pt}
\begin{picture}(20,15)
\qbezier(0,5)(5,7)(10,9)
\qbezier(10,9)(10,7)(10,5)
\qbezier(10,5)(15,7)(20,9)
\put(20,9){\vector(4,1){0}}
\end{picture}
\end{minipage}~T : \emptyset_{\epsilon_{1}}|\cdots|\emptyset_{\epsilon_{l-1}}|{\emptyset_{-}}^{\cdots A \cdots}|\emptyset_{\epsilon_{l+1}}|\cdots|\emptyset_{\epsilon_{k}}; 
\end{align} 
\begin{align}\label{diff2}
S &: \emptyset_{\epsilon_{1}}|\cdots|{\emptyset_{-}}^{\cdots A \cdots}|{\emptyset_{+}}^{\cdots A \cdots}|\emptyset_{\epsilon_{l+1}}|\cdots|\emptyset_{\epsilon_{k}}\\ \nonumber
&\begin{minipage}{20pt}
\begin{picture}(20,15)
\qbezier(0,5)(5,7)(10,9)
\qbezier(10,9)(10,7)(10,5)
\qbezier(10,5)(15,7)(20,9)
\put(20,9){\vector(4,1){0}}
\end{picture}
\end{minipage}~T : \emptyset_{\epsilon_{1}}|\cdots|\emptyset_{\epsilon_{l-1}}|{\emptyset_{+}}^{\cdots A \cdots}|\emptyset_{\epsilon_{l+1}}|\cdots|\emptyset_{\epsilon_{k}}; 
\end{align}\begin{align}\label{diff3}
S &: \emptyset_{\epsilon_{1}}|\cdots|\emptyset_{\epsilon_{l-1}}|{\emptyset_{+}}^{\cdots A \cdots}|{\emptyset_{-}}^{\cdots A \cdots}|\emptyset_{\epsilon_{l+1}}|\cdots|\emptyset_{\epsilon_{k}}\\ \nonumber
&\begin{minipage}{20pt}
\begin{picture}(20,15)
\qbezier(0,5)(5,7)(10,9)
\qbezier(10,9)(10,7)(10,5)
\qbezier(10,5)(15,7)(20,9)
\put(20,9){\vector(4,1){0}}
\end{picture}
\end{minipage}~T : \emptyset_{\epsilon_{1}}|\cdots|\emptyset_{\epsilon_{l-1}}|{\emptyset_{+}}^{\cdots A \cdots}|\emptyset_{\epsilon_{l+1}}|\cdots|\emptyset_{\epsilon_{k}}; 
\end{align}
\begin{align}\label{diff4}
S &: \emptyset_{\epsilon_{1}}|\cdots|\emptyset_{\epsilon_{l-1}}|{\emptyset_{+}}^{\cdots A \cdots}|\emptyset_{\epsilon_{l+1}}|\cdots|\emptyset_{\epsilon_{k}}\\ \nonumber
&\begin{minipage}{20pt}
\begin{picture}(20,15)
\qbezier(0,5)(5,7)(10,9)
\qbezier(10,9)(10,7)(10,5)
\qbezier(10,5)(15,7)(20,9)
\put(20,9){\vector(4,1){0}}
\end{picture}
\end{minipage}~T : \emptyset_{\epsilon_{1}}|\cdots|\emptyset_{\epsilon_{l-1}}|{\emptyset_{+}}^{\cdots A \cdots}|{\emptyset_{+}}^{\cdots A \cdots}|\emptyset_{\epsilon_{l+1}}|\cdots|\emptyset_{\epsilon_{k}}; 
\end{align}
\begin{align}\label{diff5}
S &: \emptyset_{\epsilon_{1}}|\cdots|\emptyset_{\epsilon_{l-1}}|{\emptyset_{-}}^{\cdots A \cdots}|\emptyset_{\epsilon_{l+1}}|\cdots|\emptyset_{\epsilon_{k}}\\ \nonumber
&\begin{minipage}{20pt}
\begin{picture}(20,15)
\qbezier(0,5)(5,7)(10,9)
\qbezier(10,9)(10,7)(10,5)
\qbezier(10,5)(15,7)(20,9)
\put(20,9){\vector(4,1){0}}
\end{picture}
\end{minipage}~T : \emptyset_{\epsilon_{1}}|\cdots|\emptyset_{\epsilon_{l-1}}|{\emptyset_{+}}^{\cdots A \cdots}|{\emptyset_{-}}^{\cdots A \cdots}|\emptyset_{\epsilon_{l+1}}|\cdots|\emptyset_{\epsilon_{k}}; 
\end{align}\begin{align}\label{diff6}
S &: \emptyset_{\epsilon_{1}}|\cdots|\emptyset_{\epsilon_{l-1}}|{\emptyset_{-}}^{\cdots A \cdots}|\emptyset_{\epsilon_{l+1}}|\cdots|\emptyset_{\epsilon_{k}}\\ \nonumber
&\begin{minipage}{20pt}
\begin{picture}(20,15)
\qbezier(0,5)(5,7)(10,9)
\qbezier(10,9)(10,7)(10,5)
\qbezier(10,5)(15,7)(20,9)
\put(20,9){\vector(4,1){0}}
\end{picture}
\end{minipage}~T : \emptyset_{\epsilon_{1}}|\cdots|\emptyset_{\epsilon_{l-1}}|{\emptyset_{-}}^{\cdots A \cdots}|{\emptyset_{+}}^{\cdots A \cdots}|\emptyset_{\epsilon_{l+1}}|\cdots|\emptyset_{\epsilon_{k}}; 
\end{align}
\end{definition}
For (\ref{diff1})--(\ref{diff6}), $(S : T)$ is defined as 
\begin{equation}
(S : T) := 1.  
\end{equation}

\begin{theorem}
$d \circ d =0$ \quad modulo $2$.  
\end{theorem}
\begin{proof}
Let $\epsilon_{i}$ be the $i$th marker of the $i$th letters, 
and so, $\epsilon_{i}$ is an element of $\{+,-\}$.  Consider the $k$-tuple $(\epsilon_{1}, \epsilon_{2}, \dots, \epsilon_{k})$ consisting of all the markers of a phrase.  
If card $\{j~|~\epsilon_{j} = + \}$ $\le 1$, $d^{2}(S)$ $=$ $0$.  
Thus we can assume that card $\{j~|~\epsilon_{j} = + \}$ $\ge 2$ now.\par
To prove $$d \circ d (S) = \sum_{\text{enhanced states}~T,U} (S:T)(T:U)U = 0,$$
we show $\sum_{\text{enhanced states}~T}(S:T)(T:U)=0$.\par
Let $A$ and $B$ be different parts between $S$ and $U$. 
We can assume that the other letters in the phrase have 
already been deleted by the deformation ($\ast \ast$).
We denote phrases consisting of letters replaced by the 
deformation ($\ast \ast$) as $\alpha_{j}$ ($j \in \{1, \cdots k\}$),
$x$, $y$, $z$, and $t$. We denote a state $S$ by 
$S =$ (a phrase $P$ with markers, a pseudolink given by repeating the 
deformation
($\ast \ast$) from $P$ to the end).
We verify  the following 26 cases:\\
(1)
\begin{eqnarray*}
S&=&(\alpha_{1}|\cdots|\alpha_{l-1}|A_{+}xA_{+}yB_{+}zB_{+}t|
\alpha_{l+1}|\cdots|\alpha_{k},\\
&&\emptyset_{\epsilon_{1}}|\cdots|\emptyset_{\epsilon_{l-1}}|{\emptyset_{\varepsilon_{11}}}^{\cdots A \cdots}|{\emptyset_{\varepsilon_{12}}}^{\cdots AB \cdots}|{\emptyset_{\varepsilon_{13}}}^{\cdots B \cdots}|\emptyset_{\epsilon_{l+1}}|\cdots|\emptyset_{\epsilon_{k}})
\end{eqnarray*}
(2) 
\begin{eqnarray*}
S&=&(\alpha_{1}|\cdots|\alpha_{l-1}|\overline{A}_{+}x\overline{A}_{+}yB_{+}zB_{+}t|\alpha_{l+1}|\cdots|\alpha_{k},\\
&&\emptyset_{\epsilon_{1}}|\cdots|\emptyset_{\epsilon_{l-1}}|{\emptyset_{\varepsilon_{11}}}^{\cdots AB \cdots}|{\emptyset_{\varepsilon_{12}}}^{\cdots B \cdots}|\emptyset_{\epsilon_{l+1}}|\cdots|\emptyset_{\epsilon_{k}})
\end{eqnarray*}
(3) 
\begin{eqnarray*}
S&=&(\alpha_{1}|\cdots|\alpha_{l-1}|\overline{A}_{+}x\overline{A}_{+}y
\overline{B}_{+}z\overline{B}_{+}t|\alpha_{l+1}|\cdots|\alpha_{k},\\
&&\emptyset_{\epsilon_{1}}|\cdots|\emptyset_{\epsilon_{l-1}}|{\emptyset_{\varepsilon_{11}}}^{\cdots AB \cdots}|\emptyset_{\epsilon_{l+1}}|\cdots|\emptyset_{\epsilon_{k}})
\end{eqnarray*}
(4) 
\begin{eqnarray*}
S&=&(\alpha_{1}|\cdots|\alpha_{l-1}|A_{+}xB_{+}yA_{+}zB_{+}t|
\alpha_{l+1}|\cdots|\alpha_{k},\\
&&\emptyset_{\epsilon_{1}}|\cdots|\emptyset_{\epsilon_{l-1}}|{\emptyset_{\varepsilon_{11}}^{\cdots AB \cdots}}|\emptyset_{\epsilon_{l+1}}|\cdots|\emptyset_{\epsilon_{k}})
\end{eqnarray*}
(5) 
\begin{eqnarray*}
S&=&(\alpha_{1}|\cdots|\alpha_{l-1}|\overline{A}_{+}xB_{+}y
\overline{A}_{+}zB_{+}t|\alpha_{l+1}|\cdots|\alpha_{k},\\
&&\emptyset_{\epsilon_{1}}|\cdots|\emptyset_{\epsilon_{l-1}}|{\emptyset_{\varepsilon_{11}}}^{\cdots A B \cdots}|\emptyset_{\epsilon_{l+1}}|\cdots|\emptyset_{\epsilon_{k}})
\end{eqnarray*}
(6)
\begin{eqnarray*}
S&=&(\alpha_{1}|\cdots|\alpha_{l-1}|\overline{A}_{+}x\overline{B}_{+}y
\overline{A}_{+}z\overline{B}_{+}t|\alpha_{l+1}|\cdots|\alpha_{k},\\
&&\emptyset_{\epsilon_{1}}|\cdots|\emptyset_{\epsilon_{l-1}}|{\emptyset_{\varepsilon_{11}}}^{\cdots A B \cdots}|{\emptyset_{\varepsilon_{12}}}^{\cdots A B \cdots}|\emptyset_{\epsilon_{l+1}}|\cdots|\emptyset_{\epsilon_{k}})
\end{eqnarray*}
(7) 
\begin{eqnarray*}
S&=&(\alpha_{1}|\cdots|\alpha_{l-1}|A_{+}xA_{+}y|B_{+}zB_{+}t|
\alpha_{l+1}|\cdots|\alpha_{k},\\
&&\emptyset_{\epsilon_{1}}|\cdots|\emptyset_{\epsilon_{l-1}}|{\emptyset_{\varepsilon_{11}}}^{\cdots A \cdots}|{\emptyset_{\varepsilon_{12}}}^{\cdots A \cdots}|{\emptyset_{\varepsilon_{13}}}^{\cdots B \cdots}|{\emptyset_{\varepsilon_{14}}}^{\cdots B \cdots}|\emptyset_{\epsilon_{l+1}}|\cdots|\emptyset_{\epsilon_{k}})
\end{eqnarray*}
(8) 
\begin{eqnarray*}
S&=&(\alpha_{1}|\cdots|\alpha_{l-1}|\overline{A}_{+}x
\overline{A}_{+}y|B_{+}zB_{+}t|\alpha_{l+1}|\cdots|\alpha_{k},\\
&&\emptyset_{\epsilon_{1}}|\cdots|\emptyset_{\epsilon_{l-1}}|{\emptyset_{\varepsilon_{11}}}^{\cdots A \cdots}|{\emptyset_{\varepsilon_{12}}}^{\cdots B \cdots}|{\emptyset_{\varepsilon_{13}}}^{\cdots B \cdots}|\emptyset_{\epsilon_{l+1}}|\cdots|\emptyset_{\epsilon_{k}})
\end{eqnarray*}
(9) 
\begin{eqnarray*}
S&=&(\alpha_{1}|\cdots|\alpha_{l-1}|\overline{A}_{+}x\overline{A}_{+}y|
\overline{B}_{+}z\overline{B}_{+}t|\alpha_{l+1}|\cdots|\alpha_{k},\\
&&\emptyset_{\epsilon_{1}}|\cdots|\emptyset_{\epsilon_{l-1}}|{\emptyset_{\varepsilon_{11}}}^{\cdots A \cdots}|{\emptyset_{\varepsilon_{12}}}^{\cdots B \cdots}|\emptyset_{\epsilon_{l+1}}|\cdots|\emptyset_{\epsilon_{k}})
\end{eqnarray*}
(10) 
\begin{eqnarray*}
S&=&(\alpha_{1}|\cdots|\alpha_{l-1}|A_{+}xB_{+}y|A_{+}zB_{+}t|
\alpha_{l+1}|\cdots|\alpha_{k},\\
&&\emptyset_{\epsilon_{1}}|\cdots|\emptyset_{\epsilon_{l-1}}|{\emptyset_{\varepsilon_{11}}}^{\cdots A B \cdots}|{\emptyset_{\varepsilon_{12}}}^{\cdots A B \cdots}|\emptyset_{\epsilon_{l+1}}|\cdots|\emptyset_{\epsilon_{k}})
\end{eqnarray*}
(11) 
\begin{eqnarray*}
S&=&(\alpha_{1}|\cdots|\alpha_{l-1}|\overline{A}_{+}xB_{+}y|
\overline{A}_{+}zB_{+}t|\alpha_{l+1}|\cdots|\alpha_{k},\\
&&\emptyset_{\epsilon_{1}}|\cdots|\emptyset_{\epsilon_{l-1}}|{\emptyset_{\varepsilon_{11}}}^{\cdots A B \cdots}|\emptyset_{\epsilon_{l+1}}|\cdots|\emptyset_{\epsilon_{k}})
\end{eqnarray*}
(12) 
\begin{eqnarray*}
S&=&(\alpha_{1}|\cdots|\alpha_{l-1}|\overline{A}_{+}x\overline{B}_{+}y|
\overline{A}_{+}z\overline{B}_{+}t|\alpha_{l+1}|\cdots|\alpha_{k},\\
&&\emptyset_{\epsilon_{1}}|\cdots|\emptyset_{\epsilon_{l-1}}|{\emptyset_{\varepsilon_{11}}}^{\cdots A B \cdots}|{\emptyset_{\varepsilon_{12}}}^{\cdots A B \cdots}|\emptyset_{\epsilon_{l+1}}|\cdots|\emptyset_{\epsilon_{k}})
\end{eqnarray*}
(13) 
\begin{eqnarray*}
S&=&(\alpha_{1}|\cdots|\alpha_{l-1}|A_{+}|xA_{+}yB_{+}zB_{+}t|
\alpha_{l+1}|\cdots|\alpha_{k},\\
&&\emptyset_{\epsilon_{1}}|\cdots|\emptyset_{\epsilon_{l-1}}|{\emptyset_{\varepsilon_{11}}}^{\cdots A B \cdots}|{\emptyset_{\varepsilon_{12}}}^{\cdots A B \cdots}|\emptyset_{\epsilon_{l+1}}|\cdots|\emptyset_{\epsilon_{k}})
\end{eqnarray*}
(14) 
\begin{eqnarray*}
S&=&(\alpha_{1}|\cdots|\alpha_{l-1}|\overline{A}_{+}x|
\overline{A}_{+}yB_{+}zB_{+}t|\alpha_{l+1}|\cdots|\alpha_{k},\\
&&\emptyset_{\epsilon_{1}}|\cdots|\emptyset_{\epsilon_{l-1}}|{\emptyset_{\varepsilon_{11}}}^{\cdots A B \cdots}|{\emptyset_{\varepsilon_{12}}}^{\cdots A B \cdots}|\emptyset_{\epsilon_{l+1}}|\cdots|\emptyset_{\epsilon_{k}})
\end{eqnarray*}
(15) 
\begin{eqnarray*}
S&=&(\alpha_{1}|\cdots|\alpha_{l-1}|A_{+}x|A_{+}y\overline{B}_{+}z
\overline{B}_{+}t|\alpha_{l+1}|\cdots|\alpha_{k},\\
&&\emptyset_{\epsilon_{1}}|\cdots|\emptyset_{\epsilon_{l-1}}|{\emptyset_{\varepsilon_{11}}}^{\cdots A B \cdots}|\emptyset_{\epsilon_{l+1}}|\cdots|\emptyset_{\epsilon_{k}})
\end{eqnarray*}
(16) 
\begin{eqnarray*}
S&=&(\alpha_{1}|\cdots|\alpha_{l-1}|\overline{A}_{+}x|\overline{A}_{+}y
\overline{B}_{+}z\overline{B}_{+}t|\alpha_{l+1}|\cdots|\alpha_{k},\\
&&\emptyset_{\epsilon_{1}}|\cdots|\emptyset_{\epsilon_{l-1}}|{\emptyset_{\varepsilon_{11}}}^{\cdots A B \cdots}|\emptyset_{\epsilon_{l+1}}|\cdots|\emptyset_{\epsilon_{k}})
\end{eqnarray*}
(17) 
\begin{eqnarray*}
S&=&(\alpha_{1}|\cdots|\alpha_{l-1}|A_{+}x|A_{+}y|B_{+}zB_{+}t|
\alpha_{l+1}|\cdots|\alpha_{k},\\
&&\emptyset_{\epsilon_{1}}|\cdots|\emptyset_{\epsilon_{l-1}}|{\emptyset_{\varepsilon_{11}}}^{\cdots A \cdots}|{\emptyset_{\varepsilon_{12}}}^{\cdots B \cdots}|{\emptyset_{\varepsilon_{13}}}^{\cdots B \cdots}|\emptyset_{\epsilon_{l+1}}|\cdots|\emptyset_{\epsilon_{k}})
\end{eqnarray*}
(18) 
\begin{eqnarray*}
S&=&(\alpha_{1}|\cdots|\alpha_{l-1}|\overline{A}_{+}x|
\overline{A}_{+}y|B_{+}zB_{+}t|\alpha_{l+1}|\cdots|\alpha_{k},\\
&&\emptyset_{\epsilon_{1}}|\cdots|\emptyset_{\epsilon_{l-1}}|{\emptyset_{\varepsilon_{11}}}^{\cdots A \cdots}|{\emptyset_{\varepsilon_{12}}}^{\cdots B \cdots}|{\emptyset_{\varepsilon_{13}}}^{\cdots B \cdots}|\emptyset_{\epsilon_{l+1}}|\cdots|\emptyset_{\epsilon_{k}})
\end{eqnarray*}
(19) 
\begin{eqnarray*}
S&=&(\alpha_{1}|\cdots|\alpha_{l-1}|A_{+}x|A_{+}y|\overline{B}_{+}z
\overline{B}_{+}t|\alpha_{l+1}|\cdots|\alpha_{k},\\
&&\emptyset_{\epsilon_{1}}|\cdots|\emptyset_{\epsilon_{l-1}}|{\emptyset_{\varepsilon_{11}}}^{\cdots A \cdots}|{\emptyset_{\varepsilon_{12}}}^{\cdots B \cdots}|\emptyset_{\epsilon_{l+1}}|\cdots|\emptyset_{\epsilon_{k}})
\end{eqnarray*}
(20) 
\begin{eqnarray*}
S&=&(\alpha_{1}|\cdots|\alpha_{l-1}|\overline{A}_{+}x|\overline{A}_{+}y|
\overline{B}_{+}z\overline{B}_{+}t|\alpha_{l+1}|\cdots|\alpha_{k},\\
&&\emptyset_{\epsilon_{1}}|\cdots|\emptyset_{\epsilon_{l-1}}|{\emptyset_{\varepsilon_{11}}}^{\cdots A \cdots}|{\emptyset_{\varepsilon_{12}}}^{\cdots B \cdots}|\emptyset_{\epsilon_{l+1}}|\cdots|\emptyset_{\epsilon_{k}})
\end{eqnarray*}
(21) 
\begin{eqnarray*}
S&=&(\alpha_{1}|\cdots|\alpha_{l-1}|A_{+}x|B_{+}y|A_{+}zB_{+}t|
\alpha_{l+1}|\cdots|\alpha_{k},\\
&&\emptyset_{\epsilon_{1}}|\cdots|\emptyset_{\epsilon_{l-1}}|{\emptyset_{\varepsilon_{11}}}^{\cdots A B \cdots}|\emptyset_{\epsilon_{l+1}}|\cdots|\emptyset_{\epsilon_{k}})
\end{eqnarray*}
(22) 
\begin{eqnarray*}
S&=&(\alpha_{1}|\cdots|\alpha_{l-1}|\overline{A}_{+}x|B_{+}y|
\overline{A}_{+}zB_{+}t|\alpha_{l+1}|\cdots|\alpha_{k},\\
&&\emptyset_{\epsilon_{1}}|\cdots|\emptyset_{\epsilon_{l-1}}|{\emptyset_{\varepsilon_{11}}}^{\cdots A B \cdots}|\emptyset_{\epsilon_{l+1}}|\cdots|\emptyset_{\epsilon_{k}})
\end{eqnarray*}
(23) 
\begin{eqnarray*}
S&=&(\alpha_{1}|\cdots|\alpha_{l-1}|\overline{A}_{+}x|\overline{B}_{+}y|
\overline{A}_{+}z\overline{B}_{+}t|\alpha_{l+1}|\cdots|\alpha_{k},\\
&&\emptyset_{\epsilon_{1}}|\cdots|\emptyset_{\epsilon_{l-1}}|{\emptyset_{\varepsilon_{11}}}^{\cdots A B \cdots}|\emptyset_{\epsilon_{l+1}}|\cdots|\emptyset_{\epsilon_{k}})
\end{eqnarray*}
(24) 
\begin{eqnarray*}
S&=&(\alpha_{1}|\cdots|\alpha_{l-1}|A_{+}x|A_{+}y|B_{+}z|B_{+}t|
\alpha_{l+1}|\cdots|\alpha_{k},\\
&&\emptyset_{\epsilon_{1}}|\cdots|\emptyset_{\epsilon_{l-1}}|{\emptyset_{\varepsilon_{11}}}^{\cdots A \cdots}|{\emptyset_{\varepsilon_{12}}}^{\cdots B \cdots}|\emptyset_{\epsilon_{l+1}}|\cdots|\emptyset_{\epsilon_{k}})
\end{eqnarray*}
(25) 
\begin{eqnarray*}
S&=&(\alpha_{1}|\cdots|\alpha_{l-1}|\overline{A}_{+}x|
\overline{A}_{+}y|B_{+}z|B_{+}t|\alpha_{l+1}|\cdots|\alpha_{k},\\
&&\emptyset_{\epsilon_{1}}|\cdots|\emptyset_{\epsilon_{l-1}}|{\emptyset_{\varepsilon_{11}}}^{\cdots A \cdots}|{\emptyset_{\varepsilon_{12}}}^{\cdots B \cdots}|\emptyset_{\epsilon_{l+1}}|\cdots|\emptyset_{\epsilon_{k}})
\end{eqnarray*}
(26) 
\begin{eqnarray*}
S&=&(\alpha_{1}|\cdots|\alpha_{l-1}|\overline{A}_{+}x|\overline{A}_{+}y|
\overline{B}_{+}z|\overline{B}_{+}t|\alpha_{l+1}|\cdots|\alpha_{k},\\
&&\emptyset_{\epsilon_{1}}|\cdots|\emptyset_{\epsilon_{l-1}}|{\emptyset_{\varepsilon_{11}}}^{\cdots A \cdots}|{\emptyset_{\varepsilon_{12}}}^{\cdots B \cdots}|\emptyset_{\epsilon_{l+1}}|\cdots|\emptyset_{\epsilon_{k}}).
\end{eqnarray*}
$\bullet$ Consider case (1).  \par
Let 
\begin{eqnarray*}
U&=&(\alpha_{1}|\cdots|\alpha_{l-1}|A_{-}xA_{-}yB_{-}zB_{-}t|
\alpha_{l+1}|\cdots|\alpha_{k},\\
&&\emptyset_{\epsilon_{1}}|\cdots|\emptyset_{\epsilon_{l-1}}|{\emptyset_{\varepsilon_{41}}}^{\cdots A B \cdots}|\emptyset_{\epsilon_{l+1}}|\cdots|\emptyset_{\epsilon_{k}}).
\end{eqnarray*}
It is sufficient to show that for each $(\varepsilon_{11},\varepsilon_{12},
\varepsilon_{13}) \in \{+, -\} \times \{+, -\} \times \{+, - \}$,   
the coefficient of $U$ in $d^{2}(S)$ 
is even for all $\varepsilon_{41} \in \{+, -\}$. 
Hence, for $S$ and $U$, we have to check the total
number of ways to get $U$ from $S$ is even (we denote
 the condition by ($\sharp$)).  
Let us localize the problem
of the difference parts of $S$, $A_{+}$, and $B_{+}$.
Two routes (i) and (ii)  can be found to change $A_{+}$ 
(respectively $B_{+}$) into $A_{-}$  (respectively $B_{-}$) as follows:\\
(i)
\begin{eqnarray*}
S&=&(\alpha_{1}|\cdots|\alpha_{l-1}|A_{+}xA_{+}yB_{+}zB_{+}t|
\alpha_{l+1}|\cdots|\alpha_{k},\\
&&\emptyset_{\epsilon_{1}}|\cdots|\emptyset_{\epsilon_{l-1}}|{\emptyset_{\varepsilon_{11}}}^{\cdots A \cdots}|{\emptyset_{\varepsilon_{12}}}^{\cdots A B \cdots}|{\emptyset_{\varepsilon_{13}}}^{\cdots A \cdots}|\emptyset_{\epsilon_{l+1}}|\cdots|\emptyset_{\epsilon_{k}})
\end{eqnarray*} 
\begin{eqnarray*}
\rightarrow
T&=&(\alpha_{1}|\cdots|\alpha_{l-1}|A_{-}xA_{-}yB_{+}zB_{+}t|
\alpha_{l+1}|\cdots|\alpha_{k},\\
&&\emptyset_{\epsilon_{1}}|\cdots|\emptyset_{\epsilon_{l-1}}|{\emptyset_{\varepsilon_{21}}}^{\cdots A B \cdots}|{\emptyset_{\varepsilon_{22}}}^{\cdots B \cdots}|\emptyset_{\epsilon_{l+1}}|\cdots|\emptyset_{\epsilon_{k}}) 
\rightarrow U,
\end{eqnarray*}
(ii)
\begin{eqnarray*}
S&=&(\alpha_{1}|\cdots|\alpha_{l-1}|A_{+}xA_{+}yB_{+}zB_{+}t|
\alpha_{l+1}|\cdots|\alpha_{k},\\
&&\emptyset_{\epsilon_{1}}|\cdots|\emptyset_{\epsilon_{l-1}}|{\emptyset_{\varepsilon_{11}}}^{\cdots A \cdots}|{\emptyset_{\varepsilon_{12}}}^{\cdots A B \cdots}|{\emptyset_{\varepsilon_{13}}}^{\cdots A \cdots}|\emptyset_{\epsilon_{l+1}}|\cdots|\emptyset_{\epsilon_{k}})
\end{eqnarray*} 
\begin{eqnarray*}
\rightarrow
T&=&(\alpha_{1}|\cdots|\alpha_{l-1}|A_{+}xA_{+}yB_{-}zB_{-}t|
\alpha_{l+1}|\cdots|\alpha_{k},\\
&&\emptyset_{\epsilon_{1}}|\cdots|\emptyset_{\epsilon_{l-1}}|{\emptyset_{\varepsilon_{31}}}^{\cdots A \cdots}|{\emptyset_{\varepsilon_{32}}}^{\cdots A B \cdots}|\emptyset_{\epsilon_{l+1}}|\cdots|\emptyset_{\epsilon_{k}}) 
\rightarrow U.
\end{eqnarray*}
Then the condition ($\sharp$) can also state that the
sum of the contribution of (i) to the coefficient of $U$
and the contribution of (ii) to the coefficient of $U$ is even.\\
Consider case  $(\varepsilon_{11},\varepsilon_{12},\varepsilon_{13}) = (+,+,+)$.\par
In this case $(S,T) = 0$ for all $(\varepsilon_{21},\varepsilon_{22})$ and
$(\varepsilon_{31},\varepsilon_{32})$. Thus the condition ($\sharp$) holds.\\
Consider case $(\varepsilon_{11},\varepsilon_{12},\varepsilon_{13}) = (-,+,+)$.\par
Consider route (i). $(S,T)$ is not equal to $0$ if and only if 
$(\varepsilon_{21},\varepsilon_{22}) = (+,+)$. Then for this $T$, $(T,U) = 0$ 
for all $\varepsilon_{41} \in \{ \pm \}$. On the other hand, in route (ii), 
$(S,T) = 0$ for all $\varepsilon_{31}, \varepsilon_{32} \in \{ \pm \}$.
Thus the condition ($\sharp$) holds.\\
Consider case $(\varepsilon_{11},\varepsilon_{12},\varepsilon_{13}) = (+,-,+)$.\par
Consider route (i) (respectively route (ii)).
 $(S,T)$ is not equal to $0$ if and only if 
$(\varepsilon_{21},\varepsilon_{22}) = (+,+)$ 
(respectively $(\varepsilon_{21},\varepsilon_{22}) = (+,+)$). 
Then for this $T$, $(T,U) = 0$ for all $\varepsilon_{41}$.
Thus the condition ($\sharp$) holds.\\
Consider case $(\varepsilon_{11},\varepsilon_{12},\varepsilon_{13}) = (+,+,-)$.\par
Consider route (i). in this route
$(S,T) = 0$ for all $\varepsilon_{21}, \varepsilon_{22} \in \{ \pm \}$.
On the other hand, in route (ii),
$(S,T)$ is not equal to $0$ if and only if 
$(\varepsilon_{31},\varepsilon_{32}) = (+,+)$. Then for this $T$, $(T,U) = 0$ 
for all $\varepsilon_{41} \in \{ \pm \}$. 
Thus the condition ($\sharp$) holds.\\
Consider case $(\varepsilon_{11},\varepsilon_{12},\varepsilon_{13}) = (+,-,-)$.\par
Consider route (i). $(S,T)$ is not equal to $0$ if and only if 
$(\varepsilon_{21},\varepsilon_{22}) = (+,-)$. Then for this $T$, $(T,U)$ is 
not equal to $0$ if and only if $\varepsilon_{41} = +$.  
Similarly, in route (ii), 
$(S,T)$ is not equal to $0$ for all $(\varepsilon_{31}, \varepsilon_{32})=(+,-)$.
Then for this $T$, $(T,U)$ is 
not equal to $0$ if and only if $\varepsilon_{41} = +$.
Thus the condition ($\sharp$) holds.\\
Consider case $(\varepsilon_{11},\varepsilon_{12},\varepsilon_{13}) = (-,-,+)$.\par
Consider route (i). $(S,T)$ is not equal to $0$ if and only if 
$(\varepsilon_{21},\varepsilon_{22}) = (-,+)$. Then for this $T$, $(T,U)$ is 
not equal to $0$ if and only if $\varepsilon_{41} = +$.  
Similarly, in route (ii), 
$(S,T)$ is not equal to $0$ for all $(\varepsilon_{31}, \varepsilon_{32})=(-,+)$.
Then for this $T$, $(T,U)$ is 
not equal to $0$ if and only if $\varepsilon_{41} = +$.
Thus the condition ($\sharp$) holds.\\
Consider case $(\varepsilon_{11},\varepsilon_{12},\varepsilon_{13}) = (-,-,-)$.\par
Consider route (i). $(S,T)$ is not equal to $0$ if and only if 
$(\varepsilon_{21},\varepsilon_{22}) = (-,-)$. Then for this $T$, $(T,U)$ is 
not equal to $0$ if and only if $\varepsilon_{41} = -$.  
Similarly, in route (ii), 
$(S,T)$ is not equal to $0$ for all $(\varepsilon_{31}, \varepsilon_{32})=(-,-)$.
Then for this $T$, $(T,U)$ is 
not equal to $0$ if and only if $\varepsilon_{41} = -$.
Thus the condition ($\sharp$) holds.\\
$\bullet$ Consider case (2).  \par
Let 
\begin{eqnarray*}
U&=&(\alpha_{1}|\cdots|\alpha_{l-1}|\overline{A}_{-}x\overline{A}_{-}yB_{-}zB_{-}t|
\alpha_{l+1}|\cdots|\alpha_{k},\\
&&\emptyset_{\epsilon_{1}}|\cdots|\emptyset_{\epsilon_{l-1}}|{\emptyset_{\varepsilon_{41}}}^{\cdots A \cdots}|{\emptyset_{\varepsilon_{42}}}^{\cdots A B \cdots}|\emptyset_{\epsilon_{l+1}}|\cdots|\emptyset_{\epsilon_{k}}).
\end{eqnarray*}
It is sufficient to show that for each $(\varepsilon_{11},\varepsilon_{12}) 
\in \{ (\pm, \pm) \}$ where double signs are arbitrary,
the coefficient of $U$ in $d^{2}(S)$ 
is even for all $\varepsilon_{41}, \varepsilon_{42} \in \{\pm\}$. 
Hence, for $S$ and $U$, we have to check the total
number of ways to get $U$ from $S$ is even (we denote
 the condition by ($\sharp$)).  
Let us localize the problem
of the difference parts of $S$, $A_{+}$ and $B_{+}$.
Two routes (i) and (ii)  can be found to change $A_{+}$ 
(respectively $B_{+}$) into $A_{-}$ (respectively $B_{-}$) as follows:\\
(i)
\begin{eqnarray*}
S&=&(\alpha_{1}|\cdots|\alpha_{l-1}|\overline{A}_{+}x\overline{A}_{+}yB_{+}zB_{+}t|
\alpha_{l+1}|\cdots|\alpha_{k},\\
&&\emptyset_{\epsilon_{1}}|\cdots|\emptyset_{\epsilon_{l-1}}|{\emptyset_{\varepsilon_{11}}}^{\cdots A B \cdots}|{\emptyset_{\varepsilon_{12}}}^{\cdots B \cdots}|{\emptyset_{\varepsilon_{13}}}^{\cdots A \cdots}|\emptyset_{\epsilon_{l+1}}|\cdots|\emptyset_{\epsilon_{k}})
\end{eqnarray*} 
\begin{eqnarray*}
\rightarrow
T&=&(\alpha_{1}|\cdots|\alpha_{l-1}|\overline{A}_{-}x\overline{A}_{-}yB_{+}zB_{+}t|
\alpha_{l+1}|\cdots|\alpha_{k},\\
&&\emptyset_{\epsilon_{1}}|\cdots|\emptyset_{\epsilon_{l-1}}|{\emptyset_{\varepsilon_{21}}}^{\cdots A \cdots}|{\emptyset_{\varepsilon_{22}}}^{\cdots A B \cdots}|{\emptyset_{\varepsilon_{23}}}^{\cdots B \cdots}|\emptyset_{\epsilon_{l+1}}|\cdots|\emptyset_{\epsilon_{k}}) 
\rightarrow U,
\end{eqnarray*}
(ii)
\begin{eqnarray*}
S&=&(\alpha_{1}|\cdots|\alpha_{l-1}|\overline{A}_{+}x\overline{A}_{+}yB_{+}zB_{+}t|
\alpha_{l+1}|\cdots|\alpha_{k},\\
&&\emptyset_{\epsilon_{1}}|\cdots|\emptyset_{\epsilon_{l-1}}|{\emptyset_{\varepsilon_{11}}}^{\cdots A B \cdots}|{\emptyset_{\varepsilon_{12}}}^{\cdots B \cdots}|{\emptyset_{\varepsilon_{13}}}^{\cdots A \cdots}|\emptyset_{\epsilon_{l+1}}|\cdots|\emptyset_{\epsilon_{k}})
\end{eqnarray*} 
\begin{eqnarray*}
\rightarrow
T&=&(\alpha_{1}|\cdots|\alpha_{l-1}|\overline{A}_{+}x\overline{A}_{+}yB_{-}zB_{-}t|
\alpha_{l+1}|\cdots|\alpha_{k},\\
&&\emptyset_{\epsilon_{1}}|\cdots|\emptyset_{\epsilon_{l-1}}|{\emptyset_{\varepsilon_{31}}}^{\cdots A B \cdots}|\emptyset_{\epsilon_{l+1}}|\cdots|\emptyset_{\epsilon_{k}}) 
\rightarrow U.
\end{eqnarray*}
Then the condition ($\sharp$) can also state that the
sum of the contribution of (i) to the coefficient of $U$
and the contribution of (ii) to the coefficient of $U$ is even.\\
Consider case  $(\varepsilon_{11},\varepsilon_{12}) = (+,+)$.\par
Consider route (i). Then $(S,T)$ is not equal to $0$ if and only if 
$(\varepsilon_{21},\varepsilon_{22},\varepsilon_{23})=(+,+,+)$. 
For this $T$, $(T,U)=0$ for all $\epsilon_{41}, \varepsilon_{42} \in \{ \pm \}$.
On the other hand, in route (ii), we obtain $(S,T) = 0$ for all 
$\varepsilon_{21}, \varepsilon_{22},\varepsilon_{23} \in \{ \pm \}$.
Thus the condition ($\sharp$) holds.\\
Consider case  $(\varepsilon_{11},\varepsilon_{12}) = (+,-)$.\par
Consider route (i). Then $(S,T)$ is not equal to $0$ if and only if 
$(\varepsilon_{21},\varepsilon_{22},\varepsilon_{23})=(+,+,+)$.
For this $T$, $(T,U)$ is not equal to $0$ if and only if    
$(\varepsilon_{41},\varepsilon_{42}) = (+,+)$.
Consider route (ii). Then $(S,T)$ is not equal to $0$ if and only if
$\varepsilon_{31}=+$. For this $T$,  $(T,U)$ is not equal to $0$ if and only if    
$(\varepsilon_{41},\varepsilon_{42}) = (+,+)$.
Thus the condition ($\sharp$) holds.\\ 
Consider case  $(\varepsilon_{11},\varepsilon_{12}) = (-,+)$.\par
Consider route (i). Then $(S,T)$ is not equal to $0$ if and only if 
$(\varepsilon_{21},\varepsilon_{22},\varepsilon_{23})=(-,+,+)$ or    
$(\varepsilon_{21},\varepsilon_{22},\varepsilon_{23})=(+,-,+)$.
Substitute $(\varepsilon_{21},\varepsilon_{22},\varepsilon_{23})=(-,+,+)$. 
For this $T$, we obtain $(T,U)=0$ for all $\varepsilon_{41}, \varepsilon_{42} \in \{ \pm \}$.
Substitute $(\varepsilon_{21},\varepsilon_{22},\varepsilon_{23})=(+,-,+)$.
For this $T$, $(T,U)$ is not equal to $0$ if and only if    
$(\varepsilon_{41},\varepsilon_{42}) = (+,+)$.
Consider route (ii). Then $(S,T)$ is not equal to $0$ if and only if
$\varepsilon_{31}=+$. For this $T$,  $(T,U)$ is not equal to $0$ if and only if    
$(\varepsilon_{41},\varepsilon_{42}) = (+,+)$.
Thus the condition ($\sharp$) holds. 
Consider case  $(\varepsilon_{11},\varepsilon_{12}) = (-,-)$.\par
Consider route (i). Then $(S,T)$ is not equal to $0$ if and only if 
$(\varepsilon_{21},\varepsilon_{22},\varepsilon_{23})=(-,+,-)$ or    
$(\varepsilon_{21},\varepsilon_{22},\varepsilon_{23})=(+,-,-)$.
Substitute $(\varepsilon_{21},\varepsilon_{22},\varepsilon_{23})=(-,+,-)$. 
For this $T$, $(T,U)$ is not equal to $0$ if and only if    
$(\varepsilon_{41},\varepsilon_{42}) = (-,+)$.
Substitute $(\varepsilon_{21},\varepsilon_{22},\varepsilon_{23})=(+,-,-)$.
For this $T$, $(T,U)$ is not equal to $0$ if and only if    
$(\varepsilon_{41},\varepsilon_{42}) = (+,-)$.
Consider route (ii). Then $(S,T)$ is not equal to $0$ if and only if
$\varepsilon_{31}=-$. For this $T$,  $(T,U)$ is not equal to $0$ if and only if    
$(\varepsilon_{41},\varepsilon_{42}) = (+,-)$ or $(\varepsilon_{41},\varepsilon_{42}) = (-,+)$.
Thus the condition ($\sharp$) holds. \\
$\bullet$ Consider case (3).  \par
Let 
\begin{eqnarray*}
U&=&(\alpha_{1}|\cdots|\alpha_{l-1}|\overline{A}_{-}x\overline{A}_{-}y\overline{B}_{-}z
\overline{B}_{-}t|\alpha_{l+1}|\cdots|\alpha_{k},\\
&&\emptyset_{\epsilon_{1}}|\cdots|\emptyset_{\epsilon_{l-1}}|{\emptyset_{\varepsilon_{41}}}^{\cdots A \cdots}|{\emptyset_{\varepsilon_{42}}}^{\cdots A B \cdots}|{\emptyset_{\varepsilon_{43}}}^{\cdots B \cdots}|\emptyset_{\epsilon_{l+1}}|\cdots|\emptyset_{\epsilon_{k}}).
\end{eqnarray*}
It is sufficient to show that for each $\varepsilon_{11} 
\in \{ \pm \} $ the coefficient of $U$ in $d^{2}(S)$ 
is even for all $\varepsilon_{41}, \varepsilon_{42}, \varepsilon_{43} \in 
\{\pm\}$. 
Hence, for $S$ and $U$, we have to check the total
number of ways to get $U$ from $S$ is even (we denote
 the condition by ($\sharp$)).  
Let us localize the problem
of the difference parts of $S$, $A_{+}$ and $B_{+}$.
Two routes (i) and (ii)  can be found to change $A_{+}$ 
(respectively $B_{+}$) into $A_{-}$  (respectively $B_{-}$) as follows:\\
(i)
\begin{eqnarray*}
S&=&(\alpha_{1}|\cdots|\alpha_{l-1}|\overline{A}_{+}x\overline{A}_{+}y\overline{B}_{+}z
\overline{B}_{+}t|\alpha_{l+1}|\cdots|\alpha_{k},\\
&&\emptyset_{\epsilon_{1}}|\cdots|\emptyset_{\epsilon_{l-1}}|{\emptyset_{\varepsilon_{11}}}^{\cdots A B \cdots}|\emptyset_{\epsilon_{l+1}}|\cdots|\emptyset_{\epsilon_{k}})
\end{eqnarray*} 
\begin{eqnarray*}
\rightarrow
T&=&(\alpha_{1}|\cdots|\alpha_{l-1}|\overline{A}_{-}x\overline{A}_{-}y\overline{B}_{+}z
\overline{B}_{+}t|\alpha_{l+1}|\cdots|\alpha_{k},\\
&&\emptyset_{\epsilon_{1}}|\cdots|\emptyset_{\epsilon_{l-1}}|{\emptyset_{\varepsilon_{21}}}^{\cdots A \cdots}|{\emptyset_{\varepsilon_{22}}}^{\cdots A B \cdots}|\emptyset_{\epsilon_{l+1}}|\cdots|\emptyset_{\epsilon_{k}})
\rightarrow U,
\end{eqnarray*}
(ii)
\begin{eqnarray*}
S&=&(\alpha_{1}|\cdots|\alpha_{l-1}|\overline{A}_{+}x\overline{A}_{+}y\overline{B}_{+}z
\overline{B}_{+}t|\alpha_{l+1}|\cdots|\alpha_{k},\\
&&\emptyset_{\epsilon_{1}}|\cdots|\emptyset_{\epsilon_{l-1}}|{\emptyset_{\varepsilon_{11}}}^{\cdots A B \cdots}|\emptyset_{\epsilon_{l+1}}|\cdots|\emptyset_{\epsilon_{k}})
\end{eqnarray*} 
\begin{eqnarray*}
\rightarrow
T&=&(\alpha_{1}|\cdots|\alpha_{l-1}|\overline{A}_{+}x\overline{A}_{+}y\overline{B}_{-}z
\overline{B}_{-}t|\alpha_{l+1}|\cdots|\alpha_{k},\\
&&\emptyset_{\epsilon_{1}}|\cdots|\emptyset_{\epsilon_{l-1}}|{\emptyset_{\varepsilon_{31}}}^{\cdots A B \cdots}|{\emptyset_{\varepsilon_{32}}}^{\cdots A B \cdots}|\emptyset_{\epsilon_{l+1}}|\cdots|\emptyset_{\epsilon_{k}})
\rightarrow U.
\end{eqnarray*}
Then the condition ($\sharp$) can also state that the
sum of the contribution of (i) to the coefficient of $U$
and the contribution of (ii) to the coefficient of $U$ is even.\\
Consider case $\varepsilon_{11} = +$.\par
Consider route (i). Then $(S,T)$ is not equal to $0$ if and only if 
$(\varepsilon_{21},\varepsilon_{22}) = (+,+)$. 
For this $T$, $(T,U)$ is not equal to $0$ if and only if $(\varepsilon_{41},\varepsilon_{42},
\varepsilon_{43})=(+,+,+)$.
Consider route (ii). Then $(S,T)$ is not equal to $0$ if and only if 
$(\varepsilon_{31},\varepsilon_{32}) = (+,+)$. 
For this $T$, $(T,U)$ is not equal to $0$ if and only if $(\varepsilon_{41},\varepsilon_{42},
\varepsilon_{43})=(+,+,+)$. Thus in this case the condition ($\sharp$) holds.\\
Consider case $\varepsilon_{11} = -$.\par
Consider route (i). In this case $(S,T)$ is not equal to $0$ if and only if 
$(\varepsilon_{21},\varepsilon_{22}) = (+,-)$ or $(\varepsilon_{21},\varepsilon_{22}) = (-,+)$.
Substitute $(\varepsilon_{21},\varepsilon_{22}) = (+,-)$, then for this $T$, $(T,U)$ is not 
equal to $0$ if and only if  $(\varepsilon_{41},\varepsilon_{42},
\varepsilon_{43})=(+,+,-)$. 
Substitute $(\varepsilon_{21},\varepsilon_{22}) = (-,+)$, then for this $T$, $(T,U)$ is not 
equal to $0$ if and only if  $(\varepsilon_{41},\varepsilon_{42},
\varepsilon_{43})=(-,+,+)$ or $(\varepsilon_{41},\varepsilon_{42},
\varepsilon_{43})=(+,-,+)$. On the other hand, in route (ii) 
$(S,T)$ is not equal to $0$ if and only if 
$(\varepsilon_{31},\varepsilon_{32}) = (+,-)$ or 
$(\varepsilon_{31},\varepsilon_{32}) = (-,+)$.
Substitute $(\varepsilon_{31},\varepsilon_{32}) = (+,-)$, then for this $T$, $(T,U)$ is not 
equal to $0$ if and only if  $(\varepsilon_{41},\varepsilon_{42},
\varepsilon_{43})=(+,+,-)$. 
Substitute $(\varepsilon_{31},\varepsilon_{32}) = (-,+)$, then for this $T$, $(T,U)$ is not 
equal to $0$ if and only if  $(\varepsilon_{41},\varepsilon_{42},
\varepsilon_{43})=(-,+,+)$ or $(\varepsilon_{41},\varepsilon_{42},
\varepsilon_{43})=(+,-,+)$.
Thus in this case the condition ($\sharp$) holds.\\
$\bullet$ Consider case (4).  \par
Let 
\begin{eqnarray*}
U&=&(\alpha_{1}|\cdots|\alpha_{l-1}|A_{-}xB_{-}yA_{-}zB_{-}t|\alpha_{l+1}|\cdots|\alpha_{k},\\
&&\emptyset_{\epsilon_{1}}|\cdots|\emptyset_{\epsilon_{l-1}}|{\emptyset_{\varepsilon_{41}}}^{\cdots A B \cdots}|\emptyset_{\epsilon_{l+1}}|\cdots|\emptyset_{\epsilon_{k}}).
\end{eqnarray*}
It is sufficient to show that for each $\varepsilon_{11} 
\in \{ \pm \} $ the coefficient of $U$ in $d^{2}(S)$ 
is even for all $\varepsilon_{41}, \varepsilon_{42}, \varepsilon_{43} \in 
\{\pm\}$. 
Hence, for $S$ and $U$, we have to check the total
number of ways to get $U$ from $S$ is even (we denote
the condition by ($\sharp$)).  
Let us localize the problem
of the difference parts of $S$, $A_{+}$ and $B_{+}$.
Two routes (i) and (ii)  can be found to change $A_{+}$ 
(respectively $B_{+}$) into $A_{-}$  (respectively $B_{-}$) as follows:\\
(i)
\begin{eqnarray*}
S&=&(\alpha_{1}|\cdots|\alpha_{l-1}|A_{+}xB_{+}yA_{+}zB_{+}t|\alpha_{l+1}|\cdots|\alpha_{k},\\
&&\emptyset_{\epsilon_{1}}|\cdots|\emptyset_{\epsilon_{l-1}}|{\emptyset_{\varepsilon_{11}}}^{\cdots A B \cdots}|\emptyset_{\epsilon_{l+1}}|\cdots|\emptyset_{\epsilon_{k}})
\end{eqnarray*} 
\begin{eqnarray*}
\rightarrow
T&=&(\alpha_{1}|\cdots|\alpha_{l-1}|A_{-}xB_{+}yA_{-}zB_{+}t|\alpha_{l+1}|\cdots|\alpha_{k},\\
&&\emptyset_{\epsilon_{1}}|\cdots|\emptyset_{\epsilon_{l-1}}|{\emptyset_{\varepsilon_{21}}}^{\cdots A B \cdots}|\emptyset_{\epsilon_{l+1}}|\cdots|\emptyset_{\epsilon_{k}})\rightarrow U,
\end{eqnarray*}
(ii)
\begin{eqnarray*}
S&=&(\alpha_{1}|\cdots|\alpha_{l-1}|A_{+}xB_{+}yA_{+}zB_{+}t|\alpha_{l+1}|\cdots|\alpha_{k},\\
&&\emptyset_{\epsilon_{1}}|\cdots|\emptyset_{\epsilon_{l-1}}|{\emptyset_{\varepsilon_{11}}}^{\cdots A B \cdots}|\emptyset_{\epsilon_{l+1}}|\cdots|\emptyset_{\epsilon_{k}})
\end{eqnarray*} 
\begin{eqnarray*}
\rightarrow
T&=&(\alpha_{1}|\cdots|\alpha_{l-1}|A_{+}xB_{-}yA_{+}zB_{-}t|\alpha_{l+1}|\cdots|\alpha_{k},\\
&&\emptyset_{\epsilon_{1}}|\cdots|\emptyset_{\epsilon_{l-1}}|{\emptyset_{\varepsilon_{31}}}^{\cdots A B \cdots}|\emptyset_{\epsilon_{l+1}}|\cdots|\emptyset_{\epsilon_{k}})\rightarrow U.
\end{eqnarray*}
Then the condition ($\sharp$) can also state that the
sum of the contribution of (i) to the coefficient of $U$
and the contribution of (ii) to the coefficient of $U$ is even.\\
In this case, it is clear that $(S,T)(T,U) = 0$ for all $T$ by the definition of $d$.\\
$\bullet$ Consider case (5).  \par
Let 
\begin{eqnarray*}
U&=&(\alpha_{1}|\cdots|\alpha_{l-1}|\overline{A}_{-}xB_{-}y\overline{A}_{-}zB_{-}t|\alpha_{l+1}|\cdots|\alpha_{k},\\
&&\emptyset_{\epsilon_{1}}|\cdots|\emptyset_{\epsilon_{l-1}}|{\emptyset_{\varepsilon_{41}}}^{\cdots A B \cdots}|\emptyset_{\epsilon_{l+1}}|\cdots|\emptyset_{\epsilon_{k}}).
\end{eqnarray*}
It is sufficient to show that for each $\varepsilon_{11}  
\in \{\pm\} $ the coefficient of $U$ in $d^{2}(S)$ 
is even for all $\varepsilon_{41} \in \{\pm\}$. 
Hence, for $S$ and $U$, we have to check the total
number of ways to get $U$ from $S$ is even (we denote
the condition by ($\sharp$)).  
Let us localize the problem
of the difference parts of $S$, $A_{+}$ and $B_{+}$.
Two routes (i) and (ii)  can be found to change $A_{+}$ 
(respectively $B_{+}$) into $A_{-}$  (respectively $B_{-}$) as follows:\\
(i)
\begin{eqnarray*}
S&=&(\alpha_{1}|\cdots|\alpha_{l-1}|\overline{A}_{+}xB_{+}y\overline{A}_{+}zB_{+}t|\alpha_{l+1}|\cdots|\alpha_{k},\\
&&\emptyset_{\epsilon_{1}}|\cdots|\emptyset_{\epsilon_{l-1}}|{\emptyset_{\varepsilon_{11}}}^{\cdots A B \cdots}|\emptyset_{\epsilon_{l+1}}|\cdots|\emptyset_{\epsilon_{k}})
\end{eqnarray*} 
\begin{eqnarray*}
\rightarrow
T&=&(\alpha_{1}|\cdots|\alpha_{l-1}|\overline{A}_{-}xB_{+}y\overline{A}_{-}zB_{+}t|\alpha_{l+1}|\cdots|\alpha_{k},\\
&&\emptyset_{\epsilon_{1}}|\cdots|\emptyset_{\epsilon_{l-1}}|{\emptyset_{\varepsilon_{21}}}^{\cdots A B \cdots}|\emptyset_{\epsilon_{l+1}}|\cdots|\emptyset_{\epsilon_{k}})\rightarrow U,
\end{eqnarray*}
(ii)
\begin{eqnarray*}
S&=&(\alpha_{1}|\cdots|\alpha_{l-1}|\overline{A}_{+}xB_{+}y\overline{A}_{+}zB_{+}t|\alpha_{l+1}|\cdots|\alpha_{k},\\
&&\emptyset_{\epsilon_{1}}|\cdots|\emptyset_{\epsilon_{l-1}}|{\emptyset_{\varepsilon_{11}}}^{\cdots A B \cdots}|\emptyset_{\epsilon_{l+1}}|\cdots|\emptyset_{\epsilon_{k}})
\end{eqnarray*} 
\begin{eqnarray*}
\rightarrow
T&=&(\alpha_{1}|\cdots|\alpha_{l-1}|\overline{A}_{+}xB_{-}y\overline{A}_{+}zB_{-}t|\alpha_{l+1}|\cdots|\alpha_{k},\\
&&\emptyset_{\epsilon_{1}}|\cdots|\emptyset_{\epsilon_{l-1}}|{\emptyset_{\varepsilon_{31}}}^{\cdots A B \cdots}|{\emptyset_{\varepsilon_{32}}}^{\cdots A B \cdots}|\emptyset_{\epsilon_{l+1}}|\cdots|\emptyset_{\epsilon_{k}})\rightarrow U.
\end{eqnarray*}
Then, the condition ($\sharp$) can also state that the
sum of the contribution of (i) to the coefficient of $U$
and the contribution of (ii) to the coefficient of $U$ is even.\\
Consider case $\varepsilon_{11} = +$\par
On route (i). we obtain $(S,T) = 0$ for all $T$ by the definition of $d$.
Consider route (ii). In this case $(S,T)$ is not equal to $0$ if and only 
if $(\varepsilon_{31},\varepsilon_{32}) = (+,+)$. For this $T$, $(T,U)=0$
for all $\varepsilon_{41} \in \{ \pm \}$.
Thus in this case the condition ($\sharp$) holds.\\
Consider case $\varepsilon_{11} = -$\par
Consider route (i). Then, $(S,T) = 0$ for all $T$. 
Consider route (ii). Then, $(S,T)$ is not equal to $0$ if and only if 
$(\varepsilon_{31},\varepsilon_{32})=(+,-)$ or     
$(\varepsilon_{31},\varepsilon_{32})=(-,+)$.
Substitute $(\varepsilon_{31},\varepsilon_{32})=(+,-)$. Then, for this $T$, we obtain
$(T,U)$ is not equal to $0$ if and only if $\varepsilon_{41}=+$.
Substitute $(\varepsilon_{31},\varepsilon_{32})=(-,+)$. Then, for this $T$, we obtain
$(T,U)$ is not equal to $0$ if and only if $\varepsilon_{41}=+$.
Thus the condition ($\sharp$) holds.\\
$\bullet$ Consider case (6).  \par
Let 
\begin{eqnarray*}
U&=&(\alpha_{1}|\cdots|\alpha_{l-1}|\overline{A}_{-}x\overline{B}_{-}y\overline{A}_{-}z\overline{B}_{-}t|\alpha_{l+1}|\cdots|\alpha_{k},\\
&&\emptyset_{\epsilon_{1}}|\cdots|\emptyset_{\epsilon_{l-1}}|{\emptyset_{\varepsilon_{41}}}^{\cdots A B \cdots}|\emptyset_{\epsilon_{l+1}}|\cdots|\emptyset_{\epsilon_{k}}).
\end{eqnarray*}
It is sufficient to show that for each $(\varepsilon_{11},\varepsilon_{12})  
\in \{(\pm,\pm)\}$ where double signs are arbitrary,
the coefficient of $U$ in $d^{2}(S)$ 
is even for all $\varepsilon_{41} \in \{\pm\}$. 
Hence, for $S$ and $U$, we have to check the total
number of ways to get $U$ from $S$ is even (we denote
 the condition by ($\sharp$)).  
Let us localize the problem
of the difference parts of $S$, $A_{+}$ and $B_{+}$.
Two routes (i) and (ii)  can be found to change $A_{+}$ 
(respectively $B_{+}$) into $A_{-}$  (respectively $B_{-}$) as follows:\\
(i)
\begin{eqnarray*}
S&=&(\alpha_{1}|\cdots|\alpha_{l-1}|\overline{A}_{+}x\overline{B}_{+}y\overline{A}_{+}z\overline{B}_{+}t|\alpha_{l+1}|\cdots|\alpha_{k},\\
&&\emptyset_{\epsilon_{1}}|\cdots|\emptyset_{\epsilon_{l-1}}|{\emptyset_{\varepsilon_{11}}}^{\cdots A B \cdots}|{\emptyset_{\varepsilon_{12}}}^{\cdots A B \cdots}|\emptyset_{\epsilon_{l+1}}|\cdots|\emptyset_{\epsilon_{k}})
\end{eqnarray*} 
\begin{eqnarray*}
\rightarrow
T&=&(\alpha_{1}|\cdots|\alpha_{l-1}|\overline{A}_{-}x\overline{B}_{+}y\overline{A}_{-}z\overline{B}_{+}t|\alpha_{l+1}|\cdots|\alpha_{k},\\
&&\emptyset_{\epsilon_{1}}|\cdots|\emptyset_{\epsilon_{l-1}}|{\emptyset_{\varepsilon_{21}}}^{\cdots A B \cdots}|\emptyset_{\epsilon_{l+1}}|\cdots|\emptyset_{\epsilon_{k}})\rightarrow U,
\end{eqnarray*}
(ii)
\begin{eqnarray*}
S&=&(\alpha_{1}|\cdots|\alpha_{l-1}|\overline{A}_{+}x\overline{B}_{+}y\overline{A}_{+}z\overline{B}_{+}t|\alpha_{l+1}|\cdots|\alpha_{k},\\
&&\emptyset_{\epsilon_{1}}|\cdots|\emptyset_{\epsilon_{l-1}}|{\emptyset_{\varepsilon_{11}}}^{\cdots A B \cdots}|{\emptyset_{\varepsilon_{12}}}^{\cdots A B \cdots}|\emptyset_{\epsilon_{l+1}}|\cdots|\emptyset_{\epsilon_{k}})
\end{eqnarray*} 
\begin{eqnarray*}
\rightarrow
T&=&(\alpha_{1}|\cdots|\alpha_{l-1}|\overline{A}_{+}x\overline{B}_{-}y\overline{A}_{+}z\overline{B}_{-}t|\alpha_{l+1}|\cdots|\alpha_{k},\\
&&\emptyset_{\epsilon_{1}}|\cdots|\emptyset_{\epsilon_{l-1}}|{\emptyset_{\varepsilon_{31}}}^{\cdots A B \cdots}|{\emptyset_{\varepsilon_{32}}}^{\cdots A B \cdots}|\emptyset_{\epsilon_{l+1}}|\cdots|\emptyset_{\epsilon_{k}})\rightarrow U.
\end{eqnarray*}
Then, the condition ($\sharp$) can also state that the
sum of the contribution of (i) to the coefficient of $U$
and the contribution of (ii) to the coefficient of $U$ is even.\\
In this case, we can easily check that $(S,T)(T,U) = 0$ for all $T$
by the definition of $T$.\\
$\bullet$ Consider case (7).  \par
Let 
\begin{eqnarray*}
U&=&(\alpha_{1}|\cdots|\alpha_{l-1}|A_{-}xA_{-}y|B_{-}zB_{-}t|\alpha_{l+1}|\cdots|\alpha_{k},\\
&&\emptyset_{\epsilon_{1}}|\cdots|\emptyset_{\epsilon_{l-1}}|{\emptyset_{\varepsilon_{41}}}^{\cdots A \cdots}|{\emptyset_{\varepsilon_{42}}}^{\cdots A B \cdots}|\emptyset_{\epsilon_{l+1}}|\cdots|\emptyset_{\epsilon_{k}}).
\end{eqnarray*}
It is sufficient to show that for each $(\varepsilon_{11},\varepsilon_{12},
\varepsilon_{13},\varepsilon_{14}) \in \{(\pm,\pm,\pm,\pm)\}$ 
where double signs are arbitrary,  
the coefficient of $U$ in $d^{2}(S)$ 
is even for all $\varepsilon_{41}, \varepsilon_{42} \in \{\pm\}$. 
Hence, for $S$ and $U$, we have to check the total
number of ways to get $U$ from $S$ is even (we denote
 the condition by ($\sharp$)).  
Let us localize the problem
of the difference parts of $S$, $A_{+}$ and $B_{+}$.
Two routes (i) and (ii)  can be found to change $A_{+}$ 
(respectively $B_{+}$) into $A_{-}$  (respectively $B_{-}$) as follows:\\
(i)
\begin{eqnarray*}
S&=&(\alpha_{1}|\cdots|\alpha_{l-1}|A_{+}xA_{+}y|B_{+}zB_{+}t|\alpha_{l+1}|\cdots|\alpha_{k},\\
&&\emptyset_{\epsilon_{1}}|\cdots|\emptyset_{\epsilon_{l-1}}|{\emptyset_{\varepsilon_{11}}}^{\cdots A \cdots}|{\emptyset_{\varepsilon_{12}}}^{\cdots A \cdots}|{\emptyset_{\varepsilon_{13}}}^{\cdots B \cdots}|{\emptyset_{\varepsilon_{14}}}^{\cdots B \cdots}|\emptyset_{\epsilon_{l+1}}|\cdots|\emptyset_{\epsilon_{k}})
\end{eqnarray*} 
\begin{eqnarray*}
\rightarrow
T&=&(\alpha_{1}|\cdots|\alpha_{l-1}|A_{-}xA_{-}y|B_{+}zB_{+}t|\alpha_{l+1}|\cdots|\alpha_{k},\\
&&\emptyset_{\epsilon_{1}}|\cdots|\emptyset_{\epsilon_{l-1}}|{\emptyset_{\varepsilon_{21}}}^{\cdots A \cdots}|{\emptyset_{\varepsilon_{22}}}^{\cdots B \cdots}|{\emptyset_{\varepsilon_{23}}}^{\cdots B \cdots}|\emptyset_{\epsilon_{l+1}}|\cdots|\emptyset_{\epsilon_{k}})\rightarrow U,
\end{eqnarray*}
(ii)
\begin{eqnarray*}
S&=&(\alpha_{1}|\cdots|\alpha_{l-1}|A_{+}xA_{+}y|B_{+}zB_{+}t|\alpha_{l+1}|\cdots|\alpha_{k},\\
&&\emptyset_{\epsilon_{1}}|\cdots|\emptyset_{\epsilon_{l-1}}|{\emptyset_{\varepsilon_{11}}}^{\cdots A \cdots}|{\emptyset_{\varepsilon_{12}}}^{\cdots A \cdots}|{\emptyset_{\varepsilon_{13}}}^{\cdots B \cdots}|{\emptyset_{\varepsilon_{14}}}^{\cdots B \cdots}|\emptyset_{\epsilon_{l+1}}|\cdots|\emptyset_{\epsilon_{k}})
\end{eqnarray*} 
\begin{eqnarray*}
\rightarrow
T&=&(\alpha_{1}|\cdots|\alpha_{l-1}|A_{+}xA_{+}y|B_{-}zB_{-}t|\alpha_{l+1}|\cdots|\alpha_{k},\\
&&\emptyset_{\epsilon_{1}}|\cdots|\emptyset_{\epsilon_{l-1}}|{\emptyset_{\varepsilon_{31}}}^{\cdots A \cdots}|{\emptyset_{\varepsilon_{32}}}^{\cdots A \cdots}|{\emptyset_{\varepsilon_{33}}}^{\cdots B \cdots}|\emptyset_{\epsilon_{l+1}}|\cdots|\emptyset_{\epsilon_{k}})\rightarrow U.
\end{eqnarray*}
Then, the condition ($\sharp$) can also state that the
sum of the contribution of (i) to the coefficient of $U$
and the contribution of (ii) to the coefficient of $U$ is even.\\
In this case we can easily check that the condition ($\sharp$) holds since
empty words which relates $A$ and empty words which relates $B$ are 
independent.\\
$\bullet$ Consider cases (8) and (9).\par
In this cases the condition ($\sharp$) holds similarly as the case (7).\\
$\bullet$ Consider case (10).  \par
Let 
\begin{eqnarray*}
U&=&(\alpha_{1}|\cdots|\alpha_{l-1}|A_{-}xB_{-}y|A_{-}zB_{-}t|\alpha_{l+1}|\cdots|\alpha_{k},\\
&&\emptyset_{\epsilon_{1}}|\cdots|\emptyset_{\epsilon_{l-1}}|{\emptyset_{\varepsilon_{41}}}^{\cdots A B \cdots}|{\emptyset_{\varepsilon_{42}}}^{\cdots A B \cdots}|\emptyset_{\epsilon_{l+1}}|\cdots|\emptyset_{\epsilon_{k}}).
\end{eqnarray*}
It is sufficient to show that for each $(\varepsilon_{11},\varepsilon_{12}) \in \{(\pm,\pm)\}$ where double signs are arbitrary, 
the coefficient of $U$ in $d^{2}(S)$ 
is even for all $\varepsilon_{41}, \varepsilon_{42} \in \{ \pm \}$. 
Hence, for $S$ and $U$, we have to check the total
number of ways to get $U$ from $S$ is even (we denote
the condition by ($\sharp$)).  
Let us localize the problem
of the difference parts of $S$, $A_{+}$ and $B_{+}$.
Two routes (i) and (ii)  can be found to change $A_{+}$ 
(respectively $B_{+}$) into $A_{-}$  (respectively $B_{-}$) as follows:\\
(i)
\begin{eqnarray*}
S&=&(\alpha_{1}|\cdots|\alpha_{l-1}|A_{+}xB_{+}y|A_{+}zB_{+}t|\alpha_{l+1}|\cdots|\alpha_{k},\\
&&\emptyset_{\epsilon_{1}}|\cdots|\emptyset_{\epsilon_{l-1}}|{\emptyset_{\varepsilon_{11}}}^{\cdots A B \cdots}|{\emptyset_{\varepsilon_{12}}}^{\cdots A B \cdots}|\emptyset_{\epsilon_{l-1}}|\emptyset_{\epsilon_{l+1}}|\cdots|\emptyset_{\epsilon_{k}})
\end{eqnarray*} 
\begin{eqnarray*}
\rightarrow
T&=&(\alpha_{1}|\cdots|\alpha_{l-1}|A_{-}xB_{+}y|A_{-}zB_{+}t|\alpha_{l+1}|\cdots|\alpha_{k},\\
&&\emptyset_{\epsilon_{1}}|\cdots|\emptyset_{\epsilon_{l-1}}|{\emptyset_{\varepsilon_{21}}}^{\cdots A B \cdots}|\emptyset_{\epsilon_{l+1}}|\cdots|\emptyset_{\epsilon_{k}})\rightarrow U,
\end{eqnarray*}
(ii)
\begin{eqnarray*}
S&=&(\alpha_{1}|\cdots|\alpha_{l-1}|A_{+}xB_{+}y|A_{+}zB_{+}t|\alpha_{l+1}|\cdots|\alpha_{k},\\
&&\emptyset_{\epsilon_{1}}|\cdots|\emptyset_{\epsilon_{l-1}}|{\emptyset_{\varepsilon_{11}}}^{\cdots A B \cdots}|{\emptyset_{\varepsilon_{12}}}^{\cdots A B \cdots}|\emptyset_{\epsilon_{l-1}}|\emptyset_{\epsilon_{l+1}}|\cdots|\emptyset_{\epsilon_{k}})
\end{eqnarray*} 
\begin{eqnarray*}
\rightarrow
T&=&(\alpha_{1}|\cdots|\alpha_{l-1}|A_{+}xB_{-}y|A_{+}zB_{-}t|\alpha_{l+1}|\cdots|\alpha_{k}),\\
&&(\emptyset_{\epsilon_{1}}|\cdots|\emptyset_{\epsilon_{l-1}}|{\emptyset_{\varepsilon_{31}}}^{\cdots A B \cdots}|\emptyset_{\epsilon_{l+1}}|\cdots|\emptyset_{\epsilon_{k}})\rightarrow U.
\end{eqnarray*}
Then, the condition ($\sharp$) can also state that the
sum of the contribution of (i) to the coefficient of $U$
and the contribution of (ii) to the coefficient of $U$ is even.\\
Consider case $(\varepsilon_{11}, \varepsilon_{12})=(+,+)$.\par
In this case, both in route (i) and in route (ii), 
$(S,T) = 0$ for all $T$. 
Thus the condition ($\sharp$) holds.\\  
Consider case $(\varepsilon_{11}, \varepsilon_{12})=(+,-)$.\par
Consider route (i). In this route $(S,T)$ is not equal to $0$
if and only if $\varepsilon_{21} = +$. Then, for this $T$, $(T,U)$ is 
not equal to $0$ if and only if 
$(\varepsilon_{41},\varepsilon_{42}) = (+, +)$.
Consider route (ii). In this route $(S,T)$ is not equal to $0$
if and only if $\varepsilon_{31} = +$. Then, for this $T$, $(T,U)$ is 
not equal to $0$ if and only if 
$(\varepsilon_{41},\varepsilon_{42}) = (+, +)$.
Thus the condition ($\sharp$) holds.\\
Consider case $(\varepsilon_{11}, \varepsilon_{12})=(-,+)$.\par
Consider route (i). Then, $(S,T)$ is not equal to $0$
if and only if $\varepsilon_{21} = +$. Moreover for this $T$, $(T,U)$ is 
not equal to $0$ if and only if 
$(\varepsilon_{41},\varepsilon_{42}) = (+, +)$.
Consider route (ii). Then, $(S,T)$ is not equal to $0$
if and only if $\varepsilon_{31} = +$. Moreover for this $T$, $(T,U)$ is 
not equal to $0$ if and only if 
$(\varepsilon_{41},\varepsilon_{42}) = (+, +)$.
Thus the condition ($\sharp$) holds.\\ 
Consider case $(\varepsilon_{11}, \varepsilon_{12})=(-,-)$.\par
Consider route (i). Then, $(S,T)$ is not equal to $0$
if and only if $\varepsilon_{21} = -$. Moreover for this $T$, $(T,U)$ is 
not equal to $0$ if and only if 
$(\varepsilon_{41},\varepsilon_{42}) = (+, -)$ or 
$(\varepsilon_{41},\varepsilon_{42}) = (-, +)$.
Consider route (ii). Then, $(S,T)$ is not equal to $0$
if and only if $\varepsilon_{31} = -$. Moreover for this $T$, $(T,U)$ is 
not equal to $0$ if and only if 
$(\varepsilon_{41},\varepsilon_{42}) = (+, -)$.
$(\varepsilon_{41},\varepsilon_{42}) = (-, +)$.
Thus the condition ($\sharp$) holds.\\ 
$\bullet$ Consider case (11).  \par
Let 
\begin{eqnarray*}
U&=&(\alpha_{1}|\cdots|\alpha_{l-1}|\overline{A}_{-}xB_{-}y|\overline{A}_{-}zB_{-}t|\alpha_{l+1}|\cdots|\alpha_{k},\\
&&\emptyset_{\epsilon_{1}}|\cdots|\emptyset_{\epsilon_{l-1}}|{\emptyset_{\varepsilon_{41}}}^{\cdots A B \cdots}|\emptyset_{\epsilon_{l+1}}|\cdots|\emptyset_{\epsilon_{k}}).
\end{eqnarray*}
It is sufficient to show that for each $\varepsilon_{11} \in \{ \pm \}$. 
the coefficient of $U$ in $d^{2}(S)$ 
is even for all $\varepsilon_{41} \in \{ \pm \}$. 
Hence, for $S$ and $U$, we have to check the total
number of ways to get $U$ from $S$ is even (we denote
the condition by ($\sharp$)).  
Let us localize the problem
of the difference parts of $S$, $A_{+}$ and $B_{+}$.
Two routes (i) and (ii)  can be found to change $A_{+}$ 
(respectively $B_{+}$) into $A_{-}$  (respectively $B_{-}$) as follows:\\
(i)
\begin{eqnarray*}
S&=&(\alpha_{1}|\cdots|\alpha_{l-1}|\overline{A}_{+}xB_{+}y|\overline{A}_{+}zB_{+}t|\alpha_{l+1}|\cdots|\alpha_{k},\\
&&\emptyset_{\epsilon_{1}}|\cdots|\emptyset_{\epsilon_{l-1}}|{\emptyset_{\varepsilon_{11}}}^{\cdots A B \cdots}|\emptyset_{\epsilon_{l-1}}|\emptyset_{\epsilon_{l+1}}|\cdots|\emptyset_{\epsilon_{k}})
\end{eqnarray*} 
\begin{eqnarray*}
\rightarrow
T&=&(\alpha_{1}|\cdots|\alpha_{l-1}|\overline{A}_{-}xB_{+}y|\overline{A}_{-}zB_{+}t|\alpha_{l+1}|\cdots|\alpha_{k},\\
&&\emptyset_{\epsilon_{1}}|\cdots|\emptyset_{\epsilon_{l-1}}|{\emptyset_{\varepsilon_{21}}}^{\cdots A B \cdots}|{\emptyset_{\varepsilon_{22}}}^{\cdots A B\cdots}|\emptyset_{\epsilon_{l+1}}|\cdots|\emptyset_{\epsilon_{k}})\rightarrow U,
\end{eqnarray*}
(ii)
\begin{eqnarray*}
S&=&(\alpha_{1}|\cdots|\alpha_{l-1}|\overline{A}_{+}xB_{+}y|\overline{A}_{+}zB_{+}t|\alpha_{l+1}|\cdots|\alpha_{k},\\
&&\emptyset_{\epsilon_{1}}|\cdots|\emptyset_{\epsilon_{l-1}}|{\emptyset_{\varepsilon_{11}}}^{\cdots A B \cdots}|\emptyset_{\epsilon_{l-1}}|\emptyset_{\epsilon_{l+1}}|\cdots|\emptyset_{\epsilon_{k}})
\end{eqnarray*} 
\begin{eqnarray*}
\rightarrow
T&=&(\alpha_{1}|\cdots|\alpha_{l-1}|\overline{A}_{+}xB_{-}y|\overline{A}_{+}zB_{-}t|\alpha_{l+1}|\cdots|\alpha_{k},\\
&&\emptyset_{\epsilon_{1}}|\cdots|\emptyset_{\epsilon_{l-1}}|{\emptyset_{\varepsilon_{31}}}^{\cdots A B \cdots}|{\emptyset_{\varepsilon_{32}}}^{\cdots A B\cdots}|\emptyset_{\epsilon_{l+1}}|\cdots|\emptyset_{\epsilon_{k}})\rightarrow U.
\end{eqnarray*}
In this case we can chose $\varepsilon_{21}$, $\varepsilon_{22}$,
$\varepsilon_{31}$ ,$\varepsilon_{32} \in \{ \pm \}$ so that 
$(\varepsilon_{21},\varepsilon_{22}) = (\varepsilon_{31},\varepsilon_{32})$
and $(S,T)$ is not equal to $0$. Moreover $T$s in route (i) and 
in route (ii) have same form. Thus the condition ($\sharp$) holds.\\
$\bullet$ Consider case (12).  \par
Let 
\begin{eqnarray*}
U&=&(\alpha_{1}|\cdots|\alpha_{l-1}|\overline{A}_{-}x\overline{B}_{-}y|\overline{A}_{-}z\overline{B}_{-}t|\alpha_{l+1}|\cdots|\alpha_{k},\\
&&\emptyset_{\epsilon_{1}}|\cdots|\emptyset_{\epsilon_{l-1}}|{\emptyset_{\varepsilon_{41}}}^{\cdots A B \cdots}|{\emptyset_{\varepsilon_{42}}}^{\cdots A B \cdots}|\emptyset_{\epsilon_{l+1}}|\cdots|\emptyset_{\epsilon_{k}}).
\end{eqnarray*}
Then, the condition ($\sharp$) can also state that the
sum of the contribution of (i) to the coefficient of $U$
and the contribution of (ii) to the coefficient of $U$ is even.\\
It is sufficient to show that for each $(\varepsilon_{11}, \varepsilon_{12}) 
\in \{ \pm, \pm \}$ where double signs are arbitrary, 
the coefficient of $U$ in $d^{2}(S)$ 
is even for all $\varepsilon_{41}, \varepsilon_{42} \in \{ \pm \}$. 
Hence, for $S$ and $U$, we have to check the total
number of ways to get $U$ from $S$ is even (we denote
the condition by ($\sharp$)).  
Let us localize the problem
of the difference parts of $S$, $A_{+}$ and $B_{+}$.
Two routes (i) and (ii)  can be found to change $A_{+}$ 
(respectively $B_{+}$) into $A_{-}$  (respectively $B_{-}$) as follows:\\
(i)
\begin{eqnarray*}
S&=&(\alpha_{1}|\cdots|\alpha_{l-1}|\overline{A}_{+}x\overline{B}_{+}y|\overline{A}_{+}z\overline{B}_{+}t|\alpha_{l+1}|\cdots|\alpha_{k},\\
&&\emptyset_{\epsilon_{1}}|\cdots|\emptyset_{\epsilon_{l-1}}|{\emptyset_{\varepsilon_{11}}}^{\cdots A B \cdots}|{\emptyset_{\varepsilon_{12}}}^{\cdots A B \cdots}|\emptyset_{\epsilon_{l-1}}|\emptyset_{\epsilon_{l+1}}|\cdots|\emptyset_{\epsilon_{k}})
\end{eqnarray*} 
\begin{eqnarray*}
\rightarrow
T&=&(\alpha_{1}|\cdots|\alpha_{l-1}|\overline{A}_{-}x\overline{B}_{+}y|\overline{A}_{-}z\overline{B}_{+}t|\alpha_{l+1}|\cdots|\alpha_{k},\\
&&(\emptyset_{\epsilon_{1}}|\cdots|\emptyset_{\epsilon_{l-1}}|{\emptyset_{\varepsilon_{21}}}^{\cdots A B \cdots}|\emptyset_{\epsilon_{l+1}}|\cdots|\emptyset_{\epsilon_{k}})\rightarrow U,
\end{eqnarray*}
(ii)
\begin{eqnarray*}
S&=&(\alpha_{1}|\cdots|\alpha_{l-1}|\overline{A}_{+}x\overline{B}_{+}y|\overline{A}_{+}z\overline{B}_{+}t|\alpha_{l+1}|\cdots|\alpha_{k},\\
&&\emptyset_{\epsilon_{1}}|\cdots|\emptyset_{\epsilon_{l-1}}|{\emptyset_{\varepsilon_{11}}}^{\cdots A B \cdots}|{\emptyset_{\varepsilon_{12}}}^{\cdots A B \cdots}|\emptyset_{\epsilon_{l-1}}|\emptyset_{\epsilon_{l+1}}|\cdots|\emptyset_{\epsilon_{k}})
\end{eqnarray*} 
\begin{eqnarray*}
\rightarrow
T&=&(\alpha_{1}|\cdots|\alpha_{l-1}|\overline{A}_{+}x\overline{B}_{-}y|\overline{A}_{+}z\overline{B}_{-}t|\alpha_{l+1}|\cdots|\alpha_{k},\\
&&\emptyset_{\epsilon_{1}}|\cdots|\emptyset_{\epsilon_{l-1}}|{\emptyset_{\varepsilon_{31}}}^{\cdots A B \cdots}|\emptyset_{\epsilon_{l+1}}|\cdots|\emptyset_{\epsilon_{k}})\rightarrow U.
\end{eqnarray*}
Then, the condition ($\sharp$) can also state that the
sum of the contribution of (i) to the coefficient of $U$
and the contribution of (ii) to the coefficient of $U$ is even.\\
This case is completely same as the case (ii).\\
$\bullet$The cases (13) - (23)\par
We can easily check the condition ($\sharp$) by the definition of $d$.\\
$\bullet$The cases (24) - (26)\par
In this case we can prove the condition ($\sharp$) holds 
same as the case (7).\par
Now, we have proved the theorem.
\end{proof}

\begin{definition}
We denote the mapping $d$ modulo $2$ $: C^{i, j}(P; {\mathbb{Z}}_{2}) \to C^{i+1, j}(P; {\mathbb{Z}}_{2})$ by $d_{2}^{i}$ for $i$ and $j$.  The 
Khovanov homology group $KH^{i, j}(P)$ for a pseudolink $P$ is defined as 
\begin{equation}
KH^{i, j}(P) := {\rm Ker} d_{2}^{i} \big{/} {\rm Im} d_{2}^{i-1}.  
\end{equation}
\end{definition}

\begin{remark}
$KH^{i, j}(P)$ is independent of the order of that the letters are removed from $P$ because the incidence number $(S : T)$ is always either $0$ or $1$ modulo $2$ for enhanced states $S$ and $T$.  
\end{remark}

\section{Invariance under $S_{1}$-homotopy moves}

\begin{theorem}[Manturov]
$KH^{i, j}(P)$ are $S_{1}$-homotopy invariants for pseudolinks.  
\end{theorem}

\begin{remark}
Manturov proved this result using virtual knot theory \cite{manturov2}.  
\end{remark}

\begin{proof}
From the construction of $KH^{i,j}(P)$, it is evident that $KH^{i,j}(P)$ does not depend on an arbitrary isomorphism of $P$.  Then, $KH^{i,j}$ is invariant under isomorphisms.  It remains to be proved that if a nanophrase $P$ is obtained from a nanophrase $P'$ by a homotopy move, then $KH^{i,j}(P')$ $\simeq$ $KH^{i,j}(P)$.  
The following discussion is similar to those in \cite[Subsection 5.6]{viro} and \cite[Section 2 and 3]{ito3}.  

(I) Consider the first homotopy move $xAAy$ $\rightarrow$ $xy$ and its inverse move where $|A|$ $=$ $1$.  
For $P'$ and $P$, $S_{+}(\epsilon, \eta)$ denotes the state $u|{\emptyset_{\epsilon}}^{Aw}|{\emptyset_{\eta}}^{A}|v$ of $P'$ with mark($A$) $=$ $1$ and $S_{-}(\epsilon)$ denotes the state $u|{\emptyset_{\epsilon}}^{Aw}|v$ of $P'$ with mark($A$) $=$ $-1$, where $\epsilon$, $\eta$ $\in \{+, -\}$.  
The subcomplex $C'$ of ${C}(P')$ is defined by $C'$ $:=$ ${C}( S_{+}(+, +)$, $S_{+}(+, -)$ $-$ $S_{+}(-, +))$.  

First, the retraction \[ \rho : {C}(P') \to {C}(S_{+}(+, +),~S_{+}(+, -) -  S_{+}(-, +))\] is defined by the formulas 
\begin{align*}
S_{+}(+, +) &\mapsto S_{+}(+, +), \\
S_{+}(-, +) &\mapsto S_{+}(-, +) - S_{+}(+, -), \\
{\rm{otherwise}} &\mapsto 0.  
\end{align*}

Second, the isomorphism \[ {C}(S_{+}(+, +), S_{+}(+, -) -  S_{+}(-, +)) \to {C}(P) = {C}(u|{\emptyset_{+}}^{w}|v, u|{\emptyset_{-}}^{w}|v)\] is defined by the formulas 
\begin{align*}
S_{+}(+, +) &\mapsto u|{\emptyset_{+}}^{w}|v, \\
S_{+}(+, -) - S_{+}(-, +) &\mapsto u|{\emptyset_{-}}^{w}|v.  
\end{align*}

Third, consider the following composition of this isomorphism with $\rho:$ 
\begin{equation*}
C(P') \stackrel{\rho}{\to} C' \stackrel{\mathrm{isom}}{\to} C(P).  
\end{equation*}

The map $h : {C}(P') \to {C}(P')$ such that $d \circ h$ $+$ $h \circ d$ $=$ $\operatorname{id} - \operatorname{in} \circ \rho$ is defined by the formulas 
\begin{align*}
S_{-}(+) &\mapsto S_{+}(+, -), \\
S_{-}(-) &\mapsto S_{+}(-, -), \\
{\rm{otherwise}} &\mapsto 0.  
\end{align*}

(I\!I) Consider the second homotopy move $P'$ $=$ $xAByBAz$ $\rightarrow$ $xyz$ $=$ $P$ and its inverse move where $(|A|, |B|)$ $=$ $(1, -1)$.  It is necessary to consider two distinct cases (I\!I-1), (I\!I-2) as follows.  

(I\!I--1) Consider case where the state of $P'$ with (mark($A$), mark($B$)) $=$ ($1$, $1$) is represented as $u|{\emptyset_{\epsilon}}^{ABw}|v$.  

$S_{+-}(\epsilon, \eta)$ denotes the state $u|{\emptyset_{\epsilon}}^{Aw}|{\emptyset_{\eta}}^{A B}|v$ of $P'$ with (mark($A$), mark($B$)) $=$ ($1$, $-1$), $S_{-+}(\epsilon, \eta)$ denotes the state $u|{\emptyset_{\epsilon}}^{ABw}|{\emptyset_{\eta}}^{ABt}|v$ of $P'$ with (mark($A$), mark($B$)) $=$ ($-1$, $1$), $S_{++}(\epsilon)$ denotes the state $u|{\emptyset_{\epsilon}}^{ABw}|v$ of $P'$ with (mark($A$), mark($B$)) $=$ ($1$, $1$), and $S_{--}(\epsilon)$ denotes the state $u|{\emptyset_{\epsilon}}^{ABw}|v$ of $P'$ with (mark($A$), mark($B$)) $=$ ($-1$, $-1$), where $\epsilon$, $\eta$ $\in \{+, -\}$.  
The subcomplex $C'$ of ${C}(P')$ is defined by $C'$ $:=$ ${C}\big( S_{-+}(+, +),$ $S_{-+}(+, -) + S_{+-}(+, -),$ $S_{-+}(-, +) + S_{+-}(-, -),$ $S_{-+}(-, -) + S_{+-}(-, -) \big)$.  

First, the retraction $\rho :$ ${C}(P')$ $\to$ $C'$ is defined by the formulas 
\begin{align*}
S_{-+}(+, +) &\mapsto S_{-+}(+, +), \\
S_{-+}(+, -) &\mapsto S_{-+}(+, -) + S_{+-}(+), \\
S_{-+}(-, +) &\mapsto S_{-+}(-, +) + S_{+-}(+), \\
S_{-+}(-, -) &\mapsto S_{-+}(-, -) + S_{+-}(-), \\
S_{+-}(+, +) &\mapsto S_{-+}(+, +), \\
S_{+-}(-, +) &\mapsto S_{-+}(+, -) + S_{-+}(-, +), \\
{\rm{otherwise}} &\mapsto 0.  
\end{align*}  

Second, the isomorphism \[ C' \to {C}(P) = {C}(u|{\emptyset_{\epsilon}}^{w}|{\emptyset_{\eta}}^{t}|v)\] is defined by the formulas 
\begin{align*}
S_{-+}(+, +) &\mapsto u|{\emptyset_{+}}^{w}|{\emptyset_{+}}^{t}|v, \\
S_{-+}(+, -) + S_{+-}(+) &\mapsto u|{\emptyset_{+}}^{w}|{\emptyset_{-}}^{t}|v, \\
S_{-+}(-, +) + S_{+-}(+) &\mapsto u|{\emptyset_{-}}^{w}|{\emptyset_{+}}^{t}|v, \\
S_{-+}(-, -) + S_{+-}(-) &\mapsto u|{\emptyset_{-}}^{w}|{\emptyset_{-}}^{t}|v.  
\end{align*}

Third, consider the following composition of this isomorphism with  $\rho:$ 
\begin{equation*}
C(P') \stackrel{\rho}{\to} C' \stackrel{\mathrm{isom}}{\to} C(P).  
\end{equation*}

The map $h:$ ${C}(P') \to {C}(P')$ such that $d \circ h$ $+$ $h \circ d$ $=$ $\operatorname{id} - \operatorname{in} \circ \rho$, is defined by the formulas 

\begin{align*}
S_{--}(\epsilon) &\mapsto S_{+-}(\epsilon, -), \\
S_{+-}(\epsilon, +) &\mapsto S_{++}(\epsilon), \\
{\rm{otherwise}} &\mapsto 0.  
\end{align*}

(I\!I--2) Consider the case where the state of $P'$ with (mark($A$), mark($B$)) $=$ ($1$, $1$) is represented as $u|{\emptyset_{\epsilon}}^{Aw}|{\emptyset_{\eta}}^{A Bt}|v$.  

$S_{+-}(\epsilon, \zeta, \eta)$ denotes the state $u|{\emptyset_{\epsilon}}^{Aw}|{\emptyset_{\zeta}}^{A B}|{\emptyset_{\eta}}^{Bt}|v$ of $P'$ with (mark($A$), mark($B$)) $=$ ($1$, $-1$), $S_{-+}(\epsilon)$ denotes the state $u|{\emptyset_{\epsilon}}^{A Bwt'}|v$ of $P'$ with (mark($A$), mark($B$)) $=$ ($-1$, $1$), $S_{++}(\epsilon, \eta)$ denotes the state $u|{\emptyset_{\epsilon}}^{Aw}|{\emptyset_{\eta}}^{A Bt}|v$ of $P'$ with (mark($A$), mark($B$)) $=$ ($1$, $1$), and $S_{--}(\epsilon, \eta)$ denotes the state $u|{\emptyset_{\epsilon}}^{A Bw}|{\emptyset_{\eta}}^{Bt}|v$ of $P'$ with (mark($A$), mark($B$)) $=$ ($-1$, $-1$), where $\epsilon$, $\eta$ $\in \{+, -\}$ and the word $t'$ is obtained by deleting all letters from $t$ that appear in $w$.  The subcomplex $C'$ of ${C}(P')$ is defined by $C'$ $:=$ ${C}\big($$S_{-+}(+)$ $+$ $S_{+-}(+, -,  +)$, $S_{-+}(-)$ $+$ $S_{+-}(+, -, -)$ $+$ $S_{+-}(-, -, +)$$\big)$.  

First, the retraction $\rho:$ ${C}(P')$ $\to$ $C'$ is defined by the formulas 
\begin{align*}
S_{-+}(+) &\mapsto S_{-+}(+) + S_{+-}(+, -, +), \\
S_{-+}(-) &\mapsto S_{-+}(-) + S_{+-}(+, -, -) + S_{+-}(-, -, +), \\
S_{+-}(+, +, -) &\mapsto S_{-+}(+) + S_{+-}(+, -, +), \\
S_{+-}(-, +, +) &\mapsto S_{-+}(+) + S_{+-}(+, -, +), \\
S_{+-}(-, +, -) &\mapsto S_{-+}(-) + S_{+-}(+, -, -) + S_{+-}(-, -, +), \\
{\rm{otherwise}} &\mapsto 0.  
\end{align*}  

Second, the isomorphism \[ C' \to {C}(P) = {C}(u|{\emptyset_{\epsilon}}^{wt'}|v)\] is defined by the formulas 
\begin{align*}
S_{-+}(+) + S_{+-}(+, -, +) &\mapsto u|{\emptyset_{+}}^{wt'}|v, \\
S_{-+}(-) + S_{+-}(+, -, -) + S_{+-}(-, -, +) &\mapsto u|{\emptyset_{-}}^{wt'}|v.  
\end{align*}

Third, consider the following composition of this isomorphism with  $\rho:$ 
\begin{equation*}
C(P') \stackrel{\rho}{\to} C' \stackrel{\mathrm{isom}}{\to} C(P).  
\end{equation*}

The map $h:$ ${C}(P') \to {C}(P')$ such that $d \circ h$ $+$ $h \circ d$ $=$ $\operatorname{id} - \operatorname{in} \circ \rho$, is defined by the formulas 

\begin{align*}
S_{--}(\epsilon, \eta) &\mapsto S_{+-}(\epsilon, -, \eta), \\
S_{+-}(\epsilon, +, \eta) &\mapsto S_{++}(\epsilon, \eta), \\
{\rm{otherwise}} &\mapsto 0.  
\end{align*}

By using (I\!I-1) and (I\!I-2), we proved that $KH^{i, j}(xAByBAz)$ $\simeq$ $KH^{i, j}(xyz)$ if $(|A|, |B|)$ $=$ $(1, -1)$.  In addition, (I\!I-1) and (I\!I-2) prove that $KH^{i, j}(xAByABz)$ $\simeq$ $KH^{i, j}(xyz)$ if $(|A|, |B|)$ $=$ $(-1, 1)$.  Moreover, by exchanging $A$, $B$ in the proofs above, (I\!I-1) and (I\!I-2) prove that $KH^{i, j}(xAByBAz)$ $\simeq$ $KH^{i, j}(xyz)$ if $(|A|, |B|)$ $=$ $(-1, 1)$ and $KH^{i, j}(xAByABz)$ $\simeq$ $KH^{i, j}(xyz)$ if $(|A|, |B|)$ $=$ $(1, -1)$.  

Here, consider 
\begin{align*}
xAAy &\stackrel{\text{H1}}{\sim} xABBAy \quad \text{with}\ |A| = -1, |B| = 1\\  
&\stackrel{\text{H2}}{\sim} xy.  
\end{align*}
We have already shown the invariance of $KH^{i,j}$ under the above moves and that $KH^{i,j}$ is preserved under the first homotopy move $xAAy$ $\rightarrow$ $xy$ with $|A|$ $=$ $-1$ and its inverse move.  

(I\!I\!I) Consider the third homotopy move 
\begin{equation*}
P' = xAByACzBCt \rightarrow xBAyCAzCBt = P
\end{equation*}
and its inverse move where $(|A|, |B|, |C|)$ $=$ $(-1, -1, -1)$.  For the letters $A$, $B$, and $C$, we define $w_{ABC}$, $w_{AB}$, $w_{AC}$, $w_{BC}$, $w_{A}$, $w_{B}$, and $w_{C}$ in the following.  
Let $w_{ABC}$ be a word containing $A$, $B$, and $C$.  Let $(X, Y, Z)$ $=$ $\{(A, B, C)$, $(A, C, B)$, $(B, C, A)\}$.  $w_{XY}$ denotes a word containing $X$ and $Y$ but not $Z$, and $w_{Z}$ denotes a word containing $Z$ but not $X$ and $Y$.  

(I\!I\!I--1) Consider the case where the state of $P'$ with (mark($A$), mark($B$), mark($C$ )) $=$ ($1$, $1$, $1$) is represented as $u|{\emptyset_{\epsilon}}^{w_{A B C}}|v$.  

%+++
$S_{+++}(\epsilon)$ denotes the state $u|{\emptyset_{\epsilon}}^{w_{A B C}}|v$ of $P'$ 
with (mark($A$), mark($B$), mark($C$)) $=$ ($1$, $1$, $1$), 
%-++
$S_{-++}(\epsilon, \eta)$ denotes the state $u|{\emptyset_{\epsilon}}^{w_{A B C}}|{\emptyset_{\eta}}^{A B}|v$ of $P'$ 
with (mark($A$), mark($B$), mark($C$)) $=$ ($-1$, $1$, $1$), 
%+-+
$S_{+-+}(\epsilon, \eta)$ denotes the state $u|{\emptyset_{\epsilon}}^{w_{A B C}}|{\emptyset_{\eta}}^{A B C}|v$ of $P'$ 
with (mark($A$), mark($B$), mark($C$)) $=$ ($1$, $-1$, $1$), 
%--+
$S_{--+}(\epsilon)$ denotes the state $u|{\emptyset_{\epsilon}}^{w_{A B C}}|v$ of $P'$ 
with (mark($A$), mark($B$), mark($C$)) $=$ ($-1$, $-1$, $1$), 
%++-
$S_{++-}(\epsilon, \eta)$ denotes the state $u|{\emptyset_{\epsilon}}^{w_{A B C}}|{\emptyset_{\eta}}^{w_{B C}}|v$ of $P'$ 
with (mark($A$), mark($B$), mark($C$)) $=$ ($1$, $1$, $-1$), 
%-+-
$S_{-+-}(\epsilon, \zeta, \eta)$ denotes the state $u|{\emptyset_{\epsilon}}^{w_{A C}}|{\emptyset_{\zeta}}^{w_{A B}}|{\emptyset_{\eta}}^{w_{B C}}|v$ of $P'$ 
with (mark($A$), mark($B$), mark($C$)) $=$ ($-1$, $1$, $-1$), 
%+--
$S_{+--}(\epsilon)$ denotes the state $u|{\emptyset_{\epsilon}}^{w_{A B C}}|v$ of $P'$ 
with (mark($A$), mark($B$), mark($C$)) $=$ ($1$, $-1$, $-1$), 
%---
and $S_{---}(\epsilon, \eta)$ denotes the state $u|{\emptyset_{\epsilon}}^{w_{A C}}|{\emptyset_{\eta}}^{w_{A B C}}|v$ of $P'$ with (mark($A$), mark($B$), mark($C$)) $=$ ($-1$, $-1$, $-1$).  

The subcomplex $C'$ of ${C}(P')$ is defined by $C'$ $:=$ ${C}\big($ $S_{-++}(+, +)$, $S_{-++}(+, -)$ $+$ $S_{+-+}(+, -)$, $S_{-++}(-, +)$ $+$ $S_{+-+}(+, -)$, $S_{-++}(-, -)$ $+$ $S_{+-+}(-, -)$, $S_{**-}$ $\big)$, where $S_{**-}$ denotes every state with mark($C$) $=$ $-1$.  

%T+-+
$T_{+-+}(\epsilon, \eta)$ denotes the state $u|{\emptyset_{\epsilon}}^{w_{A B C}}|{\emptyset_{\eta}}^{w_{C}}|v$ of $P$ 
with (mark($A$), mark($B$), mark( $C$)) $=$ ($1$, $-1$, $1$), 
%T-++,-
$T_{-++}(\epsilon, \zeta, \eta, -)$ denotes the state $u|{\emptyset_{\epsilon}}^{w_{A}}|{\emptyset_{\zeta}}^{w_{B}}|{\emptyset_{\eta}}^{w_{C}}|{\emptyset_{-}}^{A B C}|v$ of $P$ 
with (mark($A$), mark($B$), mark($C$)) $=$ ($-1$, $1$, $1$), 
%++-
$T_{++-}(\epsilon, \eta)$ denotes the state $u|{\emptyset_{\epsilon}}^{w_{A B C}}|{\emptyset_{\eta}}^{w_{B}}|v$ of $P$ 
with (mark($A$), mark($B$), mark($C$)) $=$ ($1$, $1$, $-1$), 
%-+-
$T_{-+-}(\epsilon, \zeta, \eta)$ denotes the state $u|{\emptyset_{\epsilon}}^{w_{A}}|{\emptyset_{\zeta}}^{w_{A B C}}|{\emptyset_{\eta}}^{w_{B}}|v$ of $P$ with (mark($A$), mark($B$), mark($C$)) $=$ ($-1$, $1$, $-1$), 
%+--
$T_{+--}(\epsilon)$ denotes the state $u|{\emptyset_{\epsilon}}^{w_{A B C}}|v$ of $P$ 
with (mark($A$), mark($B$), mark($C$)) $=$ ($1$, $-1$, $-1$), 
%---
$T_{---}(\epsilon, \eta)$ denotes the state $u|{\emptyset_{\epsilon}}^{w_{A}}|{\emptyset_{\eta}}^{w_{A B C}}|v$ of $P$ 
with (mark($A$), mark($B$), mark($C$)) $=$ ($-1$, $-1$, $-1$), 
%**-
and $T_{**-}$ denotes every state of $P$ with mark($C$) $=$ $-1$.  

The subcomplex $C$ of ${C}(P)$ is defined by $C$ $:=$ ${C}\big($$T_{+-+}(+, \eta)$ $+$ $T_{-++}(+, +, \eta,$ $-)$, $T_{+-+}(-, \eta)$ $+$ $T_{-++}(+, -, \eta, -)$ $+$ $T_{-++}(-, +, \eta, -)$, $T_{**-}$$\big)$.  

First, the retraction $\rho:$ ${C}(P')$ $\to$ $C'$ is defined by the formulas 
\begin{align*}
S_{-++}(+, +) &\mapsto S_{-++}(+, +), \\
S_{-++}(+, -) &\mapsto S_{-++}(+, -) + S_{+-+}(+, -), \\
S_{-++}(-, +) &\mapsto S_{-++}(-, +) + S_{+-+}(+, -), \\
S_{-++}(-, -) &\mapsto S_{-++}(-, -) + S_{+-+}(-, -), \\
S_{**-} &\mapsto S_{**-}, \\
S_{+-+}(+, +) &\mapsto S_{-++}(+, +) + S_{++-}(+, +), \\
S_{+-+}(-, +) &\mapsto S_{-++}(+, -) + S_{-++}(-, +) + S_{++-}(+, -) + S_{++-}(-, +), \\
S_{--+}(\epsilon) &\mapsto S_{+--}(\epsilon), \\
{\rm{otherwise}} &\mapsto 0.  
\end{align*}  

Second, consider the following composition of the following isomorphism with $\rho$
\begin{equation}\label{third-rei-eq}
C(P') \stackrel{\rho}{\to} C' \stackrel{\mathrm{isom}}{\to} C \stackrel{i}{\to} C(P).  
\end{equation}
The isomorphism $C'$ $\to$ ${C}$ is defined by the formulas
\begin{align*}
S_{-++}(+, +) &\mapsto T_{+-+}(+, +) + T_{-++}(+, +, +, -), \\
S_{-++}(+, -) + S_{+-+}(+, -) &\mapsto T_{+-+}(+, -) + T_{-++}(+, +, -, -),\\
S_{-++}(-, +) + S_{+-+}(+, -)  &\mapsto T_{+-+}(-, +) + T_{-++}(+, -, +, -)\\&\qquad\qquad\qquad\!\quad + T_{-++}(-, +, +, -), \\
S_{-++}(-, -) + S_{+-+}(-, -)  &\mapsto T_{+-+}(-, -) + T_{-++}(+, -, -, -)\\&\qquad\qquad\qquad\!\quad + T_{-++}(-, +, -, -), \\
S_{++-}(\epsilon, \eta) &\mapsto T_{++-}(\epsilon, \eta), \\
S_{-+-}(\epsilon, \zeta, \eta) &\mapsto T_{+-+}(\epsilon, \zeta, \eta), \\
S_{+--}(\epsilon) &\mapsto T_{+--}(\epsilon), \\
S_{---}(\epsilon, \eta) &\mapsto T_{---}(\epsilon, \eta).  
\end{align*}

Third, the map $h:$ ${C}(P') \to {C}(P')$ such that $d \circ h$ $+$ $h \circ d$ $=$ $\operatorname{id} - \operatorname{in} \circ \rho$ is defined by the formulas 

\begin{align*}
S_{--+}(\epsilon) &\mapsto S_{+-+}(\epsilon, -), \\
S_{+-+}(\epsilon, +) &\mapsto S_{+++}(\epsilon), \\
{\rm{otherwise}} &\mapsto 0.  
\end{align*}

(I\!I\!I--2) Consider the case where the state of $P'$ with (mark($A$), mark($B$), mark($C$ )) $=$ ($1$, $1$, $1$) is represented as $u|{\emptyset_{\epsilon}}^{w_{A C}}|{\emptyset_{\eta}}^{w_{A B C}}|v$.  

%+++
$S_{+++}(\epsilon, \eta)$ denotes the state $u|{\emptyset_{\epsilon}}^{w_{A C}}|{\emptyset_{\eta}}^{w_{A B C}}|v$ of $P'$ 
with (mark($A$), mark($B$), mark($C$)) $=$ ($1$, $1$, $1$), 
%-++
$S_{-++}(\epsilon)$ denotes the state $u|{\emptyset_{\epsilon}}^{w_{A B C}}|v$ of $P'$ 
with (mark($A$), mark($B$), mark($C$)) $=$ ($-1$, $1$, $1$), 
%+-+
$S_{+-+}(\epsilon, \zeta, \eta)$ denotes the state $u|{\emptyset_{\epsilon}}^{w_{A C}}|{\emptyset_{\zeta}}^{A B C}|{\emptyset_{\eta}}^{w_{B}}$ $|v$ of $P'$ with (mark($A$), mark($B$), mark($C$)) $=$ ($1$, $-1$, $1$), 
%--+
$S_{--+}(\epsilon, \eta)$ denotes the state $u|{\emptyset_{\epsilon}}^{w_{A B C}}|{\emptyset_{\eta}}^{w_{B}}|v$ of $P'$ 
with (mark($A$), mark($B$), mark($C$)) $=$ ($-1$, $-1$, $1$), 
%++-
$S_{++-}(\epsilon)$ denotes the state $u|{\emptyset_{\epsilon}}^{w_{A B C}}|v$ of $P'$ 
with (mark($A$), mark($B$), mark($C$)) $=$ ($1$, $1$, $-1$), 
%-+-
$S_{-+-}(\epsilon, \eta)$ denotes the state $u|{\emptyset_{\epsilon}}^{w_{A B C}}|{\emptyset_{\eta}}^{w_{A C}}|v$ of $P'$ 
with (mark($A$), mark($B$), mark($C$)) $=$ ($-1$, $1$, $-1$), 
%+--
$S_{+--}(\epsilon, \eta)$ denotes the state $u|{\emptyset_{\epsilon}}^{w_{A B C}}|{\emptyset_{\eta}}^{w_{B}}|v$ of $P'$ 
with (mark($A$), mark($B$), mark($C$)) $=$ ($1$, $-1$, $-1$), 
%---
and $S_{---}(\epsilon, \zeta, \eta)$ denotes the state $u|{\emptyset_{\epsilon}}^{w_{A B C}}|{\emptyset_{\zeta}}^{w_{A C}}|{\emptyset_{\eta}}^{w_{B}}|v$ of $P'$ 
with (mark($A$), mark($B$), mark($C$)) $=$ ($-1$, $-1$, $-1$).

The subcomplex $C'$ of ${C}(P')$ is defined by $C'$ $:=$ ${C}\big($$S_{-++}(+)$ $+$ $S_{+-+}(+, -, +)$, $S_{-++}(-)$ $+$ $S_{+-+}(+, -, -)$ $+$ $S_{+-+}(-, -, +)$, $S_{**-}$$\big)$, where $S_{**-}$ denotes every state with mark($C$) $=$ $-1$.  

%T+-+
$T_{+-+}(\epsilon)$ denotes the state $u|{\emptyset_{\epsilon}}^{w_{A B C}}|v$ of $P$ 
with (mark($A$), mark($B$), mark( $C$)) $=$ ($1$, $-1$, $1$), 
%T-++,-
$T_{-++}(\epsilon, \zeta, -)$ denotes the state $u|{\emptyset_{\epsilon}}^{w_{B C}}|{\emptyset_{\zeta}}^{w_{A}}|{\emptyset_{-}}^{A B C}|v$ of $P$ 
with (mark($A$), mark($B$), mark($C$)) $=$ ($-1$, $1$, $1$), 
%++-
$T_{++-}(\epsilon)$ denotes the state $u|{\emptyset_{\epsilon}}^{w_{A B C}}|v$ of $P$ 
with (mark($A$), mark($B$), mark($C$)) $=$ ($1$, $1$, $-1$), 
%-+-
$T_{-+-}(\epsilon, \eta)$ denotes the state $u|{\emptyset_{\epsilon}}^{w_{B C}}|{\emptyset_{\eta}}^{w_{A}}|v$ of $P$ 
with (mark($A$), mark($B$), mark($C$)) $=$ ($-1$, $1$, $-1$), 
%+--
$T_{+--}(\epsilon, \eta)$ denotes the state $u|{\emptyset_{\epsilon}}^{w_{A B C}}|{\emptyset_{\eta}}^{w_{B C}}|v$ of $P$ 
with (mark($A$), mark($B$), mark($C$)) $=$ ($1$, $-1$, $-1$), 
%---
$T_{---}(\epsilon, \zeta, \eta)$ denotes the state $u|{\emptyset_{\epsilon}}^{w_{A B C}}|{\emptyset_{\zeta}}^{w_{A}}|{\emptyset_{\eta}}^{w_{B C}}|v$ of $P$ 
with (mark($A$), mark($B$), mark($C$)) $=$ ($-1$, $-1$, $-1$), 
%**-
and $T_{**-}$ denotes every state of $P$ with mark($C$) $=$ $-1$.  

The subcomplex $C$ of ${C}(P)$ is defined by $C$ $:=$ ${C}\big($$T_{+-+}(+)$ $+$ $T_{-++}(+, +, -)$, $T_{+-+}(-)$ $+$ $T_{-++}(+, -, -)$ $+$ $T_{-++}(-, +, -)$, $T_{**-}$$\big)$.

First, the retraction $\rho:$ ${C}(P')$ $\to$ $C'$ is defined by the formulas 
\begin{align*}
S_{-++}(+) &\mapsto S_{-++}(+) + S_{+-+}(+, -, +), \\
S_{-++}(-) &\mapsto S_{-++}(-) + S_{+-+}(+, -, -) + S_{+-+}(-, -, +), \\
S_{**-} &\mapsto S_{**-}, \\
S_{+-+}(+, +, -) &\mapsto S_{-++}(+) + S_{+-+}(+, -, +) + S_{++-}(+), \\
S_{+-+}(-, +, +) &\mapsto S_{-++}(+) + S_{+-+}(+, -, +) + S_{++-}(+), \\
S_{+-+}(-, +, -) &\mapsto S_{-++}(-) + S_{+-+}(+, -, -) + S_{+-+}(-, -, -) + S_{++-}(-), \\
S_{--+}(\epsilon, \eta) &\mapsto S_{+--}(\epsilon, \eta), \\
{\rm{otherwise}} &\mapsto 0.  
\end{align*}

Second, consider the following composition (\ref{third-rei-eq}) of the following isomorphism with $\rho$.  
The isomorphism $C'$ $\to$ ${C}$ is defined by the formulas
\begin{align*}
S_{-++}(+) + S_{+-+}(+, -, +) &\mapsto T_{+-+}(+) + T_{-++}(+, +, -), \\
S_{-++}(-) + S_{+-+}(+, -, -) + S_{+-+}(-, -, +) &\mapsto T_{+-+}(-) + T_{-++}(+, -, -)\\
& \qquad\qquad\quad \ \, + T_{-++}(-, +, -), \\
S_{++-}(\epsilon) &\mapsto T_{++-}(\epsilon),  \\
S_{-+-}(\epsilon, \eta) &\mapsto T_{-+-}(\epsilon, \eta), \\
S_{+--}(\epsilon, \eta) &\mapsto T_{+--}(\epsilon, \eta), \\
S_{---}(\epsilon, \zeta, \eta) &\mapsto T_{---}(\epsilon, \zeta, \eta).  
\end{align*}

Third, the map $h:$ ${C}(P') \to {C}(P')$ such that $d \circ h$ $+$ $h \circ d$ $=$ $\operatorname{id} - \operatorname{in} \circ \rho$ is defined by the formulas 

\begin{align*}
S_{--+}(\epsilon, \eta) &\mapsto S_{+-+}(\epsilon, -, \eta), \\
S_{+-+}(\epsilon, +, \eta) &\mapsto S_{+++}(\epsilon, \eta), \\
{\rm{otherwise}} &\mapsto 0.  
\end{align*}

(I\!I\!I--3) Consider the case where the state of $P'$ with (mark($A$), mark($B$), mark( $C$)) $=$ ($1$, $1$, $1$) is represented as $u|{\emptyset_{\epsilon}}^{w_{A}}|{\emptyset_{\eta}}^{w_{A B C}}|v$.  

%+++
$S_{+++}(\epsilon, \eta)$ denotes the state $u|{\emptyset_{\epsilon}}^{w_{A}}|{\emptyset_{\eta}}^{w_{A B C}}|v$ of $P'$ 
with (mark($A$), mark($B$), mark($C$)) $=$ ($1$, $1$, $1$), 
%-++
$S_{-++}(\epsilon)$ denotes the state $u|{\emptyset_{\epsilon}}^{w_{A B C}}|v$ of $P'$ 
with (mark($A$), mark($B$), mark($C$)) $=$ ($-1$, $1$, $1$), 
%+-+
$S_{+-+}(\epsilon, \zeta, \eta)$ denotes the state $u|{\emptyset_{\epsilon}}^{w_{A}}|{\emptyset_{\zeta}}^{A B C}|{\emptyset_{\eta}}^{w_{B C}}$ $|v$ of $P'$ 
with (mark($A$), mark($B$), mark($C$)) $=$ ($1$, $-1$, $1$), 
%--+
$S_{--+}(\epsilon, \eta)$ denotes the state $u|{\emptyset_{\epsilon}}^{w_{A B C}}|{\emptyset_{\eta}}^{w_{B C}}|v$ of $P'$ 
with (mark($A$), mark($B$), mark($C$)) $=$ ($-1$, $-1$, $1$), 
%++-
$S_{++-}(\epsilon, \zeta, \eta)$ denotes the state $u|{\emptyset_{\epsilon}}^{w_{A}}|{\emptyset_{\zeta}}^{w_{B C}}|{\emptyset_{\eta}}^{w_{A B C}}|v$ of $P'$ 
with (mark($A$), mark($B$), mark($C$)) $=$ ($1$, $1$, $-1$), 
%-+-
$S_{-+-}(\epsilon, \eta)$ denotes the state $u|{\emptyset_{\epsilon}}^{w_{A B C}}|{\emptyset_{\eta}}^{w_{B C}}|v$ of $P'$ 
with (mark($A$), mark($B$), mark($C$)) $=$ ($-1$, $1$, $-1$), 
%+--
$S_{+--}(\epsilon, \eta)$ denotes the state $u|{\emptyset_{\epsilon}}^{w_{A}}$ $|{\emptyset_{\eta}}^{w_{A B C}}|v$ of $P'$ 
with (mark($A$), mark($B$), mark($C$)) $=$ ($1$, $-1$, $-1$), 
%---
and $S_{---}(\epsilon)$ denotes the state $u|{\emptyset_{\epsilon}}^{w_{A B C}}|v$ of $P'$ 
with (mark($A$), mark($B$), mark($C$)) $=$ ($-1$, $-1$, $-1$).

The subcomplex $C'$ of ${C}(P')$ is defined by $C'$ $:=$ ${C}\big($$S_{-++}(+)$ $+$ $S_{+-+}(+, -, +)$, $S_{-++}(-)$ $+$ $S_{+-+}(+, -, -)$ $+$ $S_{+-+}(-, -, +)$, $S_{**-}$$\big)$, where $S_{**-}$ denotes every states with mark($C$) $=$ $-1$.  

%T+-+
$T_{+-+}(\epsilon)$ denotes the state $u|{\emptyset_{\epsilon}}^{w_{A B C}}|v$ of $P$ 
with (mark($A$), mark($B$), mark( $C$)) $=$ ($1$, $-1$, $1$), 
%T-++,-
$T_{-++}(\epsilon, -, \eta)$ denotes the state $u|{\emptyset_{\epsilon}}^{w_{A B C}}|{\emptyset_{-}}^{A B C}|{\emptyset_{\eta}}^{w_{B}}|v$ of $P$ 
with (mark($A$), mark($B$), mark($C$)) $=$ ($-1$, $1$, $1$), 
%++-
$T_{++-}(\epsilon, \zeta, \eta)$ denotes the state $u|{\emptyset_{\epsilon}}^{w_{A C}}$ $|{\emptyset_{\zeta}}^{w_{B}}|{\emptyset_{\eta}}^{w_{A B C}}|v$ of $P$ 
with (mark($A$), mark($B$), mark($C$)) $=$ ($1$, $1$, $-1$), 
%-+-
$T_{-+-}(\epsilon, \eta)$ denotes the state $u|{\emptyset_{\epsilon}}^{w_{A B C}}|{\emptyset_{\eta}}^{w_{B}}|v$ of $P$ 
with (mark($A$), mark($B$), mark(\\$C$)) $=$ ($-1$, $1$, $-1$), 
%+--
$T_{+--}(\epsilon, \eta)$ denotes the state $u|{\emptyset_{\epsilon}}^{w_{A C}}|{\emptyset_{\eta}}^{w_{A B C}}|v$ of $P$ 
with (mark( $A$), mark($B$), mark($C$)) $=$ ($1$, $-1$, $-1$), 
%---
$T_{---}(\epsilon)$ denotes the state $u|{\emptyset_{\epsilon}}^{w_{A B C}}|v$ of $P$ with (mark($A$), mark($B$), mark($C$)) $=$ ($-1$, $-1$, $-1$), 
%**-
and $T_{**-}$ denotes every state of $P$ with mark($C$) $=$ $-1$.  

The subcomplex $C$ of ${C}(P)$ is defined by $C$ $:=$ ${C}\big($$T_{+-+}(+)$ $+$ $T_{-++}(+, -, +)$, $T_{+-+}(-)$ $+$ $T_{-++}(+, -, -)$ $+$ $T_{-++}(-, -, +)$, $T_{**-}$$\big)$.

First, the retraction $\rho:$ ${C}(P')$ $\to$ $C'$ is defined by the formulas 
\begin{align*}
S_{-++}(+) &\mapsto S_{-++}(+) + S_{+-+}(+, -, +), \\
S_{-++}(-) &\mapsto S_{-++}(-) + S_{+-+}(+, -, -) + S_{+-+}(-, -, +), \\
S_{**-} &\mapsto S_{**-}, \\
S_{+-+}(+, +, +) &\mapsto S_{++-}(+, +, +), \\
S_{+-+}(+, +, -) &\mapsto S_{-++}(+) + S_{+-+}(+, -, +) + S_{++-}(+, +, -) + S_{++-}(+, -, +), \\
S_{+-+}(-, +, +) &\mapsto S_{-++}(+) + S_{+-+}(+, -, +) + S_{++-}(-, +, +), \\
S_{+-+}(-, +, -) &\mapsto S_{-++}(-) + S_{+-+}(+, -, -) + S_{+-+}(-, -, +) + S_{++-}(-, +, -) \\
&\qquad\qquad\,\,\,\, \qquad\qquad\qquad\qquad\qquad\qquad\qquad\qquad + S_{++-}(-, -, +), \\
S_{--+}(\epsilon, \eta) &\mapsto S_{+--}(\epsilon, \eta), \\
{\rm{otherwise}} &\mapsto 0.  
\end{align*}

Second, consider the following composition (\ref{third-rei-eq}) of the following isomorphism with $\rho$.  
The isomorphism $C'$ $\to$ ${C}$ is defined by the formulas
\begin{align*}
S_{-++}(+) + S_{+-+}(+, -, +) &\mapsto T_{+-+}(+) + T_{-++}(+, -, +), \\
S_{-++}(-) + S_{+-+}(+, -, -) + S_{+-+}(-, -, +) &\mapsto T_{+-+}(-) + T_{-++}(+, -, -) \\
&\qquad\qquad\quad\ \,  + T_{-++}(-, -, +), \\
S_{++-}(\epsilon) &\mapsto T_{++-}(\epsilon),  \\
S_{-+-}(\epsilon, \eta) &\mapsto T_{-+-}(\epsilon, \eta), \\
S_{+--}(\epsilon, \eta) &\mapsto T_{+--}(\epsilon, \eta), \\
S_{---}(\epsilon, \zeta, \eta) &\mapsto T_{---}(\epsilon, \zeta, \eta).  
\end{align*}

Third, the map $h:$ ${C}(P') \to {C}(P')$ such that $d \circ h$ $+$ $h \circ d$ $=$ $\operatorname{id} - \operatorname{in} \circ \rho$ is defined by the formulas 

\begin{align*}
S_{--+}(\epsilon, \eta) &\mapsto S_{+-+}(\epsilon, -, \eta), \\
S_{+-+}(\epsilon, +, \eta) &\mapsto S_{+++}(\epsilon, \eta), \\
{\rm{otherwise}} &\mapsto 0.  
\end{align*}

(I\!I\!I--4) Consider the case where the state of $P'$ with (mark($A$), mark($B$), mark( $C$)) $=$ ($1$, $1$, $1$) is represented as $u|{\emptyset_{\epsilon}}^{w_{A B C}}|{\emptyset_{\eta}}^{w_{C}}|v$.  

%+++
$S_{+++}(\epsilon, \eta)$ denotes the state $u|{\emptyset_{\epsilon}}^{w_{A B C}}|{\emptyset_{\eta}}^{w_{C}}|v$ of $P'$ 
with (mark($A$), mark($B$), mark($C$)) $=$ ($1$, $1$, $1$), 
%-++
$S_{-++}(\epsilon, \zeta, \eta)$ denotes the state $u|{\emptyset_{\epsilon}}^{w_{A B C}}|{\emptyset_{\zeta}}^{w_{A B}}|{\emptyset_{\eta}}^{w_{C}}|v$ of $P'$ 
with (mark($A$), mark($B$), mark($C$)) $=$ ($-1$, $1$, $1$), 
%+-+
$S_{+-+}(\epsilon, \zeta, \eta)$ denotes the state $u|{\emptyset_{\epsilon}}^{w_{A B}}|{\emptyset_{\zeta}}^{A B C}|{\emptyset_{\eta}}^{w_{C}}|v$ of $P'$ 
with (mark($A$), mark($B$), mark($C$)) $=$ ($1$, $-1$, $1$), 
%--+
$S_{--+}(\epsilon, \eta)$ denotes the state $u|{\emptyset_{\epsilon}}^{w_{A B}}|{\emptyset_{\eta}}^{w_{C}}|v$ of $P'$ 
with (mark($A$), mark($B$), mark($C$)) $=$ ($-1$, $-1$, $1$), 
%++-
$S_{++-}(\epsilon)$ denotes the state $u|{\emptyset_{\epsilon}}^{w_{A B C}}|v$ of $P'$ 
with (mark($A$\\), mark($B$), mark($C$)) $=$ ($1$, $1$, $-1$), 
%-+-
$S_{-+-}(\epsilon, \eta)$ denotes the state $u|{\emptyset_{\epsilon}}^{w_{A B C}}|{\emptyset_{\eta}}^{w_{A B}}|v$ of $P'$ 
with (mark($A$), mark($B$), mark($C$)) $=$ ($-1$, $1$, $-1$), 
%+--
$S_{+--}(\epsilon, \eta)$ denotes the state $u|{\emptyset_{\epsilon}}^{w_{A B}}|{\emptyset_{\eta}}^{w_{A B C}}|v$ of $P'$ 
with (mark($A$), mark($B$), mark($C$)) $=$ ($1$, $-1$, $-1$), 
%---
and $S_{---}(\epsilon)$ denotes the state $u|{\emptyset_{\epsilon}}^{w_{A B C}}|v$ of $P'$ 
with (mark($A$), mark($B$), mark($C$)) $=$ ($-1$, $-1$, $-1$).  

The subcomplex $C'$ of ${C}(P')$ is defined by $C'$ $:=$ ${C}\big($ $S_{-++}(+, +, \eta)$, $S_{-++}(+, -, $ $\eta)$ $+$ $S_{+-+}(+, -,$ $\eta)$, $S_{-++}(-, +, \eta)$ $+$ $S_{+-+}(+, -, \eta)$, $S_{-++}(-, -, \eta)$ $+$ $S_{+-+}(-, -,$ $\eta)$, $S_{**-}$ $\big)$, where $S_{**-}$ denotes every state with mark($C$) $=$ $-1$.  

%T+-+
$T_{+-+}(\epsilon, \zeta, \eta)$ denotes the state $u|{\emptyset_{\epsilon}}^{w_{A B C}}|{\emptyset_{\zeta}}^{w_{A B}}|{\emptyset_{\eta}}^{w_{C}}|v$ of $P$ 
with (mark($A$), mark( $B$), mark($C$)) $=$ ($1$, $-1$, $1$), 
%T-++,-
$T_{-++}(\epsilon, -, \eta)$ denotes the state $u|{\emptyset_{\epsilon}}^{w_{A B}}|{\emptyset_{-}}^{A B C}|{\emptyset_{\eta}}^{w_{C}}|v$ of $P$ 
with (mark($A$), mark($B$), mark($C$)) $=$ ($-1$, $1$, $1$), 
%++-
$T_{++-}(\epsilon)$ denotes the state $u|{\emptyset_{\epsilon}}^{w_{A B C}}|v$ of $P$ 
with (mark($A$), mark($B$), mark($C$)) $=$ ($1$, $1$, $-1$), 
%-+-
$T_{-+-}(\epsilon, \eta)$ denotes the state $u|{\emptyset_{\epsilon}}^{w_{A B}}|{\emptyset_{\eta}}^{w_{A B C}}|v$ of $P$ 
with (mark($A$), mark($B$), mark($C$)) $=$ ($-1$, $1$, $-1$), 
%+--
$T_{+--}(\epsilon, \eta)$ denotes the state $u|{\emptyset_{\epsilon}}^{w_{A B C}}|{\emptyset_{\eta}}^{w_{A B}}|v$ of $P$ 
with (mark($A$), mark($B$), mark($C$)) $=$ ($1$, $-1$, $-1$), 
%---
$T_{---}(\epsilon)$ denotes the state $u|{\emptyset_{\epsilon}}^{w_{A B C}}|v$ of $P$ 
with (mark($A$), mark($B$), mark($C$)) $=$ ($-1$, $-1$, $-1$), 
%**-
and $T_{**-}$ denotes every state of $P$ with mark($C$) $=$ $-1$.

The subcomplex $C$ of ${C}(P)$ is defined by $C$ $:=$ ${C}\big($$T_{+-+}(+, +, \eta)$, $T_{+-+}(+, -, \eta)$ $+$ $T_{-++}(+, -, \eta)$, $T_{+-+}(-, +, \eta)$ $+$ $T_{-++}(+, -, \eta)$, $T_{+-+}(-, -, \eta)$ $+$ $T_{-++}(-, -, \eta)$, $T_{**-}$$\big)$.

First, the retraction $\rho:$ ${C}(P')$ $\to$ $C'$ is defined by the formulas 
\begin{align*}
S_{-++}(+, +, \eta) &\mapsto S_{-++}(+, +, \eta), \\
S_{-++}(+, -, \eta) &\mapsto S_{-++}(+, -, \eta) + S_{+-+}(+, -, \eta), \\
S_{-++}(-, +, \eta) &\mapsto S_{-++}(-, +, \eta) + S_{+-+}(+, -, \eta), \\
S_{-++}(-, -, \eta) &\mapsto S_{-++}(-, -, \eta) + S_{+-+}(-, -, \eta), \\ 
S_{**-} &\mapsto S_{**-}, \\
\end{align*}
\begin{align*}
S_{+-+}(+, +, +) &\mapsto S_{-++}(+, +, +), \\
S_{+-+}(+, +, -) &\mapsto S_{-++}(+, +, -) + S_{++-}(+), \\
S_{+-+}(-, +, +) &\mapsto S_{-++}(+, -, +) + S_{-++}(-, +, +) + S_{++-}(+), \\
S_{+-+}(-, +, -) &\mapsto S_{-++}(+, -, -) + S_{-++}(-, +, -) + S_{++-}(-), \\
S_{--+}(\epsilon, \eta) &\mapsto S_{+--}(\epsilon, \eta), \\
{\rm{otherwise}} &\mapsto 0.  
\end{align*}

Second, consider the following composition (\ref{third-rei-eq}) of the following isomorphism with $\rho$.  
The isomorphism $C'$ $\to$ ${C}$ is defined by the formulas
\begin{align*}
S_{-++}(+, +, \eta) &\mapsto T_{+-+}(+, +, \eta), \\
S_{-++}(+, -, \eta) + S_{+-+}(+, -, \eta) &\mapsto T_{+-+}(+, -, \eta) + T_{-++}(+, -, \eta), \\
S_{-++}(-, +, \eta) + S_{+-+}(+, -, \eta) &\mapsto T_{+-+}(-, +, \eta) + T_{-++}(+, -, \eta), \\
S_{-++}(-, -, \eta) + S_{+-+}(-, -, \eta) &\mapsto T_{+-+}(-, -, \eta) + T_{-++}(-, -, \eta), \\
S_{++-}(\epsilon) &\mapsto T_{++-}(\epsilon),  \\
S_{-+-}(\epsilon, \eta) &\mapsto T_{-+-}(\epsilon, \eta), \\
S_{+--}(\epsilon, \eta) &\mapsto T_{+--}(\epsilon, \eta), \\
S_{---}(\epsilon) &\mapsto T_{---}(\epsilon).  
\end{align*}

Third, the map $h:$ ${C}(P') \to {C}(P')$ such that $d \circ h$ $+$ $h \circ d$ $=$ $\operatorname{id} - \operatorname{in} \circ \rho$ is defined by the formulas 

\begin{align*}
S_{--+}(\epsilon, \eta) &\mapsto S_{+-+}(\epsilon, -, \eta), \\
S_{+-+}(\epsilon, +, \eta) &\mapsto S_{+++}(\epsilon, \eta), \\
{\rm{otherwise}} &\mapsto 0.  
\end{align*}

(I\!I\!I--5) Consider the case where the state of $P'$ with (mark($A$), mark($B$), mark( $C$)) $=$ ($1$, $1$, $1$) is represented as $u|{\emptyset_{\epsilon}}^{w_{A}}|{\emptyset_{\zeta}}^{w_{A B C}}|{\emptyset_{\eta}}^{w_{C}}|v$.  

%+++
$S_{+++}(\epsilon, \zeta, \eta)$ denotes the state $u|{\emptyset_{\epsilon}}^{w_{A}}|{\emptyset_{\zeta}}^{w_{A B C}}|{\emptyset_{\eta}}^{w_{C}}|v$ of $P'$ 
with (mark($A$), mark( $B$), mark($C$)) $=$ ($1$, $1$, $1$), 
%-++
$S_{-++}(\epsilon, \eta)$ denotes the state $u|{\emptyset_{\epsilon}}^{w_{A B C}}|{\emptyset_{\epsilon}}^{w_{C}}|v$ of $P'$ 
with (mark($A$), mark($B$), mark($C$)) $=$ ($-1$, $1$, $1$), 
%+-+
$S_{+-+}(\epsilon, \zeta, \eta, \theta)$ denotes the state $u|{\emptyset_{\epsilon}}^{w_{A}}|{\emptyset_{\zeta}}^{w_{B}}|{\emptyset_{\eta}}^{w_{C}}|{\emptyset_{\theta}}^{A B C}|v$ of $P'$ 
with (mark($A$), mark($B$), mark($C$)) $=$ ($1$, $-1$, $1$), 
%--+
$S_{--+}(\epsilon, \zeta, \eta)$ denotes the state $u|{\emptyset_{\epsilon}}^{w_{A B C}}|{\emptyset_{\zeta}}^{w_{B}}|{\emptyset_{\eta}}^{w_{C}}|v$ of $P'$ 
with (mark($A$), mark($B$), mark($C$)) $=$ ($-1$, $-1$, $1$), 
%++-
$S_{++-}(\epsilon, \eta)$ denotes the state $u|{\emptyset_{\epsilon}}^{w_{A}}|{\emptyset_{\eta}}^{w_{A B C}}|v$ of $P'$ 
with (mark($A$), mark($B$), mark($C$)) $=$ ($1$, $1$, $-1$), 
%-+-
$S_{-+-}(\epsilon)$ denotes the state $u|{\emptyset_{\epsilon}}^{w_{A B C}}|v$ of $P'$ 
with (mark($A$), mark($B$), mark($C$)) $=$ ($-1$, $1$, $-1$), 
%+--
$S_{+--}(\epsilon, \zeta, \eta)$ denotes the state $u|{\emptyset_{\epsilon}}^{w_{A}}|{\emptyset_{\zeta}}^{w_{B}}|{\emptyset_{\eta}}^{w_{A B C}}|v$ of $P'$ 
with (mark($A$), mark($B$), mark($C$)) $=$ ($1$, $-1$, $-1$), 
%---
and $S_{---}(\epsilon, \eta)$ denotes the state $u|{\emptyset_{\epsilon}}^{w_{A B C}}|{\emptyset_{\eta}}^{w_{B}}|v$ of $P'$ 
with (mark($A$), mark($B$), mark($C$)) $=$ ($-1$, $-1$, $-1$).

The subcomplex $C'$ of ${C}(P')$ is defined by $C'$ $:=$ ${C}\big($ $S_{-++}(+, \eta)$ $+$ $S_{+-+}(+, +, \eta,$ $-)$, $S_{-++}(-, \eta)$ $+$ $S_{+-+}(+, -, \eta, -)$ $+$ $S_{+-+}(-, +, \eta, -)$, $S_{**-}$ $\big)$ where $S_{**-}$ denotes every states with mark($C$) $=$ $-1$.  

%T+-+
$T_{+-+}(\epsilon, \eta)$ denotes the state $u|{\emptyset_{\epsilon}}^{w_{A B C}}|{\emptyset_{\eta}}^{w_{A B}}|v$ of $P$ 
with (mark($A$), mark($B$), mark($C$)) $=$ ($1$, $-1$, $1$), 
%T-++,-
$T_{-++}(\epsilon, -)$ denotes the state $u|{\emptyset_{\epsilon}}^{w_{A B C}}|{\emptyset_{-}}^{A B C}|v$ of $P$ 
with (mark($A$), mark($B$), mark($C$)) $=$ ($-1$, $1$, $1$), 
%++-
$T_{++-}(\epsilon, \eta)$ denotes the state $u|{\emptyset_{\epsilon}}^{w_{A C}}$ $|{\emptyset_{\eta}}^{w_{A B C}}|v$ of $P$ 
with (mark($A$), mark($B$), mark($C$)) $=$ ($1$, $1$, $-1$), 
%-+-
$T_{-+-}(\epsilon)$ denotes the state $u|{\emptyset_{\epsilon}}^{w_{A B C}}|v$ of $P$ 
with (mark($A$), mark($B$), mark($C$)) $=$ ($-1$, $1$, $-1$), 
%+--
$T_{+--}(\epsilon, \zeta, \eta)$ denotes the state $u|{\emptyset_{\epsilon}}^{w_{A C}}|{\emptyset_{\zeta}}^{w_{B C}}|{\emptyset_{\eta}}^{w_{A B}}|v$ of $P$ 
with (mark($A$), mark($B$), mark($C$)) $=$ ($1$, $-1$, $-1$), 
%---
$T_{---}(\epsilon, \eta)$ denotes the state $u|{\emptyset_{\epsilon}}^{w_{A B C}}|{\emptyset_{\eta}}^{w_{B C}}|v$ of $P$ 
with (mark($A$), mark($B$), mark($C$)) $=$ ($-1$, $-1$, $-1$), 
%**-
$T_{**-}$ denotes every state of $P$ with mark($C$) $=$ $-1$.

The subcomplex $C$ of ${C}(P)$ is defined by $C$ $:=$ ${C}\big($$T_{+-+}(+, +)$, $T_{+-+}(+, -)$ $+$ $T_{-++}(+, -)$, $T_{+-+}(-, +)$ $+$ $T_{-++}(+, -)$, $T_{+-+}(-, -)$ $+$ $T_{-++}(-, -)$, $T_{**-}$$\big)$.  

First, the retraction $\rho:$ ${C}(P')$ $\to$ $C'$ is defined by the formulas 
\begin{align*}
S_{-++}(+, \eta) &\mapsto S_{-++}(+, \eta) + S_{+-+}(+, +, \eta, -), \\
S_{-++}(-, \eta) &\mapsto S_{-++}(+, -, \eta, -) + S_{+-+}(-, +, \eta, -), \\ 
S_{**-} &\mapsto S_{**-}, \\
S_{+-+}(+, +, -, +) &\mapsto S_{++-}(+, +), \\
S_{+-+}(+, -, +, +) &\mapsto S_{-++}(+, +) + S_{+-+}(+, +, +, -) + S_{++-}(+, +), \\
S_{+-+}(+, -, -, +) &\mapsto S_{-++}(+, -) + S_{+-+}(+, +, -, -) + S_{++-}(+, -), \\
S_{+-+}(-, +, +, +) &\mapsto S_{-++}(+, +) + S_{+-+}(+, +, +, -), \\
S_{+-+}(-, +, -, +) &\mapsto S_{-++}(+, -) + S_{+-+}(+, +, -, -) + S_{++-}(-, +), \\
S_{+-+}(-, -, +, +) &\mapsto S_{-++}(-, +) + S_{+-+}(+, -, +, -) + S_{+-+}(-, +, +, -) \\
&\qquad\qquad\qquad\qquad\qquad\qquad\qquad\quad\ \ + S_{++-}(-, +), \\
S_{+-+}(-, -, -, +) &\mapsto S_{-++}(-, -) + S_{+-+}(+, -, -, -) + S_{+-+}(-, +, -, -)\\
&\qquad\qquad\qquad\qquad\qquad\qquad\qquad\quad\ \  + S_{++-}(-, -), \\
S_{--+}(\epsilon, \zeta, \eta) &\mapsto S_{+--}(\epsilon, \zeta, \eta), \\
{\rm{otherwise}} &\mapsto 0.  
\end{align*}

Second, consider the following composition (\ref{third-rei-eq}) of the following isomorphism with $\rho$.  
The isomorphism $C'$ $\to$ ${C}$ is defined by the formulas
\begin{align*}
S_{-++}(+, +) + S_{+-+}(+, +, +, -) &\mapsto T_{+-+}(+, +), \\
S_{-++}(+, -) + S_{+-+}(+, +, -, -) &\mapsto T_{+-+}(+, -) \\&\qquad\!+ T_{-++}(+, -), \\
S_{-++}(-, +) + S_{+-+}(+, -, +, -) + S_{+-+}(-, +, +, -) &\mapsto T_{+-+}(-, +) \\&\qquad\!+ T_{-++}(+, -), \\
S_{-++}(-, -) + S_{+-+}(+, -, -, -) + S_{+-+}(-, +, -, -) &\mapsto T_{+-+}(-, -) \\&\qquad\!+ T_{-++}(-, -), \\
S_{++-}(\epsilon, \eta) &\mapsto T_{++-}(\epsilon, \eta), \\
\end{align*}
\begin{align*}
S_{-+-}(\epsilon) &\mapsto T_{-+-}(\epsilon), \\
S_{+--}(\epsilon, \zeta, \eta) &\mapsto T_{+--}(\epsilon, \zeta, \eta), \\
S_{---}(\epsilon, \eta) &\mapsto T_{---}(\epsilon, \eta).  
\end{align*}

Third, the map $h:$ ${C}(P') \to {C}(P')$ such that $d \circ h$ $+$ $h \circ d$ $=$ $\operatorname{id} - \operatorname{in} \circ \rho$ is defined by the formulas 

\begin{align*}
S_{--+}(\epsilon, \zeta, \eta) &\mapsto S_{+-+}(\epsilon, \zeta, \eta, -), \\
S_{+-+}(\epsilon, \zeta, \eta, +) &\mapsto S_{+++}(\epsilon, \zeta, \eta), \\
{\rm{otherwise}} &\mapsto 0.  
\end{align*}

(I\!I\!I--1) -- (I\!I\!I--5) prove that $KH^{i, j}(xAByACzBCt)$ $\simeq$ $KH^{i, j}(xBAyCAzCBt)$ if $(|A|, |B|, |C|)$ is any among $\{(-1, -1, -1)$, $(-1, 1, 1)$, $(1, 1, -1)$ $\}$.  

Consider $P'$ $=$ $xAByACzBCt$ $\rightarrow$ $xBAyCAzCBt$ $=$ $P$, where $(|A|, |B|, |C|)$ $=$ $(1, -1, -1)$.  
\begin{align*}
xAByACzBCt &\stackrel{\nu-\text{shift}}{\sim} xBCyABzACt \quad \text{with}\ (|A|, |B|, |C|) = (1, -1, -1) \\
           &\stackrel{\text{isom}}{\simeq} xAByDAzDBt \quad \text{with}\ (|A|, |B|, |D|) = (-1, -1, 1)\\
           &\stackrel{\text{Lemma \ref{abab}}}{\sim} xAByDA\underline{CE}zDB\underline{CE}t 
           \quad \text{with}\ |C| = -1, |E| = 1\\
           &\stackrel{H3}{\sim} x\underline{BA}yD\underline{CA}EzD\underline{CB}Et \quad \text{with}\ (|A|, |B|, |C|) = (-1, -1, -1)\\
           &\stackrel{H2}{\sim} x{BA}yAEzBEt \quad \text{with}\ |C| = -1, |D| = 1\\
           &\stackrel{\text{isom}}{\simeq} xCByBAzCAt \quad \text{with}\ (|A|, |B|, |C|) = (1, -1, -1)\\
           &\stackrel{\nu-\text{shift}}{\sim} xBAyCAzCBt \quad \text{with}\ (|A|, |B|, |C|) = (1, -1, -1) 
\end{align*}
We have already shown the invariance of $KH^{i,j}$ under the above moves and that $KH^{i,j}$ is preserved under the third homotopy move H3 and its inverse move with $(|A|, |B|, |C|)$ $=$ $(1, -1, -1)$.  
In particular, in this case, we use the invariance of $KH^{i,j}$ under H3 and its inverse with $($$|A|,$ $|B|,$ $|C|)$ $=$ $(-1,$ $-1,$ $-1)$.  By using the invariance under H3 and its inverse with $($$|A|,$ $|B|,$ $|C|)$ $=$ $(1,$ $1,$ $-1)$ (resp. $(-1,$ $1,$ $1)$), we can verify the invariance of $KH^{i,j}$ under H3 and its inverse with $($$|A|,$ $|B|,$ $|C|)$ $=$ $(1,$ $1,$ $1)$ (resp. $(-1,$ $-1,$ $1)$).  

We conclude that $KH^{i,j}(P')$ $\simeq$ $KH^{i,j}(P)$ for $P'$ $\simeq_{S_{1}}$ $P$.  
\end{proof}

The following corollary is a similar to Corollary \ref{jones-s0}.  
\begin{corollary}\label{cor-s0-homotopy}
$KH^{i, j}(P)$ are $S_{0}$-homotopy invariants for nanophrases $P$ over $\alpha_{0}$.  
\end{corollary}

\section{An application of $KH^{i, j}$ via words to nanophrases over any $\alpha$}\label{sec6}

In the previous sections, we discuss $S_{1}$-homotopy invariants $\hat{J}(P)$ and $KH^{i,j}(P)$ of pseudolinks.  Here, we construct homotopy invariants of nanophrases over any $\alpha$ from $\hat{J}(P)$ and $KH^{i,j}(P)$.

Let $\alpha$ be an arbitrary alphabet, $\tau$ be $\alpha \to \alpha$; involution, $\Delta_{\alpha}$ be $\{(a, a, a)\}_{a \in \alpha}$, and $\alpha/\tau$ $:=$ $\{\tilde{a}_{1}, \dots, \tilde{a}_{m}\}$.  We consider a complete residue system $\{a_{1}, \dots, a_{m}\}$ of $\alpha/\tau$ and denote $\{a_{1}, \dots, a_{m}\}$ by ${\rm crs}(\alpha/\tau)$.  

We use the notation of Definition \ref{orbit} as in \cite[Section 4.1]{turaev1}.  
\begin{definition}\label{orbit}
An {\it orbit} of the involution $\tau:$ $\alpha \to \alpha$ is a subset of $\alpha$ consisting either of one element preserved by $\tau$ or of two elements permuted by $\tau$; in the latter case, the orbit is {\it free}.  
\end{definition}

\begin{definition}
For $A \in \mathcal{A}$, we define the ${\rm sign}$ of $A$ by 
\begin{equation}
{\rm sign}_{L} (A) := \begin{cases}
                   & 1~\text{if}~|A| \in L ; |\tilde{A}|: \text{a free orbit}\\
                   & -1~\text{if}~|A| \in \tau(L) ; |\tilde{A}|: \text{a free orbit}\\
                   & 0~\ \ \text{otherwise}
                  \end{cases} 
\end{equation}
where $L$ is a nonempty subset of ${\rm crs}(\alpha/\tau)$. 
\end{definition}

Let $\mathcal{P}_{k}(\alpha, \tau)$ be a set of nanophrases of length $k$ over $\alpha$ with $\tau$.  

\begin{definition}
For an arbitrary $(\alpha, \tau)$ and an arbitrary subset $L \subset {\rm crs}(\alpha/\tau)$, 
$\mathcal{U}_{L}:$ $\mathcal{P}_{k}(\alpha, \tau) \to \mathcal{P}_{k}(\alpha_{0}, \tau_{0});$ $P \mapsto P_{0}$ is defined by the following two steps:

(Step 1) Remove $A \in \mathcal{A}$ such that ${\rm sign}_{L}(A) = 0$ from $(\mathcal{A}, P)$ $\in \mathcal{P}_{k}(\alpha, \tau)$.  

(Step 2) Let the nanophrase be $(\mathcal{A}', P')$ after removing letters from $(\mathcal{A}, P)$ by using (Step 1).  
We consider an $\alpha_{0}$-alphabet $\mathcal{B}$ such that ${\rm card} \mathcal{B}$ $=$ ${\rm card} \mathcal{A}'$ and $\mathcal{A}' \cap \mathcal{B}$ is the empty set.  Transpose each letter of $(\mathcal{A}', P')$ and a letter in $\mathcal{B}$ as follows: 
\begin{equation}
\begin{cases}
          &\text{transform}~$A$~\text{with}~{\rm sign}(A) = 1~\text{into}~B \in \mathcal{B}~\text{with}~|B| = 1\\
          &\text{transform}~$A$~\text{with}~{\rm sign}(A) = -1~\text{into}~B \in \mathcal{B}~\text{with}~|B| = -1.  
       \end{cases}
\end{equation}
By (1) and (2), the nanophrase over $\alpha_{0}$ derived from $(\mathcal{A}, P)$ is denoted by $\mathcal{U}_{L}((\mathcal{A}, P))$ or simply $\mathcal{U}_{L}(P)$.  
\end{definition}

\begin{theorem}\label{0_to_diagonal}
For an arbitrary $L \subset {\rm crs}(\alpha/\tau)$ and for two arbitrary nanophrases $(\mathcal{A}_{1}, P_{1})$ and $(\mathcal{A}_{2}, P_{2})$, 
\[(\mathcal{A}_{1}, P_{1}) \simeq_{\Delta_{\alpha}} (\mathcal{A}_{2}, P_{2}) \Longrightarrow \mathcal{U}_{L}((\mathcal{A}_{1}, P_{1})) \simeq_{S_{0}} \mathcal{U}_{L}((\mathcal{A}_{2}, P_{2})).  \]
\end{theorem}
\begin{proof}
It is clear that isomorphisms does not change the $\mathcal{U}_{L}(P)$.\par
Consider the first homotopy move 
$$P_{1}:=(\mathcal{A}, xAAy) \longrightarrow P_{2}:=(\mathcal{A} \setminus
\{A\}, xy)$$
where $x$ and $y$ are words on $\mathcal{A}$, 
possibly including the character  ``$|$".
Suppose ${\rm sign}(A) \neq 0$. Then,
$$\mathcal{U}_{L}(P_{1})=x_{L}AAy_{L} \simeq x_{L}y_{L} = 
\mathcal{U}_{L}(P_{2})$$
where $x_{L}$ and $y_{L}$ are words that are obtained by deleting all letters
$X \in \mathcal{A}$, such that ${\rm sign}(X) = 0$, from $x$ and
$y$, respectively.\par
Suppose ${\rm sign}(A) = 0$. Then,
$$\mathcal{U}_{L}(P_{1})=x_{L}y_{L}=\mathcal{U}_{L}(P_{2}).$$
Thus the first homotopy move does not change the 
homotopy class of $\mathcal{U}_{L}(P)$.\par
Consider the second homotopy move
$$P_{1}:=(\mathcal{A}, xAByBAz) \longrightarrow 
(\mathcal{A} \setminus \{A,B\}, xyz)$$
where $|A|=\tau(|B|)$, and $x$, $y$, and $z$ are words on $\mathcal{A}$
possibly including the character ``$|$".
Suppose $|A| \in L \cup \tau(L)$ and $\tilde{|A|}$ is free orbit. 
Then, $|B| \in L \cup \tau(L)$ and $\tilde{|A|}$ is free orbit 
since $|A|= \tau(|B|)$. Thus
$$\mathcal{U}_{L}(P_{1})=x_{L}ABy_{L}BAz_{L} 
\simeq x_{L}y_{L}z_{L} = \mathcal{U}_{L}(P_{2}).$$
where $x_{L}$, $y_{L}$ and $z_{L}$ are words 
that are obtained by deleting all letters
$X \in \mathcal{A}$, such that ${\rm sign}(X) = 0$, from $x$, $y$ and 
$z$, respectively.
Suppose $|A| \not\in L \cup \tau(L)$ or $|A|$ is a fixed point of $\tau$. 
Then, $|B| \not\in L \cup \tau(L)$ or $|B|$ is a fixed point of $\tau$ 
since $|A|= \tau(|B|)$. Thus,
$$\mathcal{U}_{L}(P_{1})= x_{L}y_{L}z_{L} = \mathcal{U}_{L}(P_{2}).$$
The above equation shows that 
the second homotopy move does not change the homotopy 
class of $\mathcal{U}_{L}(P)$.\par
Consider the third homotopy move 
$$P_{1}:= (\mathcal{A}, xAByACzBCt) \rightarrow 
P_{2}:=(\mathcal{A}, xBAyCAzCBt)$$
where $|A|=|B|=|C|$, and $x$, $y$, $z$, and $t$ are words on 
$\mathcal{A}$ possibly including the character ``$|$".
Suppose ${\rm sign}(A) \neq 0$. Then, ${\rm sign}(B), {\rm sign}(C) \neq 0$ 
since $|A|=|B|=|C|$. Thus we obtain
$$\mathcal{U}_{L}(P_{1}) = x_{L}ABy_{L}ACz_{L}ACt_{L} 
\simeq x_{L}BAy_{L}CAz_{L}CBt_{L}
=\mathcal{U}_{L}(P_{2}).$$
where $x_{L}$, $y_{L}$, $z_{L}$ and $t_{L}$ 
are words that are obtained by deleting all letters
$X \in \mathcal{A}$, such that ${\rm sign}(X) = 0$, from $x$, $y$, $z$, 
and $t$ respectively.
Suppose ${\rm sign}(A) = 0$. Then, ${\rm sign}(B), {\rm sign}(C) = 0$ 
since $|A|=|B|=|C|$. Thus we obtain
$$\mathcal{U}_{L}(P_{1}) = x_{L}y_{L}z_{L}t_{L} = \mathcal{U}_{L}(P_{2}).$$
Thus the third homotopy move does not change 
the homotopy class of $\mathcal{U}_{L}(P)$.\par
The above equation shows that 
$\mathcal{U}_{L}$ is a homotopy invariant of nanophrases.   
\end{proof}

\begin{corollary}\label{f-cor}
Let $I$ be an $S_{0}$-homotopy invariant of nanophrase over $\alpha_{0}$.  
For $P$ $\in \mathcal{P}_{k}(\alpha, \tau)$, we define $I'$ as \[I'(P) := \big\{I(\mathcal{U}_{L}(P))\big\}_{L \subset {\rm crs(\alpha/\tau)}}.  \]
$I'$ is a $\Delta_{\alpha}$-homotopy invariant of $P$ $\in \mathcal{P}_{k}(\alpha, \tau)$.  
In particular, for $(\mathcal{A}, P)$ $\in \mathcal{P}_{k}(\alpha_{0}, \tau_{0})$, $I'(P)$ $=$ $\{I(P)\}$ if ${\rm crs}(\alpha_{0}/\tau_{0})$ $=$ $\{1\}$.  
\end{corollary}

Theorem \ref{0_to_diagonal} implies the following corollaries.  

\begin{corollary}\label{jones-apply}
Let $\alpha$ be an arbitrary alphabet.  
$\hat{J}(\mathcal{U}_{L}(P))$ are $\Delta_{\alpha}$-homotopy invariants for nanophrases $P$ over $\alpha$.  
\end{corollary}
\begin{corollary}\label{homology-apply}
Let $\alpha$ be an arbitrary alphabet.  
$KH^{i,j}(\mathcal{U}_{L}(P))$ are $\Delta_{\alpha}$-homotopy invariants for nanophrases $P$ over $\alpha$.  
\end{corollary}

\begin{remark}
$\hat{J}(\mathcal{U}_{L}(P))$ $=$ $\sum_{j= - \infty}^{\infty}q^{j}\sum_{i= - \infty}^{\infty}(-1)^{i}{\rm rk}KH^{i,j}(\mathcal{U}_{L}(P))$.  
\end{remark}

We present some examples of the calculation of $KH^{i,j}(P)$ or $KH^{i,j}(\mathcal{U}_{L}(P))$.  
\begin{example}
For two pseudolinks $P_{1}$ $=$ $ABCDEABCDE$ with $|A|$ $=$ $|B|$ $=$ $|C|$ $=$ $|D|$ $=$ $|E|$ $=$ $-1$ and $P_{2}$ $=$ $ABCDEFBGDHFIJEHCGAIJ$ with $|A|$ $=$ $|C|$ $=$ $|E|$ $=$ $|G|$ $=$ $|H|$ $=$ $|I|$ $=$ $|J|$ $=$ $-1$ and $|B|$ $=$ $|D|$ $=$ $|D|$ $=$ $|F|$ $=$ $1$, $\hat{J}(P_{1})$ $=$ $\hat{J}(P_{2})$.  However, $KH^{-7, 15}(P_{1})$ $\simeq$ $0$ and $KH^{-7, 15}(P_{2})$ $\simeq$ $\mathbb{Z}_{2}$ (see \cite{bar-natan, turner}).
\end{example}
\begin{theorem}
$KH^{i,j}(P)$ is a strictly stronger invariant than $\hat{J}(P)$.  
\end{theorem}

\begin{corollary}
$KH^{i, j}(\mathcal{U}_L(P))$ is a strictly stronger invariant than $\hat{J}(P)$ for nano-phrases $P$ over $\alpha$.  
\end{corollary}

In \cite{turaev1}, Turaev constructed a $\Delta_{\alpha}$-homotopy invariant $\lambda$ for nanophrases over $\alpha$.  
\begin{example}\label{lambda}
Let $a$, $b$, and $c$ be elements (possibly coinciding) of any alphabet $\alpha$ and  
$A$, $B$, and $C$ be letters with $|A|$ $=$ $a$, $|B|$ $=$ $b$, and $|C|$ $=$ $c$.  
If $a$ $=$ $c$ $=$ $\tau(b)$ $\neq$ $b$, $\lambda(ABACBC)$ $=$ $\lambda(ACAC)$ $=$ $a$ $+$ $a_{\bullet}$ $-$ $a{a_{\bullet}}^{2}$ $-$ $a^{2}a_{\bullet}$.  
However, $KH^{0, 2}(\mathcal{U}_{\{a\}}(ACAC))$ $\simeq$ $0$ and $KH^{0, 2}(\mathcal{U}_{\{a\}}(ABACBC))$ $\simeq$ $\mathbb{Z}_{2}$.  
\end{example}
\begin{remark}
Turaev comments that the invariant $\lambda$ and all the other invariants of nanowords introduced form the beginning to Section 13.2 of \cite{turaev1} do not distinguish the two nanowords in Example \ref{lambda}.  
\end{remark}
Turaev constructed a strictly stronger $\Delta_{\alpha}$-homotopy invariant $f \circ v_{+}$ than $\lambda$ for nanophrases over $\alpha$ \cite{turaev1}.  
\begin{example}
Let $a$, $b$, $c$, and $d$ be elements (possibly coinciding) of any alphabet $\alpha$ and $A$, $B$, $C$, and $D$ be letters with $|A|$ $=$ $a$, $|B|$ $=$ $b$, $|C|$ $=$ $c$, and $|D|$ $=$ $d$.  
If $a$ $=$ $b$, $c$ $=$ $\tau(b)$ $=$ $d$, $a$ $\neq$ $\tau(b)$, and $c$ $\neq$ $\tau(d)$, $f(v_{+}(ABCDCDAB))$ $=$ $f(v_{+}(\emptyset))$ $=$ $\underline{1}$.  
However, $KH^{0, 3}(\mathcal{U}_{\{c\}}(\emptyset))$ $\simeq$ $0$ and $KH^{0, 3}(\mathcal{U}_{\{c\}}(ABCDCDAB))$ $\simeq$ $\mathbb{Z}_{2}$.  
\end{example}

\begin{theorem}
Let $\alpha$ be an arbitrary alphabet and $S$ be $\Delta_{\alpha}$.  
$KH^{i,j}(\mathcal{U}_{L}(P))$ is independent of $f \circ v_{+}$ for nanophrases $P$ over $\alpha$.  
\end{theorem}

\section*{Acknowledgments}
The authors would like to express their gratitude to Professors Toshiyuki Akita, Goo Ishikawa and Jun Murakami for their support.  The authors also wish to thank Andrew Gibson, Professor Kokoro Tanaka for their useful comments.  The authors are Research Fellows of the Japan Society for the Promotion of Science.  This work was partly supported by KAKENHI.

\vspace{0.3cm}
Department of Mathematics

Hokkaido University

Sapporo 060-0810, Japan

e-mail: fukunaga@math.sci.hokudai.ac.jp

\vspace{0.5cm}

Department of Pure and Applied Mathematics

Waseda University

Tokyo 169-8555, Japan

e-mail: noboru@moegi.waseda.jp
\end{document}